\documentclass[a4paper]{article}

\usepackage{graphicx}
\usepackage{latexsym}
\usepackage{amsmath}
\usepackage{amssymb}

\allowdisplaybreaks

\newenvironment{itemise}{\begin{itemize}}{\end{itemize}}

\newtheorem{proposition}{Proposition}

\newcommand{\vblock}{\rule[-4ex]{0ex}{9ex}}

\newcommand{\brak}[1]{\ensuremath{\left[#1\right]}}
\newcommand{\bbrak}[1]{\ensuremath{\Bigl[#1\Bigr]}}

\newcommand{\paren}[1]{\ensuremath{\left(#1\right)}}
\newcommand{\bparen}[1]{\ensuremath{\Bigl(#1\Bigr)}}

\newcommand{\mat}[1]{\begin{bmatrix}#1\end{bmatrix}}

\newcommand{\re}[1]{\mathrm{Re}\paren{#1}}
\newcommand{\im}[1]{\mathrm{Im}\paren{#1}}
\newcommand{\adj}{\mathop{\mathrm{adj}}}
\newcommand{\sgn}{\mathop{\mathrm{sgn}}}
\DeclareMathOperator{\arccot}{arccot}

\newcommand{\cplx}{\mathbb{C}}

\newcommand{\e}{\mathrm{e}}
\newcommand{\goesto}{\rightarrow}

\newcommand{\reals}{\mathbb{R}}

\newcommand{\vgnd}[1]{\ensuremath{v^{[#1]}}}
\newcommand{\ygnd}[2]{\ensuremath{Y^{[{#2},e_{#1}]}}}
\newcommand{\ygndc}[4]{\ensuremath{Y^{[{#2},e_{#1}][e_{#3},{#4}]}}}
\newcommand{\ycol}[2]{\ensuremath{Y^{[{#2},{#1}]}}}

\newcommand{\tz}[2]{\ensuremath{\mathcal{T}_{{#1},{#2}}}}
\newcommand{\nw}{\ensuremath{\mathcal{N}}}
\newcommand{\tr}{\ensuremath{\mathfrak{T}}}
\newcommand{\path}{\ensuremath{\mathcal{P}}}
\newcommand{\ivl}{\ensuremath{\mathcal{I}}}
\newcommand{\gfull}{\overline{G}}
\newcommand{\ifull}{\overline{I}}
\newcommand{\vfull}{\overline{V}}

\begin{document}

\title{The Admittance Matrix and Network Solutions}
\author{Anthony B.\ Morton}
\date{January 2017---additions January 2022}
\maketitle

\section{Introduction}

This note collects some results from classical network theory concerning properties of the network admittance matrix, and the relationship between electrical characteristics of the network and various mathematical properties of the admittance matrix.

An admittance matrix exists for any linear electrical network, whose branches (excluding sources) satisfy Ohmic relations.
That is, the current $i_{\alpha}$ on any branch $\alpha$ is presumed equal to $y_{\alpha} v_{\alpha} = y_{\alpha} (v_{\alpha+} - v_{\alpha-})$, where $y_{\alpha}$ is a scalar quantity (the \emph{admittance} of the branch) and $v_{\alpha+}$, $v_{\alpha-}$ are the voltages at the network nodes linked by the branch.
Each branch is assigned an arbitrary orientation in order to impart a sign convention to the current $i_{\alpha}$ (and thus the node-branch topology of the network is described by an underlying \emph{directed graph}).
For the purpose of determining the admittance matrix, all independent sources in the network are suppressed---voltage sources being replaced with short circuits and current sources with open circuits.

For a network with $n$ nodes, the admittance matrix $Y$ is a $n \times n$ matrix with rows and columns indexed by the network nodes (hence it is alternatively known as the \emph{node-admittance matrix}).
It is defined as follows:
\begin{itemise}
\item
Each diagonal element $Y_{kk}$ is equal to the sum of admittances of all branches incident with node $k$ (whether oriented toward or away from node $k$).
\item
Each element $Y_{jk}$ with $j \neq k$ is the negative sum of admittances of all branches whose endpoints are node $j$ and node $k$ (regardless of orientation), or zero if there are no such branches.
\end{itemise}
Mathematically, the admittance matrix is equivalent to a \emph{weighted graph Laplacian}.
The ordinary graph Laplacian is the matrix as just defined with all admittances equal to 1---so $Y_{kk}$ is simply the number of branches incident with node $k$ in the network graph and $Y_{jk}$ the (negative) number of branches linking nodes $j$ and $k$, where branch orientation is ignored.

Maintaining the graph-theoretic setting, the admittance matrix can alternatively be defined as $Y = G \Upsilon G^T$, where $G$ is the (directed) \emph{graph incidence matrix} and $\Upsilon$ a diagonal matrix of individual branch admittances.
For a network with $n$ nodes and $m$ branches, $G$ is a $n \times m$ matrix whose rows and columns are indexed by the nodes and branches of the network graph respectively, and where $G_{jk}$ is equal to 1 if branch $k$ is directed into node $j$, to $-1$ if branch $k$ is directed out of node $j$, and to 0 otherwise.
The $m \times m$ diagonal matrix $\Upsilon$ has $\Upsilon_{kk} = y_k$, the admittance of branch $k$, and all other elements zero.
(More generally, off-diagonal elements can be used to represent \emph{mutual coupling}, where current on one branch $j$ induces potential on another branch $k$.)
If all $y_k = 1$, one obtains the ordinary graph Laplacian as $Y = G G^T$.
Observe that although the branch orientations figure in the definition of $G$, the matrix $Y$ itself remains independent of these orientations.

Several properties of the admittance matrix are immediately evident from the definition:
\begin{itemise}
\item
$Y$ is symmetric.
\item
All row and column sums of $Y$ are equal to zero.
Thus $Y$ has at least one zero eigenvalue with $1_n$, an $n$-vector of all ones, as an eigenvector.
(It will turn out that for a connected network with all admittances real and positive, $Y$ has rank $n - 1$, so this completely characterises the nullspace of $Y$.)
\item
If the network admittances are real and positive, $Y$ is weakly diagonally dominant: that is, $|Y_{jj}| \geq |Y_{jk}|$ for all $j$ and $k$.
\end{itemise}
Much of classical network theory focusses on the case of resistive DC networks, where all the $y_k$ are positive real numbers.
This imparts to $Y$ a wider range of mathematical properties relative to the general case, where $y_k$ may take arbitrary complex values in AC networks, or even be differential operators describing an implicit dynamical system.

\section{Determinants and Cofactors}

`Classical' electrical network analysis, as developed in the 19th and early 20th century by Kirchhoff, Heaviside, Steinmetz, Rayleigh, Jeans \cite{j:mtoeam} and others including the `Trinity Four' \cite{bsst:doris}, draws heavily on the theory of determinants.
This theory has ancient antecedents, and long predates the formal theory of matrices: originally, a determinant was simply a property of a system of linear equations.
Leibniz, Cramer, Vandermonde, Laplace, Binet, Cauchy and Jacobi developed key results on determinants, which were eventually subsumed into Cayley and Sylvester's matrix theory from the mid-19th century.

Some basic elements of the theory of determinants are collected here for reference.
In what follows, $n$ denotes the number of rows or columns in square matrix $A$ with elements $\{a_{jk}\}$.
The set of integers from $1$ to $n$ is denoted $[1, n]$.

\begin{enumerate}
\item\label{det:transpose}
$\det A = \det A^T$.
Thus, any property of the determinant stated in terms of columns is equally valid for rows, and vice versa.
\item\label{det:multiply}
The determinant is multiplicative: $\det A B = \det A \det B$.
This implies in particular that $\det I = 1$, where $I$ is the $n \times n$ identity matrix, and that $\det A^{-1} = 1 / \det A$ when $A^{-1}$ exists.
\item\label{det:nlinear}
$\det A$ is an $n$-linear function of the columns of $A$.
This specifically implies that:
\begin{enumerate}
\item\label{det:scaling}
Scaling any one column by a constant factor scales the determinant by the same amount.
More generally, if all columns are scaled by arbitrary factors, say $c_k$ for column $k$, then $\det A$ is scaled by $\prod_{k=1}^n c_k$.
In particular, $\det c A = c^n \det A$.
\item\label{det:zerocol}
If any column is all zeros, then $\det A = 0$.
(This follows immediately from (\ref{det:scaling}): consider a scaling factor of zero.)
\item\label{det:addcol}
If $A$ and $B$ differ only in column $k$, and $A +_k B$ denotes the matrix formed from $A$ by adding column $k$ from $B$ to column $k$ of $A$, then $\det A + \det B = \det(A +_k B)$.
\end{enumerate}
\item\label{det:alternate}
$\det A$ is an \emph{alternating form} on the columns of $A$.
That is to say, if any two columns of $A$ are identical, or differ only by a scaling, or more generally if any column $k$ can be written as a linear combination of columns other than $k$, then $\det A = 0$.
(In fact, the determinant is the \emph{unique} $n$-linear alternating function on the columns of $n \times n$ matrices such that $\det I = 1$.)
\item\label{det:permute}
Permuting the columns of $A$ multiplies $\det A$ by the sign of the permutation.
(This is another general property of alternating forms.)
In particular:
\begin{enumerate}
\item\label{det:exchange}
Exchanging any two columns of $A$ multiplies $\det A$ by $-1$.
\item\label{det:cycle}
Rotating through $m$ columns of $A$ multiplies $\det A$ by $(-1)^{m-1}$.
\end{enumerate}
\item\label{det:combcol}
The value of the determinant is unchanged if any multiple of one column is added to another.
More generally, if to column $k$ one adds any linear combination of columns \emph{other} than $k$, the value of the determinant is unchanged.
This is a consequence of properties (\ref{det:addcol}) and (\ref{det:alternate}) above: by (\ref{det:addcol}) it amounts to perturbing $\det A$ by the determinant of a matrix whose columns are linearly dependent, which is zero by (\ref{det:alternate}).
\end{enumerate}

The determinants of square submatrices of $A$ also play an important part of the theory, and in the electrical network analysis derived from it.
If from a $n \times n$ square matrix $A$ one deletes some row $j$ and some column $k$, the $(n - 1) \times (n - 1)$ matrix that remains is the \emph{$(j,k)$-minor} $A_{(j,k)}$ of $A$.
This is known as a \emph{principal minor} if $j = k$.
Note that in a principal minor, the diagonal elements correspond to diagonal elements of $A$.

The \emph{first cofactor} associated with the $(j,k)$-minor of $A$ is $C_{j,k}(A) = (-1)^{j+k} \det A_{(j,k)}$.
(The explicit dependence on $A$ may be dropped depending on context.)
It is useful to observe that the additional sign $(-1)^{j+k}$ in the definition is also the sign of any permutation of rows and columns of $A$ that transforms the minor $A_{(j,k)}$ into a principal minor of the transformed matrix.

First cofactors are often simply called \emph{cofactors}, and figure in the well-known formula (due to Laplace) for recursive evaluation of determinants by `cofactor expansion':
\begin{equation}
\det A = \sum_{k=1}^n a_{jk} C_{j,k}(A) \qquad \text{($j \in [1, n]$ arbitrary).}
\label{eq:laplace}
\end{equation}
Lsplace's formula, stated in (\ref{eq:laplace}) for expansion along an arbitrary row of $A$, is equally valid for expansion along an arbitrary column.
The \emph{classical adjugate} (sometimes \emph{classical adjoint}) $\adj A$ of $A$ is the transpose of the matrix formed from all first cofactors---thus, $(\adj A)_{jk} = C_{k,j}(A)$.
It has the property that $(\adj A) A = A (\adj A) = (\det A) I$, allowing (in principle) the explicit construction of matrix inverses by means of determinants.
Notice also that if $A$ is symmetric it is readily shown that $C_{j,k}(A) = C_{k,j}(A)$, so the adjugate matrix (and hence $A^{-1}$) is also symmetric.

The concept of minors and cofactors extends further, by considering the deletion of multiple rows and columns.
Thus, the \emph{$(jp,kq)$-minor} $A_{(jp,kq)}$ of $A$ is the $(n - 2) \times (n - 2)$ matrix obtained by deleting the rows $j$ and $p$, and the columns $k$ and $q$, where $j \neq p$ and $k \neq q$.
This can also be defined recursively from the $(j,k)$-minor, by deleting from $A_{(j,k)}$ one additional row and column corresponding to row $p$ and column $q$ respectively in the original matrix $A$.
(Naturally, one must pay careful attention to the indexing: if $j < p$ then row $p$ of $A$ corresponds to row $p - 1$ of $A_{(j,k)}$, not row $p$.
Fortunately, the correct adjustments to the row and column indices can always be determined by inspection of the index values.)
$A_{(jp,kq)}$ is a principal minor if $j = k$ and $p = q$.

The \emph{second cofactor} of $A$ associated with row indices $j$ and $p$, and column indices $k$ and $q$ (where $j \neq p$ and $k \neq q$), is defined as
\begin{equation}
C_{jp,kq}(A) = (-1)^{\sigma(j, p) + \sigma(k, q)} \det A_{(jp,kq)}, \quad \text{where} \quad
\sigma(a, b) = \begin{cases} a + b - 1, & a < b \\ a + b, & \text{otherwise.} \end{cases}
\label{eq:cof2}
\end{equation}
Just as with first cofactors, the additional sign in the definition of the second cofactor can be characterised as that of a permutation of rows and columns of $A$---the permutation being any one where $A_{(jp,kq)}$ becomes a principal minor of the transformed matrix, and where the excluded rows and columns appear in the same order in the transformed matrix as they appear in the indices.
(So for example, whenever $j$ and $p$ are consecutive integers $p = j + 1$, there is no sign change associated with the row indices.)
This scheme generalises in a natural way to cofactors of any order.

By definition, $C_{jp,kq}$ is sensitive to the ordering of the row and column indices, and indeed is antisymmetric: $C_{pj,kq} = C_{jp,qk} = -C_{jp,kq}$.
It is also found that if $A$ is symmetric, then $C_{kq,jp} = C_{jp,kq}$: interchanging row and column indices (while preserving the order) leaves cofactors unchanged.
The following additional conventions are assumed: $C_{jp,kq} = 0$ when $j = p$ or $k = q$ (ensuring the definition applies to all index sets without restriction); and $C_{12,12} = 1$ for any $2 \times 2$ matrix.

Importantly, the definition of $C_{jp,kq}$ by (\ref{eq:cof2}) also ensures consistency with the alternative recursive definition, as the (recursively sign-modified) first cofactor of the minor $A_{(j,k)}$: one has
\begin{equation}
C_{jp,kq}(A) = (-1)^{j + k + p' + q'} \det A_{(j,k)(p',q')} = (-1)^{j + k} C_{p',q'}\!\paren{A_{(j,k)}},
\label{eq:cof2a}
\end{equation}
where $p' = \sigma(j,p) - j$, $q' = \sigma(k,q) - k$, are the modified indices for obtaining the second minor from the minor $A_{(j,k)}$.

So then, if one develops the cofactor expansion (\ref{eq:laplace}) to second order, by expanding the determinant in each cofactor along an arbitrary row $p \neq j$ of $A$, one obtains
\begin{equation}
\det A = \sum_{k=1}^n a_{jk} (-1)^{j + k}
      \mathop{\sum_{q =1}^n}_{q \neq k} a_{pq} C_{p',q'}\!\paren{A_{(j,k)}}
   = \sum_{k=1}^n \mathop{\sum_{q=1}^n}_{q \neq k} a_{jk} a_{pq} C_{jp,kq}(A),
\label{eq:laplace2}
\end{equation}
where the formula is valid for any distinct row indices $j, p \in [1, n]$.
Note again that while the elements $a_{jk}$ and $a_{pq}$ are identified by the original row and column indices in $A$, the cofactor of $A_{(j,k)}$ must be identified by the modified indices $p'$, $q'$ in that matrix.
The overall result is a double sum over all second cofactors of the form $C_{jp,**}$ of the original matrix $A$.

Observe that this double sum is over all \emph{ordered} pairs $(k, q)$ with $q \neq k$.
If one considers some arbitrary pair $(a, b)$ of column indices without regard to order, it is seen from (\ref{eq:laplace2}) that a cofactor involving both these column indices appears exactly twice: as the cofactor $C_{jp,ab}$ with weighting $a_{ja} a_{pb}$, and as the cofactor $C_{jp,ba} = - C_{jp,ab}$ with weighting $a_{jb} a_{pa}$.
This leads to the equivalent formula
\begin{equation}
\det A = \sum_{k=1}^{n-1} \sum_{q=k+1}^n \left|A_{jp,kq}\right| C_{jp,kq}(A), \qquad
   j, p \in [1, n], \quad j \neq p \text{ arbitrary},
\label{eq:laplace2a}
\end{equation}
where the sum is over the $n (n - 1) / 2$ \emph{unordered} distinct pairs from $[1, n]$ (identified as $(k,q)$ with $k < q$ without loss of generality), and where $|A_{jp,kq}| = a_{jk} a_{pq} - a_{jq} a_{pk}$ is the $2 \times 2$ determinant obtained from the intersection of the respective columns with rows $j$ and $p$.

One may similarly proceed to express $\det A$ as a multiple sum of weighted cofactors of higher orders up to $n$.
At order $n - 1$ they become single elements of $A$ (up to sign), and the `$n$th cofactors' are all just $+1$ or $-1$, depending whether the permutation on the $n$ row and $n$ column indices is even or odd.
By carrying the recursive process through to $n$th order one recovers one of the oldest general formulae for the determinant, due to Leibniz:
\begin{equation}
\det A = \sum_{\sigma \in S_n} \sgn\sigma \prod_{k=1}^n a_{k \sigma(k)},
\label{eq:leibniz}
\end{equation}
where (in modern terminology) $S_n$ represents the group of all permutations of $[1, n]$, and $\sgn \sigma$ denotes the sign of the permutation $\sigma$ ($-1$ if it is obtained from the identity via an odd number of transpositions, and $+1$ if by an even number).
Looked at another way, (\ref{eq:leibniz}) expresses $\det A$ as a signed sum of $n$-fold products of elements of $A$, where the $n$ elements making up each product occupy distinct rows and columns in every possible way.
Coming full circle, it is possible to take (\ref{eq:leibniz}) as the \emph{definition} of the determinant, and derive all the above properties from this.

\section{Cofactors in Network Theory}

The significance of the admittance matrix and its cofactors for electrical network analysis stems from the fact that the matrix $Y$ expresses the \emph{node-voltage equations} for the network.
For each $k$, let $v_k$ denote the voltage (potential) at node $k$, and $i_k$ an external current injected into node $k$ from outside the network (with a negative sign if the current flow is outward).
Assume for now that no other sources are present.
Then the $n$-vector $v$ of all node voltages is related to the $n$-vector $i$ of all current injections by the matrix-vector equation
\begin{equation}
Y v = i.
\label{eq:nodev}
\end{equation}
Equation (\ref{eq:nodev}) follows in a straightforward way from Ohm's and Kirchhoff's laws, with row $k$ of (\ref{eq:nodev}) expressing Kirchhoff's current law at node $k$ by virtue of the construction of $Y$.

Of course as $Y$ is singular, equation (\ref{eq:nodev}) will not be solvable for $v$ for arbitrary current sources $i$.
In particular, premultiplying both sides by the row vector $1_n^T$ results in $0 = 1_n^T i$, implying that any combination of current injections satisfying (\ref{eq:nodev}) must sum to zero.
Intuitively, this expresses the fact that the net flow of charge across the boundary of the circuit must be zero.

So then, the best one can hope for getting meaningful solutions of (\ref{eq:nodev}) is that $Y$ is not `very' singular---specifically, that its rank is $n - 1$, the largest possible given its construction.
In that case, one has the following:
\begin{proposition}
\label{prop:unique}
Let $Y$ be of maximal rank $n - 1$.
Then for any vector $i$ of current injections satisfying $1_n^T i = 0$, there is a unique voltage vector $\hat{v}$ with $1_n^T \hat{v} = 0$ that satisfies (\ref{eq:nodev}).
Further, all solutions $v$ to (\ref{eq:nodev}) with the same $i$ can be expressed in the form $v = \hat{v} + \lambda 1_n$ with $\lambda$ a scalar.
\end{proposition}
The last part of this proposition expresses mathematically the common understanding of voltage as a relative concept.
One may be precise about the \emph{difference} in potential between two points, but there is no absolute zero of potential: one may add the same constant amount to every voltage in an electrical configuration without consequence.
Here, it reduces to a basic principle in linear algebra: $v$ and $v'$ are two solutions to (\ref{eq:nodev}) with the same $i$ if and only if $v - v'$ belongs to the nullspace of $Y$; and if $Y$ is of rank $n - 1$, its nullspace is precisely the linear variety $\lambda 1_n$ with $\lambda$ a free scalar.

The remainder of Proposition \ref{prop:unique} follows from a fundamental fact concerning the first cofactors of $Y$, which will provide the second key result here:
\begin{proposition}
\label{prop:cof1}
Let $Y$ be any square matrix having all row and column sums equal to zero.
Then all first cofactors $\{C_{jk}(Y)\}$ are equal to one another.
\end{proposition}
The proof is by construction, using elementary row and column operations and the properties of determinants to show that any one minor, say $Y_{(1,1)}$, may be transformed into an arbitrary minor $Y_{(j,k)}$ without altering its determinant, except in a defined manner with regard to sign.

Begin then with $Y_{(1,1)}$, which is the submatrix formed from rows 2 through $n$, and columns 2 through $n$ of $Y$.
The first step is to transform this to $Y_{(j,1)}$, comprising all rows of $Y$ except row $j$.
If $j = 1$, there is nothing to do.
Otherwise, consider what happens if to row $j - 1$ of $Y_{(1,1)}$ (corresponding to row $j$ of $Y$) one applies the following operations:
\begin{enumerate}
\item
Change the sign of all elements in row $j - 1$.
(This changes the sign of the determinant, by property \ref{det:scaling} applied to rows.)
\item
Subtract every \emph{other} row in $Y_{(1,1)}$ in turn from row $j - 1$.
(This leaves the determinant unchanged, by property \ref{det:combcol} applied to rows.)
\end{enumerate}
This has the overall effect of replacing row $j - 1$ with the negative sum of all rows of $Y_{(1,1)}$---that is, with row 1 of $Y$ (except the first element), since the column sums of $Y$ are all zero.
This gives a matrix that has the same rows and columns as $Y_{(j,1)}$, but not quite in the right order if $j > 2$, since row 1 of $Y$ is now in row position $j - 1$ of this matrix.
But this is easily remedied by applying a cyclic permutation to the first $j - 1$ rows of the matrix, bringing row $j - 1$ to the top and moving the other rows down.
So now the matrix has rows 1 through $j - 1$ matching the equivalent rows of $Y$ (minus the first elements) and rows $j$ through $n - 1$ likewise matching rows $j + 1$ through $n$ of $Y$---that is to say, the matrix is identical to $Y_{(j,1)}$.
Since permuting the rows will have applied $j - 2$ sign changes to the determinant, and there was one further sign change involved in replacing row $j - 1$, it has been shown that
\begin{equation}
\det Y_{(j,1)} = (-1)^{j - 1} \det Y_{(1,1)}.
\label{eq:detj1}
\end{equation}
It is easily verified from the above that the rule (\ref{eq:detj1}) applies also in the special cases $j = 1$ and $j = 2$.
By now applying the equivalent procedure to the \emph{columns} of $Y_{(j,1)}$ when $k > 1$---that is, replacing column $k - 1$ of $Y_{(j,1)}$ with column $1$ of $Y$ using column operations, then permuting the first $k - 1$ columns cyclically when $k > 2$---one shows likewise that
\begin{equation}
\det Y_{(j,k)} = (-1)^{k - 1} \det Y_{(j,1)} = (-1)^{j + k - 2} \det Y_{(1,1)} = (-1)^{j + k} \det Y_{(1,1)}.
\label{eq:detjk}
\end{equation}
Notice that all steps leading up to the rule (\ref{eq:detjk}) assumed only that $Y$ is a square matrix with row and column sums equal to zero.
Since by definition of cofactors one has $C_{1,1}(Y) = \det Y_{(1,1)}$, (\ref{eq:detjk}) shows immediately that $C_{j,k}(Y) = C_{1,1}(Y)$, and since $j$ and $k$ are arbitrary, this proves Proposition \ref{prop:cof1}.

Using Proposition \ref{prop:cof1}, one can now prove the remainder of Proposition \ref{prop:unique}: that there is a unique $\hat{v}$ with $1_n^T \hat{v} = 0$ satisfying (\ref{eq:nodev}), as long as $1_n^T i = 0$ and $Y$ is of rank $n - 1$.
First, observe that if $Y$ is of this rank, then its first cofactors cannot be zero.
Thus by Proposition \ref{prop:cof1}, all first cofactors of $Y$ are equal to some constant $c \neq 0$.
Now consider the linear system obtained from (\ref{eq:nodev}) by replacing the redundant $n$th equation with the alternative equation $1_n^T v = 0$.
This is equivalent to replacing the $n$th row of $Y$ with a row of all ones, obtaining a new matrix $Y'$, but leaves the cofactors $C_{n,k}(Y')$ unchanged and equal to $c$.
It follows by cofactor expansion along the $n$th row that $\det Y' = nc \neq 0$, and so the modified equations have a unique solution $\hat{v}$ with the required properties.

Notice now that in Proposition \ref{prop:unique}, one may replace the condition $1_n^T \hat{v} = 0$ on the unique solution with any other condition of the form $a^T \hat{v} = b$, where $a$ is an $n$-vector and $b$ a scalar, the only requirement being that $a^T 1_n \neq 0$.
For then one may replace row $n$ (or any other single row) of $Y$ with the vector $a^T$ and the equivalent element of $i$ with $b$, and (\ref{eq:nodev}) will still yield a unique solution, since the determinant of the modified $Y$ matrix is equal to $(a^T 1_n) c \neq 0$, where $c$ is the common value of all first cofactors of $Y$.

In the special case where $a = e_k$, the $k$th elementary $n$-vector, and $b = E$ is an arbitrary voltage, one obtains the following:
\begin{proposition}
\label{prop:voltsrc}
Let $Y$ be of maximal rank $n - 1$.
Then for any vector $i$ of current injections satisfying $1_n^T i = 0$, and any arbitrary assignment of voltage $E$ to node $k$ of the network (where $k \in [1, n]$), there is a unique voltage vector $v$ satisfying (\ref{eq:nodev}) and whose $k$th element is $v_k = E$.
This solution $v$ is related to $\hat{v}$ in Proposition \ref{prop:unique} by the elementwise formula
\begin{equation}
v_j = E + \hat{v}_j - \hat{v}_k, \qquad j \in [1, n].
\label{eq:voltsrcsoln}
\end{equation}
\end{proposition}
To see why the formula (\ref{eq:voltsrcsoln}) holds, observe that by the latter part of Proposition \ref{prop:unique}, the solution $v$ under Proposition \ref{prop:voltsrc} will also be of the form $\hat{v} + \lambda 1_n$ for some value of $\lambda$.
If $\hat{v}_k$ denotes the $k$th element of $\hat{v}$, it is immediate that $\lambda = E - \hat{v}_k$, and (\ref{eq:voltsrcsoln}) follows.

Proposition \ref{prop:voltsrc} has an immediate corollary that proves useful in developing practical network solutions.
Given any distinguished circuit node $k \in [1,n]$, let $\vgnd{k}$ denote the unique solution $v$ to (\ref{eq:nodev}) satisfying $v_k = 0$.
This solution $\vgnd{k}$ is said to be \emph{grounded at node $k$}.
Now consider any other solution $v$ with the same current vector $i$: let $E$ denote the voltage $v_k$ in this solution.
Again by Proposition \ref{prop:unique}, $v$ can be expressed in the form $\vgnd{k} + \lambda 1_n$ for some constant $\lambda$, and a quick comparison of the $k$th elements of $v$ and $\vgnd{k}$ reveals that in fact $\lambda = E$.
This of course implies that any other element $v_j$ of the solution $v$ is obtained in the same way from the corresponding element of $\vgnd{k}$, and thus
\begin{equation}
v_j = \vgnd{k}_j + E = \vgnd{k}_j + v_k.
\label{eq:vgnddiff}
\end{equation}
This leads to the following important property of grounded solutions:
\begin{proposition}
\label{prop:grounded}
Let $Y$ be of maximal rank $n - 1$ and let $k \in [1,n]$ be an arbitrary distinguished node in the circuit.
For a given vector $i$ of current injections satisfying $1_n^T i = 0$, let $\vgnd{k}$ denote the solution grounded at node $k$, that is, the unique voltage solution satisfying $\vgnd{k}_k = 0$.
Then an arbitrary voltage element $\vgnd{k}_j$ of this solution gives the voltage difference $v_j - v_k$ in \emph{all} other network voltage solutions having the same current injection $i$.
\end{proposition}
Proposition \ref{prop:grounded} justifies restricting one's attention to grounded solutions, since any other solution with the same input currents may be derived from this simply by applying the new voltage $E$ to node $k$ and treating the other voltages in the solution as the voltage \emph{differences} to node $k$.

When it comes to the actual calculation of network solutions, it is the second cofactors of $Y$ that play a key role.

Consider the simplest case for the injected currents $i$, with just two nonzero entries: one equal to $+I$ (at node $p$, say) and the other equal to $-I$ (at node $q \neq p$).
Note that this is equivalent to connecting an independent current source, with value $I$, directly between nodes $p$ and $q$, and directed toward node $p$.
(Such sources function as elementary generators of more general current injections $i$, in the sense that any vector of injected currents with $1_n^T i = 0$ can be realised as a sum of single sources connected between pairs of nodes; this remains true even if one requires that all sources have their negative terminal at a fixed node $k$.)

When an independent source $I$ is connected as above between nodes $p$ and $q$, the current injection vector is $i = I (e_p - e_q)$, where $e_k$ is the $k$th elementary $n$-vector.
The network solution for $v$ can therefore be found by considering the solution to simplified node-voltage equations (\ref{eq:nodev}), of the form $Y v = e_p$ where the right hand side has a single nonzero element.
By virtue of Proposition \ref{prop:grounded} and by linearity, one may further restrict to the case of grounded solutions having $v_k = 0$ for some arbitrary node $k$.
Imposing the condition $v_k = 0$ at once makes the equations nonsingular and simplifies the left hand side.
More precisely, one has
\begin{equation}
\ygnd{k}{r} \vgnd{k} = e_p,
\label{eq:nodevelem}
\end{equation}
where $\ygnd{k}{r}$ represents the matrix obtained by choosing some $r \in [1,n]$ distinct from $p$, and replacing row $r$ of $Y$ with the elementary row vector $e_k^T$.
The condition $r \neq p$ is needed for the time being, since if $r = p$ the right hand side of the equation $e_k^T \vgnd{k} = 0$ for the grounded solution replaces the nonzero element in $e_p$, resulting in the trivial solution $\vgnd{k} = 0$.

\begin{proposition}
\label{prop:gndsoln}
Let $Y$ be any $n \times n$ matrix and $r, k \in [1,n]$ such that $C_{r,k}(Y) \neq 0$.
Then for any $p \in [1,n]$ distinct from $r$, equation (\ref{eq:nodevelem}) has a unique solution $\vgnd{k}$, the elements of which are given by cofactors of $Y$ as follows:
\begin{equation}
\vgnd{k}_j = \frac{C_{pr,jk}(Y)}{C_{r,k}(Y)}.
\label{eq:gndsoln}
\end{equation}
\end{proposition}
The solution (\ref{eq:gndsoln}) is found using Cramer's rule for linear equations.
Applied to (\ref{eq:nodevelem}), this asserts the $j$th element of the solution is
\begin{equation}
\vgnd{k}_j = \frac{\det \ygndc{k}{r}{p}{j}}{\det \ygnd{k}{r}},
\label{eq:cramer}
\end{equation}
where the matrix $\ygndc{k}{r}{p}{j}$ is obtained from $\ygnd{k}{r}$, by replacing column $j$ with the vector $e_p$.

The denominator of (\ref{eq:cramer}) is immediately seen to be equal to $C_{r,k}(Y)$, by expanding along the replaced row $r$.
The numerator of (\ref{eq:cramer}) is evaluated through two cofactor expansions, first down column $j$ and then along row $r$.
The exception is the trival case $j = k$: in this case (\ref{eq:cramer}) evaluates to zero, since with $p \neq r$ a zero element of $e_p$ replaces the single nonzero element in row $r$.

Assuming then that $j \neq k$, expanding $\det \ygndc{k}{r}{p}{j}$ down column $j$ results in
\begin{equation}
\det \ygndc{k}{r}{p}{j} = C_{p,j}\!\paren{\ygndc{k}{r}{p}{j}} = (-1)^{p + j} \det \ygndc{k}{r}{p}{j}_{(p,j)}.
\label{eq:ygndc1}
\end{equation}
Observe now that the $(p,j)$-minor on the right hand side of (\ref{eq:ygndc1}) is nearly identical to the minor $Y_{(p,j)}$ of the original matrix $Y$: the only difference is that one row---with index $r' = \sigma(p,r) - p$ where $\sigma(p,r)$ was defined in (\ref{eq:cof2})---is replaced with an elementary row vector.
This row vector has a 1 in position $k' = \sigma(j,k) - j$, so a cofactor expansion along this row gives
\begin{equation}
\det \ygndc{k}{r}{p}{j} = (-1)^{p + j} (-1)^{r' + k'} \det \ygndc{k}{r}{p}{j}_{(pr,jk)}.
\label{eq:ygndc2}
\end{equation}
But now, since the $(pr,jk)$-minor of the modified matrix $\ygndc{k}{r}{p}{j}$ is in fact identical to the $(pr,jk)$-minor of the original matrix $Y$, the expression (\ref{eq:ygndc2}) reduces to a second cofactor of $Y$:
\begin{equation}
\det \ygndc{k}{r}{p}{j} = (-1)^{p + j + r' + k'} \det Y_{(pr,jk)} = C_{pr,jk}(Y).
\label{eq:ygndc3}
\end{equation}
This (along with the convention $C_{pr,kk} = 0$) establishes the formula (\ref{eq:gndsoln}).

(One may note that the above calculation of $\det \ygndc{k}{r}{p}{j}$ could just as well have been carried out by expanding first along the row $r$, and then down column $j$.
Carrying this through gives, in place of (\ref{eq:ygndc3}), the cofactor $C_{rp,kj}(Y)$.
However, the antisymmetry of the second cofactors ensures that this is in fact equal to $C_{pr,jk}(Y)$ as required.)

With the aid of Proposition \ref{prop:gndsoln}, one may now write the general network solution with a single current source $I$ connected between the nodes $p$ and $q$.
The equation for the grounded solution in this case is
\begin{equation}
\ygnd{k}{r} \vgnd{k} = I \paren{e_p - e_q},
\label{eq:nodevipq}
\end{equation}
where the row index $r$ is (for the time being) chosen to be distinct from both $p$ and $q$.
By Proposition \ref{prop:gndsoln} and linearity the solution for the voltage at any node $j$ is
\begin{equation}
\vgnd{k}_j = I \cdot \frac{C_{pr,jk}(Y) - C_{qr,jk}(Y)}{C_{r,k}(Y)}
   = I \cdot \frac{C_{pr,jk}(Y) - C_{qr,jk}(Y)}{c(Y)},
\label{eq:isolnpq}
\end{equation}
where $c(Y)$ is the common value of all first cofactors of $Y$.
(This notation emphasises that $c(Y)$ is a property intrinsic to the matrix $Y$, and to the underlying network.)

This formula was obtained with $r$ distinct from $p$ and $q$.
But if one lets $r = q \neq p$, the equation thus removed is still linearly dependent on the others in (\ref{eq:nodevipq}), so the solution is unaffected.
But now, the zero right hand side of the equation $e_k^T \vgnd{k} = 0$ in row $q$ replaces the 1 in the vector $e_q$ in (\ref{eq:nodevipq}), and the equation that results is
\begin{equation}
\ygnd{k}{q} \vgnd{k} = I \cdot e_p.
\label{eq:nodevi}
\end{equation}
Since $p \neq q$, Proposition \ref{prop:gndsoln} again applies and the solution is
\begin{equation}
\vgnd{k}_j = I \cdot \frac{C_{pq,jk}(Y)}{c(Y)}.
\label{eq:isoln}
\end{equation}
Comparing (\ref{eq:isoln}) with (\ref{eq:isolnpq}) reveals the following transitivity property of the second cofactors:
\begin{equation}
C_{pq,jk}(Y) = C_{pr,jk}(Y) - C_{qr,jk}(Y) = C_{pr,jk}(Y) + C_{rq,jk}(Y).
\label{eq:transitive}
\end{equation}
This will be revisited shortly.

Taking into account Proposition \ref{prop:grounded} regarding grounded solutions, (\ref{eq:isoln}) leads to the following general result:
\begin{proposition}
\label{prop:transz}
Let $Y$ be of maximal rank $n - 1$ and let $c(Y)$ denote the common value of all first cofactors of $Y$.
Then for any node indices $p,q,j,k \in [1,n]$, the ratio
\begin{equation}
\tz{pq}{jk} = \frac{C_{pq,jk}(Y)}{c(Y)}
\label{eq:transz}
\end{equation}
when multiplied by the source current $I$ of an independent current source directed from node $q$ to node $p$, gives the voltage difference $v_j - v_k$ between the nodes $j$ and $k$ produced in the network due to that source.
$\tz{pq}{jk}$ is denoted the \emph{transfer impedance} or \emph{transpedance} between the node pair $(p,q)$ and the node pair $(j,k)$.
\end{proposition}

A number of well-known properties of linear electrical networks follow as straightforward consequences of Proposition \ref{prop:transz}.
\begin{description}
\item[Reciprocity.]
Since $Y$ is symmetric, second cofactors with transposed row and column indices are equal: $C_{pq,jk}(Y) = C_{jk,pq}(Y)$.
This means that the same transfer impedance $\tz{pq}{jk}$ is also the ratio of the voltage $v_p - v_q$ to the current $I$ of a source directed from node $k$ to node $j$.
\item[Transitivity.]
The transitivity property (\ref{eq:transitive}), put in the language of Proposition \ref{prop:transz}, states that placing two current sources of the same value $I$ between nodes $q$ and $r$, and $r$ and $p$ respectively, is equivalent to placing a single source between nodes $q$ and $p$.
Under reciprocity, it expresses the addition principle for voltages in the network solution: $v_p - v_q = (v_p - v_r) + (v_r - v_q)$.
\item[Superposition.]
The superposition property of solutions is a consequence of the linearity of the original node-voltage equations.
As mentioned above, any admissible current injection vector $i$ can be exprssed as a sum of influences of individual sources between specific pairs of nodes.
Proposition \ref{prop:transz} provides a method of calculating the influence of each individual source, and these influences are summed to provide the overall solution.
\item[Symmetry.]
This is another consequence of linearity, but is made explicit in the relation of transfer impedances to cofactors.
The symmetry relation $C_{qp,kj}(Y) = C_{pq,jk}(Y)$, thus $\tz{qp}{kj} = \tz{pq}{jk}$, is equivalent to stating that reversing the polarity of the current source driving the network reverses the signs of all consequent voltage differences.
\item[Driving-point impedances.]
The transfer impedance from a node pair $(j,k)$ to itself is the \emph{driving-point impedance}, \emph{self-impedance} or \emph{effective impedance} $Z_{jk} = Z_{kj} = \tz{jk}{jk}$.
Driving-point impedances are obtained by (\ref{eq:transz}) from the `principal' second cofactor $C_{jk,jk}(Y)$ for each node pair $(j,k)$ (after dividing by $c(Y)$), and by symmetry, are independent of the polarity of the source connected between nodes $j$ and $k$.
One also has the following identity, by which any transfer impedance can be expressed as sums and differences of driving-point impedances:
\begin{equation}
2 \tz{pq}{jk} = Z_{pk} + Z_{qj} - Z_{pj} - Z_{qk}.
\label{eq:transimp}
\end{equation}
This is readily established using the (anti)symmetry and transitivity properties of cofactors:
\begin{equation}
\begin{split}
\lefteqn{C_{pk,pk} + C_{qj,qj} - C_{pj,pj} - C_{qk,qk}} \\
   &= C_{pk,pk} - C_{qk,pk} + C_{qk,pk} - C_{qk,qk} + C_{qj,qj} - C_{pj,qj} + C_{pj,qj} - C_{pj,pj} \\
   &= C_{pq,pk} + C_{qk,pq} + C_{qp,qj} + C_{pj,qp} = C_{pq,jp} - C_{pq,kp} + C_{pq,jq} - C_{pq,kq}
   = 2 C_{pq,jk}.
\end{split}
\label{eq:transcof}
\end{equation}
\item[The Jacobi identity.]
Applying the identity (\ref{eq:transimp}) three times with permuted indices leads to the following `Jacobi identity' for transfer impedances (and indeed for cofactors more generally):
\begin{equation}
\tz{pq}{jk} + \tz{jp}{qk} + \tz{qj}{pk} = 0.
\label{eq:jacobi}
\end{equation}
Observe that in this formula the index $k$ remains in a fixed position while the others undergo a cyclic permutation.
Thus in electrical terms, (\ref{eq:jacobi}) asserts a property of arbitrary triplets of nodes $(p, q, j)$ in a network grounded at node $k$.
It states that when a source with fixed current $I$ is connected between successsive pairs of nodes in the triplet, a voltage is produced at each remaining node such that the three voltages produced in this manner sum to zero.
If any two indices in (\ref{eq:jacobi}) are equal (say $j = k$ or $p = q$) it reduces to the reciprocity property. \item[Foster's mean impedance theorem.]
An elegant formula of R.\ M.\ Foster \cite{f:aioen} gives the weighted average of driving-point impedances across all branches in a network.
(For Foster's better known `reactance theorem', see Proposition \ref{prop:reactance} in Section \ref{sec:acspecial} below.)
For all node pairs $j, k \in [1,n]$ with $j \neq k$, let $y_{jk} = - Y_{jk}$ denote the total admittance of all branches with nodes $j$ and $k$ as endpoints (or zero), and let $y_{kk} = 0$ for all $k$.
Fix any $p \in [1,n]$ and any $j \in [1,n]$ distinct from $p$, and evaluate $c(Y) = \det Y_{(p,p)}$ by cofactor expansion along row $j$ of $Y$ (recalling this is row $j'$ of $Y_{(p,p)}$, and will involve second cofactors of $Y$ according to (\ref{eq:cof2a})).
This leads to
\begin{multline}
c(Y) = \mathop{\sum_{k=1}^n}_{k \neq p} Y_{jk} C_{pj,pk}(Y)
   = \paren{y_{j1} + y_{j2} + \ldots + y_{jp} + \ldots} C_{pj,pj}(Y)
      - \mathop{\sum_{k=1}^n}_{k \neq p,j} y_{jk} C_{pj,pk}(Y) \\
   = y_{jp} C_{pj,pj}(Y) + \mathop{\sum_{k=1}^n}_{k \neq p,j} y_{jk} \paren{C_{pj,pj}(Y) - C_{pj,pk}(Y)}
   = \sum_{k=1}^n y_{jk} C_{pj,kj}(Y),
\label{eq:yz1}
\end{multline}
hence if $Y$ is of full rank $n - 1$ (and invoking symmetry),
\begin{equation}
\sum_{k=1}^n y_{jk} \tz{jp}{jk} = 1 \qquad \forall j, p \in [1,n] \text{ with } p \neq j.
\label{eq:yz}
\end{equation}
(This identity (\ref{eq:yz}) should not be particularly surprising: by virtue of the definition of $\tz{jp}{jk}$, it is Kirchhoff's current law at node $j$ in disguise.)
For a fixed $p \in [1,n]$ (\ref{eq:yz}) gives $n - 1$ identities, one for each $j$ distinct from $p$.
If one adds together all of these and uses (\ref{eq:transimp})---observing $Z_{jj} = 0$ for all $j$ by definition---there results
\begin{multline}
n - 1
   = \frac{1}{2} \mathop{\sum_{j=1}^n}_{j \neq p} \sum_{k=1}^n y_{jk} \paren{Z_{jk} + Z_{pj} - Z_{pk}} \\
   = \frac{1}{2} \mathop{\sum_{j,k = 1}^n}_{j \neq p} y_{jk} Z_{jk}
      + \frac{1}{2} \mathop{\sum_{j=1}^n}_{j \neq p} Z_{pj} \sum_{k=1}^n y_{jk}
      - \frac{1}{2} \sum_{k=1}^n Z_{pk} \mathop{\sum_{j=1}^n}_{j \neq p} y_{jk} \\
   \displaybreak \\
   = \frac{1}{2} \mathop{\sum_{j,k = 1}^n}_{j \neq p} y_{jk} Z_{jk}
      + \frac{1}{2} \mathop{\sum_{j=1}^n}_{j \neq p} Y_{jj} Z_{pj}
      - \frac{1}{2} \sum_{k=1}^n \paren{Y_{kk} - y_{pk}} Z_{pk}
   = \frac{1}{2} \sum_{j,k = 1}^n y_{jk} Z_{jk},
\label{eq:foster1}
\end{multline}
observing that the third summation contributes the `missing' terms $y_{pk} Z_{pk}$ in the first, while all other terms in the second and third summations cancel (since $Z_{pp} = 0$).
The resulting formula is independent of the chosen value of $p$, and so produces a single identity for the network.
Since $y_{jk} Z_{jk} = y_{kj} Z_{kj}$, it is conveniently restated in the form
\begin{equation}
\mathop{\sum_{j,k = 1}^n}_{j < k} y_{jk} Z_{jk} = n - 1 \qquad \text{or} \qquad
\sum_{\alpha} y_{\alpha} Z_{\alpha} = n - 1,
\label{eq:foster}
\end{equation}
where the second sum is over all network branches $\alpha$, and $Z_{\alpha}$ is the effective impedance between the endpoints of branch $\alpha$.
This is Foster's (1948) `first theorem' as usually stated.
\end{description}
A comment is also called for regarding Tellegen's theorem, the other key property of electrical networks, which asserts that the sum of $v_{\alpha} i_{\alpha}$ over all branches $\alpha$, \emph{including} sources, is zero.
This does not rely on the solution (\ref{eq:nodev}) but is a pure consequence of Kirchhoff's laws.
Indeed, let $\gfull$ be the graph incidence matrix as defined in the Introduction, but for the network including source branches, and let vectors $\ifull$ and $\vfull$ denote the full set of branch currents and branch voltage drops respectively.
Then Kirchhoff's current law asserts that $\gfull \ifull = 0$, while the definition of $\gfull$ entails that $\vfull = \gfull^T v$ where $v$ is the vector of node voltages as above.
Tellegen's theorem is
\begin{equation}
\sum_{\alpha} v_{\alpha} i_{\alpha} = \vfull^T \ifull = v^T \gfull \ifull = 0.
\label{eq:tellegen}
\end{equation}
Recast into the above framework based on $Y$ for the source-free network, and sources represented by external current injections $i$, Tellegen's theorem asserts that the quantity
\begin{equation}
P_{\mathrm{tot}} = v^T i = v^T Y v
\label{eq:ptot}
\end{equation}
represents both the net power injected into the network from external sources, and the net power absorbed within the network described by $Y$.
To see why (\ref{eq:ptot}) holds in general, observe that it certainly holds in the case where the network solution is grounded at some arbitrary node $k$ and all sources have their negative terminal at $k$.
Then observe that any source connected between nodes $p, q \neq k$ can be replaced with equivalent sources connected between $p, k$ and $q, k$ and having source current $i_{\alpha}, -i_{\alpha}$ respectively.
The combined source power $P = (v_p - v_q) i_{\alpha}$ is the same as that for the original source.
The general result follows by superposition, and the fact that adding the same amount to all voltages in $v$ leaves the quantity $v^T i$ unchanged, so the restriction to a grounded solution $v$ can likewise be relaxed.

Under the assumption of a solution $\vgnd{k}$ grounded at some node $k$, (\ref{eq:ptot}) becomes
\begin{equation}
P_{\mathrm{tot}} = \hat{v}^T Y_{(k,k)} \hat{v},
\label{eq:ptotk}
\end{equation}
where $\hat{v}$ is the truncated form of $\vgnd{k}$ omitting the (null) $k$th element.
Provided $Y$ is of maximum rank, this is a positive definite quadratic form with determinant $c(Y)$.
The reasoning above also shows that the elements of the solution $\vgnd{k}$ are given by (\ref{eq:isoln}) when $I$ is a source current connected between nodes $p$ and $q = k$, leading to the general formula for a grounded solution
\begin{equation}
\hat{v} = T \hat{\imath} \qquad \text{where} \qquad T_{pq} = T_{qp} = \tz{pk}{qk}
\qquad \text{and} \qquad T = Y_{(k,k)}^{-1},
\label{eq:gsoln}
\end{equation}
and to the alternative formula
\begin{equation}
P_{\mathrm{tot}} = \hat{v}^T \hat{\imath} = \hat{\imath}^T T \hat{\imath},
\label{eq:ptoti}
\end{equation}
where $\hat{\imath}$ is the truncated source current vector conforming to $\hat{v}$, that is to say with the $k$th element excluded.
Observe that the diagonal elements $T_{pp}$ of $T$ are equal to the driving-point impedances $Z_{pk}$ for each node $p \neq k$, and that an arbitrary driving-point impedance $Z_{pj}$ is (after substituting $q = k$ in (\ref{eq:transimp})) given from $T$ by
\begin{equation}
Z_{pj} = T_{pp} + T_{jj} - 2 T_{pj}.
\label{eq:gsolnz}
\end{equation}
A development of the formal theory of transfer impedances based on Proposition \ref{prop:transz} forms part of the famous paper of Brooks, Smith, Stone and Tutte on `squared rectangles' \cite{bsst:doris}, in which most of the above results can be found.
Brooks et al.\ considered the enumeration of `nets' (networks in which all admittances are equal to 1) in which nodes $p$ and $q$ are fixed and the source current set to $I = c(Y) / \rho$, where $\rho$ is the greatest common divisor of all second cofactors $C_{pq,**}$.
This choice of $I$ and the use of unit admittances ensures that all currents and voltages are integers, `reduced' in the sense of having no common divisor.
As $c(Y)$ is an inconsequential scaling factor in this theory, Brooks et al.\ defined a `transpedance' as a second cofactor $C_{pq,jk}$, which they denoted $[pq,jk]$.

\section{Modification of Networks}

The above theory, connecting cofactors of the admittance matrix $Y$ with properties of the network, has at least one significant gap.
Throughout, it has been necessary to assume that the $Y$ matrix has the maximal rank ($n - 1$), ensuring that the common value $c(Y)$ of all first cofactors is nonzero.
In fact, this turns out to be generically true provided the network is connected, in the graph-theoretic sense---that is, if a path exists between any pair of network nodes.

Establishing the \emph{converse} of this claim is straightforward.
If the network is \emph{not} connected, then without loss of generality one may find some $m < n$ and (re)label the nodes so that the first $m$ belong to one or more connected components separate from the remaining $n - m$ nodes (in the sense that there is no branch connecting a node in the first group to one in the second).
Then by construction of $Y$ it is seen that $Y$ itself separates into a block diagonal matrix with at least two components, viz
\begin{equation}
Y = \mat{Y_1 & 0 \\ 0 & Y_2},
\label{eq:sep}
\end{equation}
where $Y_1$ is of dimensions $m \times m$ and $Y_2$ of dimensions $(n - m) \times (n - m)$.
Considering the respective subnetworks, it is seen that the matrices $Y_1$ and $Y_2$ are of rank no greater than $m - 1$ and $n - m - 1$ respectively, so that $Y$ is of maximum rank $n - 2$.
Accordingly, any network on $n$ nodes that is not connected fails to have an admittance matrix of rank $n - 1$.

Establishing the original claim, on the other hand, requires understanding in some detail how the matrix $Y$ and its cofactors change when various operations are performed on the network.
The following results, stated in \cite{bsst:doris}, prove sufficient to describe arbitrary network changes in terms of a small number of elementary operations.
Two of these---termed \emph{expansion} and \emph{contraction}---are treated in detail, while the third---termed \emph{augmentation}---can be viewed in its own right or as derived from the first two.

\subsection{Expansion: Adding One Node and Branch}

The first elementary operation, that of expansion, involves choosing some $k \in [1,n]$, creating a new node with index $\nu = n + 1$, and creating a new branch between nodes $k$ and $\nu$ with a freely chosen admittance, denoted $y_+$.

If the original network has admittance matrix $Y$, then the admittance matrix $Y^+$ of the modified network differs from $Y$ precisely as follows:
\begin{enumerate}
\item
$Y^+$ has an additional row and column with index $\nu$.
\item
If $y_+ \neq 0$ then row $\nu$ has exactly two nonzero elements, $Y_{\nu\nu}^+ = y_+$ and $Y_{\nu k}^+ = - y_+$ (and the equivalent holds for column $\nu$).
\item
The minor $Y_{(\nu,\nu)}^+$ is identical to $Y$ except for the $k$th element on the diagonal, where one has $Y_{kk}^+ = Y_{kk} + y_+$.
\end{enumerate}
The following results are consequences of these properties of $Y^+$.
\begin{proposition}
\label{prop:expcof1}
Let a network with admittance matrix $Y$ be expanded as above, with $Y^+$ being the new admittance matrix and $y_+$ being the admittance of the new branch.
Then the common value of all first cofactors of $Y^+$ is
\begin{equation}
c(Y^+) = y_+ \cdot c(Y).
\label{eq:expcof1}
\end{equation}
In particular, if $Y$ is of maximum rank $(n - 1)$ and $y_+ \neq 0$, then $Y^+$ is of maximum rank $(\nu - 1)$.
\end{proposition}
This is seen by considering the cofactor $C_{\nu, k}(Y^+)$: this is $(-1)^{\nu + k}$ times a determinant whose $n$th column is empty save one element in position $k$ equal to $-y_+$.
Hence it is $(-1)^{\nu + k}$ times $(-y_+)$, times the equivalent cofactor which evaluates as $(-1)^{k + n} \det Y_{(k,k)}$.
Given that $n$ and $\nu$ differ by 1 and that $\det Y_{(k,k)} = C_{k,k}(Y) = c(Y)$, one has $C_{\nu,k}(Y^+) = y_+ \cdot c(Y)$, and this is the common value of all first cofactors of $Y^+$, since $Y^+$ is itself an admittance matrix (it satisfies Proposition \ref{prop:cof1}).

\begin{proposition}
\label{prop:expz}
In a network expanded from node $k$ as above, all second cofactors of $Y^+$ associated with the new node and branch are equal to the first cofactor of $Y$:
\begin{equation}
C_{k\nu,j\nu}(Y^+) = C_{j\nu,k\nu}(Y^+) = c(Y), \qquad j \in [1,n].
\label{eq:expz}
\end{equation}
Hence in the new network, the effective impedance $Z_{k\nu}$ across the added branch is identical to the branch impedance $1 / y_+$.
\end{proposition}
Equation (\ref{eq:expz}) follows from observing that the excluded rows and columns for the cofactor include the only ones where $Y^+$ and $Y$ differ, hence $C_{k\nu,j\nu}(Y^+)$ is the same as $(-1)^{k+j} \det Y_{(k,j)} = c(Y)$.
In the case $j = k$, the effective impedance by (\ref{eq:transz}) and (\ref{eq:expcof1}) is $Z_{k\nu} = c(Y) / c(Y^+) = 1 / y_+$, as required.
Of course, Proposition \ref{prop:expz} is intuitively correct given that if a source is connected with one terminal at node $\nu$, there is no other path for the current than through the admittance $y_+$, giving a voltage drop $I / y_+$ between nodes $k$ and $\nu$ where $I$ is the source current.

With a little more work, one may relate other second cofactors of $Y^+$ to those of $Y$.
\begin{proposition}
\label{prop:expzij}
In an expanded network as above, all second cofactors $C_{ij,pq}$ of $Y^+$ with $i,j,p,q \in [1,n]$ are given by
\begin{equation}
C_{ij,pq}(Y^+) = y_+ \cdot C_{ij,pq}(Y), \qquad i, j, p, q \in [1,n].
\label{eq:expzij}
\end{equation}
Consequently, effective impedances $Z_{pq}$ calculated from $Y^+$ are identical to those calculated from $Y$, whenever $p,q \in [1,n]$.
\end{proposition}
In the cofactor $C_{ij,pq}(Y^+)$, neither row $\nu$ nor column $\nu$ are eliminated.
So in the last row of the minor $Y_{(ij,pq)}^+$ there is either a single nonzero element or two nonzero elements, depending whether either of $p$ or $q$ are equal to $k$.
In case either equal $k$, the sole nonzero element corresponds to $Y_{\nu\nu}^+ = y_+$, and expanding by cofactors yields $\det Y_{(ij,pq)}^+ = y_+ \cdot \det Y_{(ij,pq)}$.
Identity (\ref{eq:expzij}) is an immediate consequence.
Similarly, if either of the row indices $i$ or $j$ equals $k$, expanding down the last column yields the same result.

This leaves the case where none of $i,j,p,q$ are equal to $k$.
For brevity, let $U = Y_{(ij,pq)}$ and $U^+ = Y_{(ij,pq)}^+$ denote the relevant minors of $Y$ and $Y^+$ respectively, and observe that $C_{ij,pq}(Y) = (-1)^s \det U$ and $C_{ij,pq}(Y^+) = (-1)^s \det U^+$ where $s = \sigma(i,j) + \sigma(p,q)$.
It therefore, again, suffices to show that
\begin{equation}
\det U^+ = y_+ \cdot \det U.
\label{eq:expdetu}
\end{equation}
This may be done from first principles in several equivalent ways.
Perhaps the most succinct is to expand $\det U^+$ using the second-order Laplace formula (\ref{eq:laplace2}).
For convenience, let $n' = n - 2$ denote the size of the square matrix $U$, and $\nu' = n - 1$ the size of $U^+$.
The row/column index corresponding to $k$ in both $U$ and $U^+$ will be denoted $k'$.
Now, setting $A = U^+$, $j = \nu'$ and $p = k'$ in formula (\ref{eq:laplace2}) gives
\begin{equation}
\det U^+
   = \sum_{r=1}^{\nu'} \mathop{\sum_{t=1}^{\nu'}}_{t \neq r} u_{\nu' r}^+ u_{k' t}^+ C_{\nu' k', r t}(U^+).
\label{eq:expdetup}
\end{equation}
Observe that in this double sum, the only nonzero terms will be those corresponding to $r = \nu'$ (with $u_{\nu' r}^+ = y_+$) and to $r = k'$ (with $u_{\nu' r} = -y_+$).
Writing these out explicitly gives
\begin{equation}
\det U^+ = y_+ \sum_{t=1}^{n'} u_{k' t}^+ C_{\nu' k', \nu' t}(U^+)
   - y_+ \mathop{\sum_{t=1}^{\nu'}}_{t \neq k'} u_{k' t}^+ C_{\nu' k', k' t}(U^+).
\label{eq:expdetup2}
\end{equation}
In the first sum in (\ref{eq:expdetup2}), the elements $u_{k' t}^+$ of $U^+$ are all identical to the corresponding elements $u_{k' t}$ of $U$, with the exception of $u_{k' k'}^+ = u_{k' k'} + y_+$.
In the second sum, all terms other than that with $t = \nu'$ are zero (because $C_{\nu' k', k' t}(U^+)$ with $t < \nu'$ involves a minor whose rightmost column is all zeros), and the $t = \nu'$ term involves $u_{k' \nu'}^+ = -y_+$.
Taking all this into account, (\ref{eq:expdetup2}) is equivalent to
\begin{equation}
\det U^+ = y_+ \sum_{t=1}^{n'} u_{k' t} C_{\nu' k', \nu' t}(U^+) + y_+^2 \cdot C_{\nu' k', \nu' k'}(U^+)
   + y_+^2 \cdot C_{\nu' k', k' \nu'}(U^+).
\label{eq:expdetup3}
\end{equation}
The two terms in $y_+^2$ immediately cancel, owing to the antisymmetry of the second cofactors.
Lastly, one observes that
\begin{equation}
C_{\nu' k', \nu' t}(U^+) = (-1)^{2\nu' + k' + t} \det U_{(\nu' k', \nu' t)}^+ = (-1)^{k' + t} \det U_{(k',t)}
   = C_{k',t}(U)
\label{eq:expcofu}
\end{equation}
and so (\ref{eq:expdetup3}) is revealed as $y_+$ times an ordinary Laplace expansion of $\det U$, as required:
\begin{equation}
\det U^+ = y_+ \sum_{t=1}^{n'} u_{k' t} C_{k',t}(U) = y_+ \cdot \det U.
\label{eq:expdetupu}
\end{equation}

\begin{proposition}
\label{prop:expzz}
In a network expanded from node $k$ as above, all second cofactors of $Y^+$ of the form $C_{p\nu,q\nu}$ with $p,q \in [1,n]$ are given by
\begin{equation}
C_{p\nu,q\nu}(Y^+) = c(Y) + y_+ \cdot C_{pk,qk}(Y), \qquad p, q \in [1,n].
\label{eq:expzz}
\end{equation}
In particular, for any $p \in [1,n]$ one has $Z_{p\nu} = 1/y_+ + Z_{pk}$.
\end{proposition}
Note that in the case $p = k$ or $q = k$ this reduces to the result in Proposition \ref{prop:expz}, which is consistent with (\ref{eq:expzz}) by virtue of the convention $C_{kk,**} = 0$.
Assume then that both $p$ and $q$ are distinct from $k$.
The cofactor $C_{p\nu,q\nu}$ then reduces to an equivalent cofactor of the principal minor $Y_{(\nu,\nu)}^+$:
\begin{equation}
C_{p\nu,q\nu}(Y^+) = C_{\nu p,\nu q}(Y^+) = C_{p,q}(Y_{(\nu,\nu)}^+).
\label{eq:expzz1}
\end{equation}
Observe now that $Y_{(\nu,\nu)}^+$, differing from $Y$ in just one element, can be written using an earlier notation as
\begin{equation}
Y_{(\nu,\nu)}^+ = Y +_k \ycol{k}{y_+ \! e_k}
\label{eq:ynn}
\end{equation}
where \ycol{k}{y_+ \! e_k} is the matrix formed from $Y$ by replacing the $k$th column with an elementary vector $y_+ \! e_k$ having admittance $y_+$ in the $k$th position and zeroes elsewhere.
This vector, in column $k$, represents the difference between $Y_{(\nu,\nu)}^+$ and $Y$.

It is now readily verified that any $(p,q)$-minor of the matrix in (\ref{eq:ynn}), where $p \neq k$ and $q \neq k$, is also equal to a composition analogous to (\ref{eq:ynn}).
In strict terms, this composition is
\begin{equation}
Y_{(\nu,\nu)(p,q)}^+ = Y_{(p,q)} +_{k'} \ycol{k'}{y_+ \! e_{k''}}_{(p,q)}
\label{eq:ynnpq}
\end{equation}
where $k' = \sigma(q,k) - q$ is the column index in $Y_{(p,q)}$ corresponding to $k$, and $k'' = \sigma(p,k) - p$ the row index.
Taking determinants, using property \ref{det:addcol}, and expanding down the column $y_+ \! e_{k''}$, one has
\begin{multline}
C_{p,q}(Y_{(\nu,\nu)}^+)
   = (-1)^{p + q} \paren{\det Y_{(p,q)} + (-1)^{k'' + k'} y_+ \det Y_{(p,q)(k'',k')}} \\
   = (-1)^{p + q} \det Y_{(p,q)} + (-1)^{\sigma(p,k) + \sigma(q,k)} y_+ \det Y_{(pk,qk)}
   = C_{p,q}(Y) + y_+ C_{pk,qk}(Y),
\label{eq:yqqdet}
\end{multline}
and the result (\ref{eq:expzz}) follows.

\begin{proposition}
\label{prop:expzk}
In a network expanded from node $k$, all second cofactors of $Y^+$ of the form $C_{k\nu,pq}$ or $C_{pq,k\nu}$ with $p,q \in [1,n]$ are zero:
\begin{equation}
C_{k\nu,pq}(Y^+) = C_{pq,k\nu}(Y^+) = 0, \qquad p, q \in [1,n].
\label{eq:expzk}
\end{equation}
\end{proposition}
The proof here is straightforward: consider the minor $Y_{(k\nu,pq)}^+$ and observe this is formed by eliminating rows $k$ and $\nu$, and \emph{not} eliminating column $\nu$.
Hence it contains a column of all zeros, and consequently the associated cofactor $C_{k\nu,pq}$ is zero.
The result for $C_{pq,k\nu}$ follows by symmetry.

\begin{proposition}
\label{prop:expzn}
In a network expanded from node $k$, all second cofactors of $Y^+$ of the form $C_{j\nu,pq}$ or $C_{pq,j\nu}$ with $j,p,q \in [1,n]$ are given by
\begin{equation}
C_{j\nu,pq}(Y^+) = C_{pq,j\nu}(Y^+) = y_+ \cdot C_{jk,pq}(Y), \qquad j, p, q \in [1,n].
\label{eq:expzn}
\end{equation}
\end{proposition}
Once again, this proposition has an earlier one as a special case, reducing to Proposition \ref{prop:expzk} when $j = k$ by virtue of the identity $C_{kk,**}(Y) = 0$.
In all other cases with $j \neq k$, the minor $Y_{(j\nu,pq)}^+$ has a single entry $-y_+$ in the rightmost column, in row $k' = \sigma(j,k) - j$..
Forming the cofactor strikes out row $k'$, leaving the minor $Y_{(jk,pq)}$ of the original matrix $Y$.
Thus
\begin{multline}
C_{j\nu,pq}(Y^+)
   = (-1)^{\sigma(j,\nu) + \sigma(p,q)} \cdot (-y_+) \cdot (-1)^{k' + \nu - 2} \det Y_{(jk,pq)} \\
   = y_+ \cdot (-1)^{j + k' + \sigma(p,q)} \det Y_{(jk,pq)}
   = y_+ \cdot (-1)^{\sigma(j,k) + \sigma(p,q)} \det Y_{(jk,pq)} = y_+ \cdot C_{jk,pq}(Y).
\label{eq:expzn1}
\end{multline}
The result for $C_{pq,j\nu}(Y^+)$ also follows by symmetry.

Propositions \ref{prop:expcof1} through \ref{prop:expzn} allow all first and second cofactors of $Y^+$ to be determined from those of $Y$  thus characterising the properties of the expanded network in terms of the original network.
A key consequence is that when the original network has a $Y$ matrix of maximum rank, Proposition \ref{prop:expcof1} guarantees expansion with a nonzero admittance always increases the rank by one.

\subsection{Contraction: Identifying Two Nodes}

The second elementary operation, that of contraction, involves choosing distinct values $j, k \in [1,n]$ (where $j < k$ without loss of generality) and identifying the nodes with indices $j$ and $k$, deleting any branches with $(j,k)$ as endpoints and leaving node $j$ connected to all other branches that were incident to node $j$ or $k$ in the original network.

If the original network has admittance matrix $Y$, then the admittance matrix $Y^-$ of the modified network, of size $n - 1$, is formally constructed from $Y$ by means of the following recipe:
\begin{enumerate}
\item
Form matrix $Y^a$ from $Y$ by replacing row $j$ with the sum of rows $j$ and $k$ (thus $Y_{jq}^a = Y_{jq} + Y_{kq}$ for $q \in [1,n]$).
Note that this operation eliminates any contribution in row $j$ from branches connected between nodes $j$ and $k$.
\item
Form matrix $Y^b$ from $Y^a$ by replacing column $j$ with the sum of columns $j$ and $k$ (thus $Y_{qj}^b = Y_{qj}^a + Y_{qk}^a$ for $q \in [1,n]$).
Observe that $Y_{jj}^b = Y_{jj} + Y_{kk} + Y_{jk} + Y_{kj}$ is the sum of admittances of branches incident with either node $j$ or node $k$, but omitting any branch between $j$ and $k$.
Observe also that $Y^b$ is symmetric, and $Y$ and $Y^b$ differ only in row and column $j$.
\item
Delete the $k$th row and column from $Y^b$, forming $Y^-$ as the principal minor $Y_{(k,k)}^b$.
\end{enumerate}
The next proposition follows immediately from the above recipe and property \ref{det:combcol} of determinants---that adding one row or column of an arbitrary square matrix $A$ to another leaves $\det A$ unchanged.
\begin{proposition}
\label{prop:conminor}
For some $m \in [1,n]$, let $\alpha$ and $\beta$ denote (ordered) $m$-tuples of indices in $[1,n]$.
Let the network on $n$ nodes with admittance matrix $Y$ be contracted on nodes $j$ and $k$ as above.
Suppose $m$-tuple $\alpha$ either includes $j$ or excludes $k$ (or both).
Suppose the same is separately true of $\beta$.
Then
\begin{equation}
\det Y_{(\alpha,\beta)}^b = \det Y_{(\alpha,\beta)} \qquad \text{and} \qquad
C_{\alpha,\beta}(Y^b) = C_{\alpha,\beta}(Y)
\label{eq:conminor}
\end{equation}
where $Y_{(\alpha,\beta)}$ denotes the (generalised) minor of $Y$ formed by deleting rows with indices in $\alpha$ and columns with indices in $\beta$, $C_{\alpha,\beta}(Y)$ is the associated cofactor, and $Y^b$ is defined as above.
\end{proposition}
As a corollary one directly obtains the following:
\begin{proposition}
\label{prop:concof1}
Let the network with admittance matrix $Y$ be contracted on arbitrary nodes $j$ and $k$, and let $Y^-$ be the admittance matrix of the contracted network.
Then all first cofactors of $Y^-$ are equal to the principal second cofactor of $Y$ associated with indices $j$ and $k$:
\begin{equation}
C_{p,q}(Y^-) = C_{jk,jk}(Y) = C_{kj,kj}(Y), \qquad p, q \in [1, n - 1].
\label{eq:concof1}
\end{equation}
\end{proposition}
Proposition \ref{prop:concof1} follows from Proposition \ref{prop:conminor} upon setting $\alpha = \beta = (jk)$, and noting that the minor $Y_{(j,j)}^-$ is identical to $Y_{(jk,jk)}^b$.
This shows that $C_{j,j}(Y^-) = C_{jk,jk}(Y^b) = C_{jk,jk}(Y)$; other first cofactors also take this value owing to the structure of $Y^-$ as an admittance matrix.

The admittance matrix $Y^-$ has maximum possible rank $n - 2$, where $n$ is the size of the original network.
Proposition \ref{prop:concof1} states that $Y^-$ has this maximum rank if and only if $C_{jk,jk}(Y) \neq 0$.

Turning to the second cofactors of $Y^-$, certain of these---those with $j$ as both row and column indices---can with the aid of Proposition \ref{prop:conminor} be recognised almost immediately as third cofactors of the original matrix $Y$.
\begin{proposition}
\label{prop:concof2j}
Let the network with admittance matrix $Y$ be contracted on arbitrary nodes $j$ and $k$ with $j < k$, forming a network with admittance matrix $Y^-$.
Then all second cofactors of $Y^-$ involving row and column $j$ (with $j$ the index of the coalesced node) are given by equivalent third cofactors of $Y$ involving, in addition, row and column $k$.
That is, for any $p,q \in [1,n]$ not equal to $j$ or $k$, and $p' = \sigma(k,p) - k$, $q' = \sigma(k,q) - k$ the equivalent indices in $Y^-$, one has
\begin{equation}
C_{jp',jq'}(Y^-) = C_{kjp,kjq}(Y), \qquad p,q \in [1,n] \text{ other than } j,k.
\label{eq:concof2j}
\end{equation}
\end{proposition}
The third cofactor in (\ref{eq:concof2j}) is defined in the usual recursive manner as $C_{kjp,kjq}(Y) = C_{jp',jq'}(Y_{(k,k)})$, with the sign unchanged as the eliminated row and column have the same index.
(Recall the assumption that $j < k$, so that index $j$ represents the same node in $Y$ and $Y_{(k,k)}$.)
This reveals Proposition \ref{prop:concof2j} to be a direct corollary of Proposition \ref{prop:conminor} with $m = 3$:
\begin{equation}
C_{jp',jq'}(Y^-) = C_{jp',jq'}(Y_{(k,k)}^b) = C_{kjp,kjq}(Y^b) = C_{kjp,kjq}(Y).
\label{eq:concof2a}
\end{equation}

Obtaining a more general formula for the second cofactors of $Y^-$ requires the following identity, which is a refinement (of a particular case) of Sylvester's identity for determinants.
\begin{proposition}
\label{prop:sylvcof}
Let $p, q, r, s \in [1,n]$ and $m \in [0, n - 2]$, and let $\alpha$ and $\beta$ denote ordered $m$-tuples (empty if $m = 0$) of indices in $[1,n]$ such that $p, q \not\in \alpha$ and $r, s \not\in \beta$.
Then the following is an identity for the cofactors of any $n \times n$ matrix $A$:
\begin{equation}
C_{p\alpha,r\beta}(A) \, C_{q\alpha,s\beta}(A) - C_{p\alpha,s\beta}(A) \, C_{q\alpha,r\beta}(A)
   = C_{pq\alpha,rs\beta}(A) \, C_{\alpha,\beta}(A).
\label{eq:sylvcof}
\end{equation}
By convention, when $m = 0$ the cofactor $C_{\alpha,\beta}(A)$ with empty indices is equal to $\det A$, while when $m = n - 2$ the cofactor $C_{pq\alpha,pq\beta}(A)$ with full indices is equal to 1.
\end{proposition}
Proposition \ref{prop:sylvcof} is a consequence of Sylvester's identity \cite[\S 0.8.6]{hj:ma}, which is generally stated in terms of minors of a matrix $A$ that include (not exclude) specific rows and columns.
Accordingly, for index sets $\alpha$ and $\beta$ let $A_{[\alpha,\beta]}$ (with square brackets) denote the minor of $A$ obtained by intersecting rows having indices in $\alpha$ with columns having indices in $\beta$.
Given an $m$-tuple $\alpha$, Sylvester's identity is
\begin{equation}
\det B = \det A_{[\alpha,\alpha]}^{n - m - 1} \det A \qquad \text{where} \qquad
b_{jk} = \det A_{[\alpha j', \alpha k']}, \quad j, k \in [1, n - m],
\label{eq:sylvester}
\end{equation}
in which $j'(j)$ and $k'(k)$ range over the $n - m$ indices not included in $\alpha$.
A concise proof of (\ref{eq:sylvester}) can be found in \cite{bs:di}.

Specialising to the case $n - m = 2$, with $j,k$ now the two indices excluded from $\alpha$, and restating in terms of excluded rows and columns, (\ref{eq:sylvester}) becomes
\begin{equation}
\det A_{(j,j)} \det A_{(k,k)} - \det A_{(j,k)} \det A_{(k,j)} = \det A_{(jk,jk)} \det A.
\label{eq:sylv2}
\end{equation}
Now, as a further refinement, suppose the matrix $A$ is itself extended by $m$ `silent' rows and columns, all of which are omitted from the determinants in (\ref{eq:sylv2}).
If $\gamma$ denotes the set of indices of these non-participating rows and columns, (\ref{eq:sylv2}) extends to the following identity for the larger matrix $A'$:
\begin{equation}
\det A_{(j\gamma,j\gamma)}' \det A_{(k\gamma,k\gamma)}'
   - \det A_{(j\gamma,k\gamma)}' \det A_{(k\gamma,j\gamma)}'
   = \det A_{(jk\gamma,jk\gamma)}' \det A_{(\gamma,\gamma)}'.
\label{eq:sylv2m}
\end{equation}
Since all the minors in (\ref{eq:sylv2m}) are principal apart from the symmetric pair $A_{(j\gamma,k\gamma)}'$ and $A_{(k\gamma,j\gamma)}'$, one may restate it immediately in terms of cofactors:
\begin{equation}
C_{j\gamma,j\gamma}(A') \, C_{k\gamma,k\gamma}(A')
   - C_{j\gamma,k\gamma}(A') \, C_{k\gamma,j\gamma}(A')
   = C_{jk\gamma,jk\gamma}(A') \, C_{\gamma,\gamma}(A').
\label{eq:sylv2cof}
\end{equation}
To obtain the identity (\ref{eq:sylvcof}) from this, it only remains to show that (\ref{eq:sylv2cof}) is invariant under arbitrary permutations of the rows and columns of $A'$.
Thus, let $\pi$ denote some permutation of the $n$ rows of $A'$, and $\rho$ another permutation on the $n$ columns, such that the resulting rearrangement maps $A'$ onto the (new) matrix $A$.
Under the postulated invariance, identity (\ref{eq:sylvcof}) results with
\begin{equation}
p = \pi(j), \quad q = \pi(k), \quad r = \rho(j), \quad s = \rho(k), \quad
\alpha = \pi(\gamma), \quad \beta = \rho(\gamma).
\label{eq:pirho}
\end{equation}
The required invariance is readily shown using the transformation rule for cofactors
\begin{equation}
C_{\pi(\mu), \rho(\nu)}(A) = \sgn \pi | \mu \cdot \sgn \rho | \nu \cdot C_{\mu,\nu}(A')
\label{eq:trcof}
\end{equation}
where $\mu$ and $\nu$ are arbitrary index sets, and the notation $\pi | \mu$ denotes the restriction of $\pi$ to a permutation on the elements of $\mu$ alone.
(So for example, if $\pi$ is a transposition of rows $j$ and $k$, then a cofactor involving both rows $j$ and $k$ changes sign when going from $A'$ to $A$, but a cofactor not involving both rows will be unchanged.)
Rule (\ref{eq:trcof}) maps formula (\ref{eq:sylvcof}) onto identity (\ref{eq:sylv2cof}) given the restricted permutations have the property that
\begin{equation}
\sgn \pi | j \gamma \cdot \sgn \pi | k \gamma = \pm \sgn \pi | j k \gamma \cdot \sgn \pi | \gamma
\label{eq:permsgn}
\end{equation}
(and likewise for $\rho$), where the negative sign on the right hand side applies if $\pi$ exchanges the order of $j$ and $k$, and the positive sign applies in all other cases.
(To verify (\ref{eq:permsgn}) it suffices to verify it for any single transposition, given any permutation is a product of transpositions, and sign is a multiplicative property.)
The cases where $j$ and $k$ are reversed by exactly one of $\pi$ and $\rho$, producing a sign change on the right hand side of (\ref{eq:sylvcof}), are precisely the cases where there is a compensating sign change due to the reversal of terms in the expression on the left hand side of (\ref{eq:sylvcof}).
This completes the proof of Propostion \ref{prop:sylvcof}.

Applying Proposition \ref{prop:sylvcof} to carefully chosen minors of the modified admittance matrix $Y^b$ provides the desired general formula, as stated in \cite{bsst:doris}, that expresses the second cofactors of $Y^-$ for a contracted network.
\begin{proposition}
\label{prop:concof2}
Let the network with admittance matrix $Y$ be contracted on arbitrary nodes $j$ and $k$, forming a network with admittance matrix $Y^-$.
Then the second cofactors of $Y^-$ are determined from those of $Y$ by the formula
\begin{equation}
C_{p'q',r's'}(Y^-) = \frac{1}{c(Y)} \bbrak{C_{pq,rs}(Y) \, C_{jk,jk}(Y) - C_{pq,jk}(Y) \, C_{rs,jk}(Y)}
\label{eq:concof2}
\end{equation}
where $c(Y)$ is the common value of all first cofactors of $Y$, and $p',q',r',s'$ are indices of $Y^-$ respectively equivalent to indices $p,q,r,s \neq k$ of $Y$ (thus if $j < k$ then $p' = \sigma(k,p) - k$, etc.).
\end{proposition}
To see how (\ref{eq:concof2}) follows from (\ref{eq:sylvcof}), observe first that (\ref{eq:concof2}) satisfies the antisymmetry property of second cofactors.
Using this property, one may assume without loss of generality that neither $q$ nor $s$ is equal to $j$, hence that the minor $Y_{(q,s)}^b$ includes both rows and columns with index $j$ and $k$.
Apply Proposition \ref{prop:sylvcof} to the matrix $Y^b$ with $\alpha = (q)$, $\beta = (s)$, and $k$ in place of indices $q$ and $s$ in (\ref{eq:sylvcof}): this formula then becomes
\begin{equation}
C_{pkq,rks}(Y^b) \, C_{q,s}(Y^b) = C_{pq,rs}(Y^b) \, C_{kq,ks}(Y^b) - C_{pq,ks}(Y^b) \, C_{kq,rs}(Y^b).
\label{eq:conscof}
\end{equation}
On the left hand side of (\ref{eq:conscof}), $C_{pkq,rks}(Y^b) = C_{kpq,krs}(Y^b)$ is readily seen to be equal to $C_{p'q',r's'}(Y^-)$ after exclusion of row and column $k$; meanwhile, $C_{q,s}(Y^b) = c(Y)$ by Proposition \ref{prop:conminor} since neither $q$ nor $s$ is equal to $k$.
On the right hand side, $C_{pq,rs}(Y^b) = C_{pq,rs}(Y)$ again by Proposition \ref{prop:conminor} since none of $p,q,r,s$ is equal to $k$, and $C_{kq,ks}(Y) = C_{q',s'}(Y^-) = C_{jk,jk}(Y)$ by Proposition \ref{prop:concof1}.

This leaves the last two cofactors on the right of (\ref{eq:conscof}), both of which are expressible in the form $C_{**,ku}(Y^b)$ where neither row index equals $k$ and where column index $u$ does not equal $j$.
Because $u \neq j$, the associated minor $Y_{(**,ku)}^b$ will include a column corresponding to $j$, the entries of which are obtained as the sum of columns $j$ and $k$ in the matrix $Y^a$.
By property \ref{det:addcol} of determinants, it follows that
\begin{equation}
\det Y_{(**,ku)}^b = \det Y_{(**,ku)}^a + (-1)^{\sigma(j,u) - \sigma(k,u) - 1} \det Y_{(**,ju)}^a
\label{eq:concofa1}
\end{equation}
where the last term is found by forming the minor with column $k$ in place of column $j$ from $Y^a$, then notionally moving this column to the position occupied by column $j$---the expression in the sign accounting for how the relative values of $u$, $j$ and $k$ change the sign of the determinant taking into account the deletion of columns $j$ and $u$.
Now, since neither deleted row index is equal to $k$, the determinants involving $Y^a$ in (\ref{eq:concofa1}) will be equal to the corresponding determinants involving $Y$ (either by property \ref{det:combcol} applied to rows, or else because one of the excluded rows is actually row $j$).
Finally, multiplying through by $(-1)^{\sigma(k,u)}$ and by the appropriate sign $(-1)^{\sigma(*,*)}$ based on the row indices, one obtains for the equivalent cofactors
\begin{equation}
C_{**,ku}(Y^b) = C_{**,ku}(Y) + (-1) C_{**,ju}(Y) = C_{**,kj}(Y).
\label{eq:concofa2}
\end{equation}
Using (\ref{eq:concofa2}) for the cofactors $C_{pq,ks}(Y^b)$ and $C_{kq,rs}(Y^b) = C_{rs,kq}(Y^b)$ in (\ref{eq:conscof}), and transposing indices one last time (with two cancelling changes of sign), now yields formula (\ref{eq:concof2}).

Note that in the language of transfer impedances, using Proposition \ref{prop:transz} and Proposition \ref{prop:concof1}, Proposition \ref{prop:concof2} becomes
\begin{equation}
\tz{p'q'}{r's'}(Y^-) = \tz{pq}{rs}(Y) - \tz{pq}{jk}(Y) \frac{\tz{rs}{jk}(Y)}{Z_{jk}(Y)},
\label{eq:contz}
\end{equation}
and for the self-impedances,
\begin{equation}
Z_{p'q'}(Y^-) = Z_{pq}(Y) - \frac{\tz{pq}{jk}^2(Y)}{Z_{jk}(Y)}.
\label{eq:conz}
\end{equation}
Thus, in the contracted network, the effect on voltages of a source in the original network is weakened by the effect of an equivalent `shadow source' connected between the contracted nodes in the original network.
The relative magnitude of this weakening is determined by the ratio of voltage induced by the shadow source at the actual source location, relative to that at the shadow source location.
The weakening can also be interpreted as a reduction in the equivalent network impedances.

\subsection{Augmentation: Adding a Branch Between Existing Nodes}

Another elementary operation that readily suggests itself is that of adding a new branch, with admittance $y_+$, between two existing nodes $j$ and $k$.
If $\bar{Y}$ denotes the admittance matrix of the augmented network, the process of going from $Y$ to $\bar{Y}$ can be described in either of two ways:
\begin{enumerate}
\item
As a direct operation on $Y$, that adds $y_+$ to the diagonal elements $Y_{jj}$ and $Y_{kk}$, and subtracts $y_+$ from the elements $Y_{jk}$ and $Y_{kj}$.
\item
As the successive application of the two elementary operations alraedy described: first expand the network from node $k$ with the new admittance $y_+$, then contract the new network by identifying nodes $j$ and $\nu$.
\end{enumerate}
Either method may be used to characterise the cofactors of $\bar{Y}$ in terms of those of $Y$.
For example:
\begin{proposition}
\label{prop:augcof1}
Let the network with admittance matrix $Y$ be augmented by adding a new branch with admittance $y_+$ between nodes $j$ and $k$, where $j,k \in [1,n]$.
Then all first cofactors of the admittance matrix $\bar{Y}$ of the augmented network are given by
\begin{equation}
C_{p,q}(\bar{Y}) = c(Y) + y_+ \cdot C_{jk,jk}(Y), \qquad p, q \in [1,n]
\label{eq:augcof1}
\end{equation}
where $c(Y)$ is the common value of all first cofactors of $Y$.
\end{proposition}
Using method 2 above, viewing the augmentation as an expansion from $k$ followed by a contraction on $j$ and $\nu$, Proposition \ref{prop:augcof1} is quickly obtained as a consequence of Proposition \ref{prop:expzz} and Proposition \ref{prop:concof1}.
It may also be derived by method 1 from direct inspection of $C_{k,k}(\bar{Y})$, by an argument almost identical to that in the proof of Proposition \ref{prop:expzz}.

As a consequence of Proposition \ref{prop:augcof1}, if a network admittance matrix $Y$ is of full rank $n - 1$, then augmenting the network preserves the full-rank property if and only if the expression (\ref{eq:augcof1}) is nonzero.
A sufficient condition is that the quantities $y_+$ and $Z_{jk}$, if nonzero, agree in sign.
(For pure resistive networks, both will normally be positive.)

\begin{proposition}
\label{prop:augimp}
Let a network be augmented by adding a branch between two nodes $j$ and $k$.
Then in the admittance matrix $\bar{Y}$ of the augmented network, second cofactors involving the node pair $(j,k)$ as rows or columns are identical to the corresponding cofactors of the admittance matrix $Y$ of the original network:
\begin{equation}
C_{jk,pq}(\bar{Y}) = C_{pq,jk}(\bar{Y}) = C_{pq,jk}(Y), \qquad p, q \in [1,n].
\label{eq:augimp}
\end{equation}
\end{proposition}
The equivalence of this result for rows and columns follows once again from the fact that admittance matrices are symmetric.
Using method 1 above, the equality of the second cofactors for $Y$ and $\bar{Y}$ is a simple consequence of the construction of $\bar{Y}$: eliminating rows $j$ and $k$ (or columns $j$ and $k$) also eliminates all elements where $\bar{Y}$ and $Y$ differ.

While in this case method 1 leads to the more direct proof of Proposition \ref{prop:augimp}, a proof is also straightforward (albeit longer) using method 2.
If $Y^+$ denotes the expanded admittance matrix before contracting nodes $j$ and $\nu$, the required cofactor of $\bar{Y}$ is given by Proposition \ref{prop:concof2} as
\begin{equation}
C_{pq,jk}(\bar{Y}) = \frac{1}{y_+ \cdot c(Y)}
   \bbrak{C_{pq,jk}(Y^+) \, C_{j\nu,j\nu}(Y^+) - C_{pq,j\nu}(Y^+) \, C_{jk,j\nu}(Y^+)}
\label{eq:augimpp}
\end{equation}
to which Propositions \ref{prop:expzij}, \ref{prop:expzz} and \ref{prop:expzn} may now be applied to express the cofactors of $Y^+$ in terms of cofactors of $Y$:
\begin{equation}
\begin{split}
C_{pq,jk}(\bar{Y}) &= \frac{1}{y_+ \cdot c(Y)}
   \bbrak{y_+ C_{pq,jk}(Y) \paren{c(Y) + y_+ C_{jk,jk}(Y)} - y_+^2 C_{pq,jk}(Y) \, C_{jk,jk}(Y)} \\
   &= C_{pq,jk}(Y).
\end{split}
\label{eq:augimpp2}
\end{equation}

\begin{proposition}
\label{prop:augcof2}
Let the network on $n$ nodes with admittance matrix $Y$ be augmented by adding a branch with admittance $y_+$ between tow nodes $j$ and $k$, giving a network with admittance matrix $\bar{Y}$.
Then for arbitrary indices $p,q,r,s \in [1,n]$ the second cofactors of $\bar{Y}$ are given by
\begin{equation}
C_{pq,rs}(\bar{Y}) = C_{pq,rs}(Y)
   + \frac{y_+}{c(Y)} \bbrak{C_{pq,rs}(Y) \, C_{jk,jk}(Y) - C_{pq,jk}(Y) \, C_{rs,jk}(Y)}
\label{eq:augcof2}
\end{equation}
where $c(Y)$ is the common value of all first cofactors of the original matrix $Y$.
\end{proposition}
Formula (\ref{eq:augcof2}) is readily demonstrated by method 2, by essentially the same procedure as above.
Indeed, Proposition \ref{prop:augimp} is simply the special case $(r,s) = (j,k)$, where the terms in square brackets cancel.
A proof by method 1 in this case is considerably more difficult, and requires resort to additional tools such as Proposition \ref{prop:sylvcof}.

As suggested in \cite{bsst:doris}, a more concise statement of (\ref{eq:augcof2}) for augmented networks is as the formula
\begin{equation}
C_{pq,rs}(\bar{Y}) = C_{pq,rs}(Y) + y_+ \cdot C_{pq,rs}^{[jk]}(Y)
\label{eq:augcof2a}
\end{equation}
where the notation $C_{pq,rs}^{[jk]}(Y)$ stands for the expression on the right hand side of (\ref{eq:concof2}) for a contracted network.

\section{The Kirchhoff Characteristic}
\label{sec:kirchhoff}

\subsection{Definition and Relation to $c(Y)$}

The preceding network analysis based on the admittance matrix $Y$ highlighted the existence of a key property specific to $Y$ and hence to the underlying network.
Denoted $c(Y)$ above, it is the common value of all first cofactors of $Y$, and if nonzero, certifies that $Y$ is of the maximum rank $(n - 1)$ required for the existence of unique network solutions.

It turns out that $c(Y)$ is equivalent to a property of networks that can be defined directly from the network topology itself (that is, from the underlying graph) together with the branch admittances.
The essence of this was described by Kirchhoff, so it will be termed the \emph{Kirchhoff characteristic} here.
(Other terms for this include the `complexity' as used in \cite{bsst:doris}, and the `determinant' reflecting its definition as a cofactor of $Y$.)

As at the outset, consider the underlying graph of the network, with incidence matrix $G$.
(The directed structure on the graph is not required for present purposes, so the elements of $G$ may be considered without regard to sign.)
A \emph{loop} in this graph (or \emph{ring main} in engineering parlance) is a closed path formed by branches in the graph without repetitions, and a \emph{tree} is any connected subgraph containing no loops.
A \emph{spanning tree} is any tree that connects all $n$ nodes; any connected graph contains at least one spanning tree, and any spanning tree contains exactly $n - 1$ branches.

In matrix terms, the spanning trees are in one-to-one correspondence with the column submatrices of $G$ formed from $n - 1$ columns and having maximum rank $n - 1$ (so that excluding any one row leaves a nonsingular submatrix).
The $n - 1$ selected columns correspond to the $n - 1$ branches in the spanning tree.
If the columns are selected from the matrix $G \Upsilon$ in place of $G$, the columns will be weighted by the respective branch admittances.
The network defined by $G$ is connected if and only if $G$ is of rank $n - 1$ (the maximum rank possible), in which case it also has at least one nonsingular submatrix with $n - 1$ rows and columns, the columns of which define a spanning tree.
By definition the rows of $G$ sum to zero (modulo 2 if signs are disregarded), from which it follows that the one row eliminated to form a nonsingular submatrix may be chosen arbitrarily.

Since the admittance matrix itself is given as $Y = G \Upsilon G^T$, it may readily be conjectured that $Y$ (assuming real positive admittances) attains its maximum rank $n - 1$ whenever the underlying graph is connected, and thus has at least one spanning tree ensuring $G$ has rank $n - 1$.
The Kirchhoff characteristic establishes this more precisely, with the help of the results obtained earlier on modified networks.

The Kirchhoff characteristic of a network $\nw$ is defined as the sum-of-products form
\begin{equation}
\kappa(\nw) = \sum_{\tr \subseteq \nw} \prod_{\alpha \in \tr} y_{\alpha}
\label{eq:kirchhoff}
\end{equation}
where the sum is over all spanning trees $\tr$ in $\nw$, and each tree in the sum contributes the product of all branch admittances making up the tree.
If $\nw$ has no spanning trees (is disconnected) then $\kappa(\nw) = 0$.
If $\nw$ is a single isolated node, one sets $\kappa(\nw) = 1$ by convention.

Now take any (connected) network $\nw$ on $n \geq 2$ nodes, and consider the following step-by-step procedure for obtaining $\kappa(\nw)$.
Since the network is connected, it has a spanning tree $\tr$.
Commence building the network with a single branch of $\tr$ and its two endpoint nodes, then add branches one at a time (the expansion operation discussed previously) utnil all $n$ nodes and the entirety of $\tr$ are obtained.
Then, add the remaining network branches (if any) in any desired order, using the augmentation operation discussed previously, until all branches are added and $\nw$ is fully constructed.
Observe that at each step of this process the network remains connected.

One may now derive recursive formulae for both $\kappa$ and the admittance matrix $Y$ at each step.
For the initial network with one branch, $\kappa = y_1$ is just the admittance of this branch, while for the admittance matrix one has
\begin{equation}
Y = \mat{\phantom{-}y_1 & -y_1 \\ -y_1 & \phantom{-}y_1}, \qquad
\text{with} \qquad c(Y) = \kappa = y_1.
\label{eq:twonode}
\end{equation}
In the process of expanding this to the tree $\tr$, the network obtained at each step is a tree, and so $\kappa$ is a product of all admittances added to the network thus far.
Meanwhile $Y$ expands by one row and column at each step, and by Proposition \ref{prop:expcof1}, $c(Y)$ also gets multiplied by each successive added admittance.
Once the last of the $n$ nodes is added, one has
\begin{equation}
\kappa(\tr) = c(Y_{\tr}) = \prod_{\alpha \in \tr} y_{\alpha}.
\label{eq:ktree}
\end{equation}

The remaining branches of $\nw$ are introduced by augmentation.
This alters $c(Y)$ at each step according to Proposition \ref{prop:augcof1} and formula (\ref{eq:augcof1}), which can also be expressed as
\begin{equation}
c(\bar{Y}) = c(Y) + y_+ \cdot c(Y^{[jk]})
\label{eq:augcof1c}
\end{equation}
where $y_+$ is the admittance added between nodes $j$ and $k$ to produce the augmented network with admittance matrix $\bar{Y}$, and $Y^{[jk]}$ is the admittance matrix of the network contracted on nodes $j$ and $k$, for which $c(Y^{[jk]}) = C_{jk,jk}(Y)$ by Proposition \ref{prop:concof1}.

Now consider what happens to the Kirchhoff characteristic when network $\nw$ is augmented to network $\bar{\nw}$ by adding a branch with admittance $y_+$ between nodes $j$ and $k$.
Clearly, any spanning tree in $\nw$ is also a spanning tree of $\bar{\nw}$, so $\kappa(\bar{\nw})$ will by definition be the sum of $\kappa(\nw)$ plus additional terms corresponding to new spanning trees in $\bar{\nw}$ that are not in $\nw$.
Any such `new' tree must contain the new branch with admittance $y_+$ (since otherwise it would belong to $\nw$).
Further, when this branch is removed it leaves two smaller trees, necessarily subtrees of $\nw$, one of which contains node $j$ but not node $k$, and the other of which contains node $k$ but not node $j$.
(Either of these may degenerate to a single node with no branches.)
It is not hard to see that any such pair of trees when `glued back together' uniquely defines a spanning tree of the network $\nw^{[jk]}$ obtained from $\nw$ by contracting nodes $j$ and $k$.
Conversely if one takes any spanning tree in $\nw^{[jk]}$, one may break this at node $j$ and insert the new branch with impedance $y_+$ to obtain a spanning tree of $\bar{\nw}$ that is not in $\nw$ (since it contains the new branch).
It follows that $\kappa(\nw^{[jk]})$ consists of a sum of terms that, when each is multiplied by $y_+$, gives a sum of terms corresponding one-to-one with the `new' trees of $\bar{\nw}$.
It has thus been established that
\begin{equation}
\kappa(\bar{\nw}) = \kappa(\nw) + y_+ \cdot \kappa(\nw^{[jk]}).
\label{eq:kaug}
\end{equation}
The isomorphic relationship between linear electrical networks $\nw$ and admittance matrices $Y$ was outlined in the Introduction.
Now it is seen that the common cofactor $c(Y)$ of the admittance matrix and the Kirchhoff characteristic $\kappa(\nw)$ of the underlying network satisfy the same formulae at each step in the network's construction: the same starting values by (\ref{eq:twonode}), the same values on trees by (\ref{eq:ktree}), and the same recursively defined formulae (\ref{eq:augcof1c}) and (\ref{eq:kaug}) under augmentation with additional branches.
Since any connected linear network may be recursively assembled in this manner, this furnishes an inductive proof of the equality of $c(Y)$ and $\kappa(\nw)$:
\begin{proposition}
\label{prop:kcof}
Let the electrical network $\nw$ have admittance matrix $Y$, and suppose the network is connected with at least two nodes.
Then the common value $c(Y)$ of all first cofactors of $Y$ is equal to the Kirchhoff characteristic of $\nw$ defined by (\ref{eq:kirchhoff}):
\begin{equation}
c(Y) = \kappa(\nw).
\label{eq:kcof}
\end{equation}
This equality remains true when $\nw$ is not connected, in which case $c(Y)$ and $\kappa(\nw)$ are both zero.
\end{proposition}
Note that for disconnected networks it was already shown in the previous section that $c(Y) = 0$, since then $Y$ is of rank no greater than $n - 2$.
$\kappa(\nw) = 0$ will also be true for such a network by definition, since $\nw$ has no spanning trees making $\kappa(\nw)$ an empty sum.

The following is an immediate corollary of Proposition \ref{prop:kcof}:
\begin{proposition}
\label{prop:fullrank}
Suppose the network with admittance matrix $Y$ is connected with $n \geq 2$ nodes, and all the admittances of the network are nonnegative real numbers, such that when all branches with zero admittance are deleted the network remains connected.
Then each first cofactor $c(Y)$ is a positive real number, ensuring in turn that $Y$ is of maximal rank $n - 1$.
More generally, if the network after deletion of zero-admittance branches has exactly $k$ separate connected components (with an isolated node counted as a component in its own right) then $Y$ is of rank at most $n - k$, and $c(Y) = 0$ if $k > 1$.
A sufficient condition for $Y$ to attain rank $n - k$ is that all remaining branch admittances are positive real numbers.
\end{proposition}
Proposition \ref{prop:fullrank} follows directly from the equality of $c(Y)$ and $\kappa(\nw)$.
When $\nw$ is connected (after omitting any zero-admittance branches) it has at least one spanning tree, and all spanning trees make strictly positive contributions to $\kappa(\nw)$ since each is a product of positive real numbers.
Thus $\kappa(\nw)$ is a sum of one or more positive real numbers, hence a positive real number itself.

The generalisation follows by considering the explicit partitioning of $Y$ as in (\ref{eq:sep}) in accordance with the $k$ connected components, from which it follows the maximum rank is $n - k$.
Applying Proposition \ref{prop:kcof} separately to each connected component (other than isolated nodes) establishes that positive admittances are again sufficient to obtain this rank.
Any isolated node adds a row and column of zeros to $Y$ and increases both $n$ and $k$ by 1, but does not affect the rank of $Y$.

A famous result in algebraic graph theory (originally due to Kirchhoff) also follows as a corollary of Proposition \ref{prop:kcof}.
\begin{proposition}[Matrix--Tree Theorem]
\label{prop:matrixtree}
Let $\Lambda = G G^T$ be the Laplacian of the graph with incidence matrix $G$.
Then any first cofactor of $\Lambda$ equals the number of spanning trees in the graph.
\end{proposition}
Indeed, let the electrical network $\nw$ be defined from the graph with incidence matrix $G$ by assigning admittance $y = 1$ to every branch.
Then the admittance matrix of $\nw$ is $Y = G G^T = \Lambda$, and every spanning tree contributes exactly 1 to the Kirchhoff characteristic $\kappa(\nw)$.
It immediately follows that $\kappa(\nw) = c(\Lambda)$ counts the number of spanning trees.

\subsection{Identities Involving Kirchhoff Characteristics}

A number of useful identities may be derived connecting the Kirchhoff characteristic $\kappa(\nw)$ to other properties of $\nw$ and related networks.

Firstly, the Kirchhoff characteristic provides a topological characterisation of the driving-point impedances $Z_{jk}$ between pairs of nodes in the network $\nw$.
By Proposition \ref{prop:concof1}, the principal cofactor $C_{jk,jk}(Y)$ in the definition of $Z_{jk}$ is equal to the first cofactors of $Y^{[jk]}$, the admittance matrix of the network $\nw^{[jk]}$ formed by contracting on nodes $j$ and $k$.
This is the same as the Kirchhoff characteristic $\kappa(\nw^{[jk]})$, by Proposition \ref{prop:kcof}.
It follows that
\begin{equation}
Z_{jk}(\nw) = \frac{C_{jk,jk}(Y)}{c(Y)} = \frac{c(Y^{[jk]})}{c(Y)} = \frac{\kappa(\nw^{[jk]})}{\kappa(\nw)}.
\label{eq:kzjkc}
\end{equation}
This quantity $\kappa(\nw^{[jk]})$ featured in the argument leading up to formula (\ref{eq:kaug}).
As noted there, the spanning trees of $\nw^{[jk]}$ are in one-to-one correspondence with pairs $(\tr_j, \tr_k)$ of disjoint trees in $\nw$, such that $\tr_j$ includes node $j$ but not $k$, $\tr_k$ includes $k$ but not $j$, and the two jointly span all nodes (though not all branches) in $\nw$.
The pairs include those where $\tr_j$ is an isolated node $j$, or $\tr_k$ an isolated node $k$.
The above may be summed up as the following result:
\begin{proposition}
\label{prop:kcon}
Let $\nw$ be a network on at least two nodes and $j$, $k$ any distinct nodes in $\nw$.
If $\nw^{[jk]}$ denotes the network formed from $\nw$ by contracting nodes $j$ and $k$, then the driving-point impedance $Z_{jk}$ in $\nw$ is given by the ratio (\ref{eq:kzjkc}) of Kirchhoff characteristics $Z_{jk}(\nw) = \kappa(\nw^{[jk]}) / \kappa(\nw)$.
Furthermore, $\kappa(\nw^{[jk]})$ is given by the formula
\begin{equation}
\kappa(\nw^{[jk]})
   = \sum_{\tr_j, \tr_k \subseteq \nw} \prod_{\alpha \in \tr_j} y_{\alpha} \prod_{\beta \in \tr_k} y_{\beta}
\label{eq:kzjk}
\end{equation}
where the sum is over all disjoint tree pairs $(\tr_j, \tr_k) \subseteq \nw$ such that $j \in \tr_j$, $k \in \tr_k$ and $\tr_j \cup \tr_k$ spans all nodes in $\nw$.
In (\ref{eq:kzjk}) an isolated node for $\tr_j$ or $\tr_k$ contributes the empty product 1.
\end{proposition}
Note that $\kappa(\nw^{[jk]})$ may alternatively be expressed as a restricted Kirchhoff characteristic $\kappa_{jk}(\bar{\nw}^{jk})$, where
\begin{itemise}
\item
$\bar{\nw}^{jk}$ refers to the network augmented with a fictitious unit admittance branch between nodes $j$ and $k$; and
\item
the restricted characteristic $\kappa_{jk}$ includes only those spanning trees of $\bar{\nw}^{jk}$ that contain the new branch between $j$ and $k$.
\end{itemise}
As a corollary to Proposition \ref{prop:kcon} one obtains the intuitive formula for the driving-point impedances $Z_{jk}(\nw)$ when $\nw$ itself is a tree.
In this case, all disjoint spanning tree-pairs are formed from $\nw$ by deleting one arbitrary branch, and nodes $j$ and $k$ belong to distinct trees of such a pair precisely when the removed branch lies on the unique path between $j$ and $k$.
Denoting this removed branch by $\gamma$, the contribution to $\kappa(\nw^{[jk]})$ in (\ref{eq:kzjk}) from the tree-pair with $\gamma$ removed is the product of all admittances except $y_{\gamma}$; observing that $\kappa(\nw) = \sum_{\alpha \in \nw} y_{\alpha}$, one further sees that the contribution to $Z_{jk}(\nw)$ by (\ref{eq:kzjkc}) is precisely $y_{\gamma}^{-1} = z_{\gamma}$.
Repeating for all branches on the path from $j$ to $k$, one has:
\begin{proposition}
\label{prop:ztree}
Let $\nw$ be a network whose nonzero admittances form a tree, and set $z_{\alpha} = y_{\alpha}^{-1}$ for each branch $\alpha$.
Then for all distinct nodes $j, k \in \nw$ one has
\begin{equation}
Z_{jk}(\nw) = \sum_{\alpha \in \path_{jk}} z_{\alpha}
\label{eq:ztree}
\end{equation}
where $\path_{jk}$ is the unique path in $\nw$ connecting nodes $j$ and $k$.
\end{proposition}
Another useful principle involving Kirchhoff characteristics is in essence a restatement of the formula (\ref{eq:kaug}) derived above for a network augmentation.
\begin{proposition}[Deletion--Contraction Principle]
\label{prop:kdelc}
Let $\nw$ be any connected network with $n \geq 2$ nodes and let $\alpha$ be any branch in $\nw$, with admittance $y_{\alpha}$.
Denote by $\nw - \alpha$ the subnetwork of $\nw$ (not necessarily connected) obtained by deleting the branch $\alpha$, and by $\nw \circ \alpha$ the network obtained by contracting the endpoints of $\alpha$.
Then
\begin{equation}
\kappa(\nw) = \kappa(\nw - \alpha) + y_{\alpha} \cdot \kappa(\nw \circ \alpha),
\label{eq:kdelc}
\end{equation}
where one recalls the convention that $\kappa(\nw - \alpha) = 0$ if $\nw - \alpha$ is disconnected ($\alpha$ is a `bridge' in $\nw$), and $\kappa(\nw \circ \alpha) = 1$ if $\nw \circ \alpha$ is a single isolated node.
\end{proposition}
The rule (\ref{eq:kdelc}) follows immediately from (\ref{eq:kaug}) in the case where $\nw - \alpha$ is connected: take the latter as the network $\nw$ in (\ref{eq:kaug}) and $\bar{\nw}$ as the network formed by adding back the branch $\alpha$ between nodes $j$ and $k$.
If $\nw - \alpha$ is disconnected, then $\kappa(\nw) = y_{\alpha} \cdot \kappa(\nw \circ \alpha)$ also follows by a similar argument, considering the spanning trees of the (necessarily connected) network $\nw \circ \alpha$.

Observing that by definition neither $\kappa(\nw - \alpha)$ nor $\kappa(\nw \circ \alpha)$ involves the admittance of the branch $\alpha$, the following is an immediate corollary:
\begin{proposition}
\label{prop:kderiv}
If the Kirchhoff characteristic $\kappa(\nw)$ is regarded as a function of the admittance $y_{\alpha}$ of a specified branch $\alpha \in \nw$, then
\begin{equation}
\frac{\partial}{\partial y_{\alpha}} \kappa(\nw) = \kappa(\nw \circ \alpha).
\label{eq:kderiv}
\end{equation}
\end{proposition}
By virtue of the construction of $\kappa(\nw)$, the derivative in (\ref{eq:kderiv}) can be understood analytically where appropriate, or else as a formal derivative of a multilinear sum-of-products form in $m$ indeterminates (where $m$ is the number of branches in $\nw$).

A useful identity that connects Kirchhoff characteristics of networks under contraction of two distinct node pairs can be derived from Proposition \ref{prop:concof2}.
\begin{proposition}
\label{prop:kcon2}
Let $\nw$ be a network on $n \geq 3$ nodes with admittance matrix $Y$, and let $j,k,p,q$ be nodes in $\nw$ such that the pairs $(j,k)$ and $(p,q)$ are distinct (have at most one node in common).
Let $\nw^{[jk]}$ and $\nw^{[pq]}$ denote networks formed from $\nw$ by contracting node pairs $(j,k)$ and $(p,q)$ respectively, and let $\nw^{[jk][pq]}$ denote the network formed from $\nw$ by contracting both pairs simultaneously.
Then
\begin{equation}
\kappa(\nw^{[jk]}) \cdot \kappa(\nw^{[pq]}) - \kappa(\nw) \cdot \kappa(\nw^{[jk][pq]}) = C_{pq,jk}(Y)^2,
\label{eq:kcon2}
\end{equation}
noting that $\kappa(\nw^{[jk][pq]}) = 1$ if $\nw^{[jk][pq]}$ is an isolated node.
In particular, when the network admittances are real numbers, the expression on the left hand side of (\ref{eq:kcon2}) is nonnegative, and is zero only when $C_{pq,jk}(Y) = 0$.
In the latter case, $(j,k)$ and $(p,q)$ are said to be \emph{conjugate pairs} in $\nw$.
\end{proposition}
Strictly speaking, Proposition \ref{prop:kcon2} also allows the possibility that $j = k$ or $p = q$ or both, in which case the respective contractions do not occur; but (\ref{eq:kcon2}) is then vacuous, with both sides immediately yielding zero values.
To show (\ref{eq:kcon2}) holds when $j \neq k$ and $p \neq q$, let $Y^{[jk]}$ and $Y^{[pq]}$ denote the admittance matrices of $\nw^{[jk]}$ and $\nw^{[pq]}$ respectively.
Both are obtained from $Y$ by contraction, and in particular
\begin{equation}
c(Y^{[jk]}) = C_{jk,jk}(Y), \qquad c(Y^{[pq]}) = C_{pq,pq}(Y).
\label{eq:ccon2}
\end{equation}
Now let $Y^{[jk][pq]}$ denote the admittance matrix of the doubly contracted network.
By assumption, node indices $p \neq q$ in $Y$ also correspond to distinct indices $p' \neq q'$ in the contracted admittance matrix $Y^{[jk]}$.
It follows that
\begin{equation}
c(Y^{[jk][pq]}) = C_{p'q',p'q'}(Y^{[jk]}).
\label{eq:cconpq}
\end{equation}
This is a second cofactor of an admittance matrix for a contracted network, and its value can therefore be written in terms of cofactors of $Y$ using Proposition \ref{prop:concof2}.
Indeed, using (\ref{eq:concof2}) and premultiplying by $c(Y)$ for convenience,
\begin{equation}
c(Y) \cdot C_{p'q',p'q'}(Y^{[jk]}) = C_{pq,pq}(Y) \cdot C_{jk,jk}(Y) - C_{pq,jk}(Y)^2.
\label{eq:cconpqpq}
\end{equation}
But this is no more or less than formula (\ref{eq:kcon2}) once the Kirchhoff characteristics are written in terms of equivalent cofactors:
\begin{multline}
\kappa(\nw^{[jk]}) \cdot \kappa(\nw^{[pq]}) - \kappa(\nw) \cdot \kappa(\nw^{[jk][pq]}) \\
   = C_{jk,jk}(Y) \cdot C_{pq,pq}(Y) - c(Y) \cdot C_{p'q',p'q'}(Y_{[jk]}) = C_{pq,jk}(Y)^2.
\label{eq:cconk}
\end{multline}
While Proposition \ref{prop:kcon2} links Kirchhoff characteristics of networks obtained by contraction of two distinct branches, a close variant on this links networks obtained by both deletion and contraction of branches.
It is readily derived by combining (\ref{eq:kcon2}) with two successive applications of the deletion-contraction principle (\ref{eq:kdelc}), and noting that deletion or contraction of branches between \emph{distinct} node pairs are commutative operations.
\begin{proposition}
\label{prop:kdelc2}
Let $\alpha$ and $\beta$ be two branches in a connected network $\nw$, connecting distinct node pairs $(j,k)$ and $(p,q)$.
Let $\nw - \alpha \beta$, $\nw \circ \alpha \beta$, $\nw \circ \alpha - \beta$ and $\nw \circ \beta - \alpha$ denote in an evident (and unique) manner the four networks formed from $\nw$ by either deletion or contraction of these two branches.
If $\nw$ has admittance matrix $Y$, then
\begin{equation}
\kappa(\nw \circ \alpha - \beta) \cdot \kappa(\nw \circ \beta - \alpha)
   - \kappa(\nw - \alpha \beta) \cdot \kappa(\nw \circ \alpha \beta) = C_{pq,jk}(Y)^2,
\label{eq:kdelc2}
\end{equation}
keeping in mind the Kirchhoff characteristic of a disconnected network is zero.
In particular, when the network admittances are real numbers, (\ref{eq:kdelc2}) is positive unless $(j,k)$ and $(p,q)$ are conjugate pairs, in which case it is zero.
\end{proposition}

\section{Monotonicity Properties of Passive DC Networks}
\label{sec:monopos}

This and the following section use previous results to derive certain well-known---and some lesser known---general properties of electrical networks with `realistic' (passive) admittances.
These are characterised by algebraic or differential current-voltage relationships with positive coefficients, or in AC networks, a well-defined increasing dependence on AC frequency $\omega$.
When admittances are prescribed in this manner thay impart certain monotonicity and positivity properties to characteristic functions of the network, additional to those described in previous sections.

In direct-current (DC) networks, where voltages and currents are constant in the steady state and transient behaviour is disregarded, the network is characterised by ordinary Ohm's law relationships between branch currents and voltage differences.
Accordingly, the branch admittances $y_{\alpha}$ in a DC network are assumed to be nonnegative real numbers.
Any admittance may be set to zero to omit that branch from the network, provided that doing so does not cause $\nw$ to become disconnected.
Henceforth, a network is called \emph{nonzero-connected} when it is connected and remains connected after deleting all zero-admittance branches.
If a branch admittance $y_{\alpha}$ is nonzero, then the branch impedance $z_{\alpha} = 1 / y_{\alpha}$ is also defined and is a positive real number.

\subsection{Orientability of DC networks, and the impedance metric}
\label{sec:dcmetric}

An important property of flows in DC networks is expressed informally in the principle that `power always flows downhill'.
As a consequence of Ohm's law with positive admittances, the direction of current and power flow is always from nodes at higher voltage to those at lower voltage.

From this observation one may also infer a a useful extremal principle: that when a DC network $\nw$ is energised from a single source connected between nodes $p$ and $q$, the voltages $v_p$ and $v_q$ establish bounds for all other node voltages in $\nw$.
\begin{proposition}
\label{prop:dcextrema}
Let $\nw$ be a nonzero-connected network whose admittances are nonnegative real numbers.
Let $p$ and $q$ be any two nodes in $\nw$ and let a DC source with current $I > 0$ be connected from node $q$ to node $p$.
Then $v_p$ (respectively $v_q$) is maximal (respectively minimal) among the node voltages in $\nw$: that is, $v_p \geq v \geq v_q$ for all node voltages $v$.
Consequently, the difference $v_p - v_q$ is also maximal among all node-pair voltage differences in $\nw$.
Further, for any node $k \not\in \{p,q\}$ the voltage $v_k$ equals $v_p$ (respectively $v_q$) if and only if every path of nonzero admittances in $\nw$ from $k$ to $q$ (respectively $p$) includes node $p$ (respectively $q$); equivalently, if in every spanning tree $\tr \subseteq \nw$ of nonzero admittances, the deletion of node $p$ (respectively $q$) separates $k$ from the other of $\{p, q\}$.
\end{proposition}
The underlying idea in this proposition is quite intuitive and indeed often taken for granted.
Together with the basic `power flows downhill' principle it expresses the idea that voltages in a DC network are \emph{orientable} based on the source connections.

Nonetheless, Proposition \ref{prop:dcextrema} is not quite trivial to establish.
To do so, let the source $I$ be connected as stated, and assign each branch in $\nw$ having nonzero current flow a direction corresponding to that of the current.
Observe that this direction is necessarily from a node at higher voltage to one at strictly lower voltage, as a consequence of Ohm's law.
Now as $I$ is the only source, it follows immediately from Kirchhoff's current law that any node other than $p$ or $q$ with an outgoing directed branch must also have an incoming directed branch, and vice versa.
Suppose there is a node $k \neq p$ whose voltage $v_k$ is maximal among nodes in $\nw$.
Since $\nw$ is nonzero-connected, node $k$ has a positive-admittance branch linking it to some other node in $\nw$, which by hypothesis has voltage no greater than $v_k$.
But nor can this voltage be less than $v_k$: otherwise, Ohm's law forces this to be an outgoing directed branch from $k \neq p$, which presupposes an incoming directed branch from some other node $j$ with $v_j > v_k$, contradicting the hypothesis that $v_k$ is maximal.

It follows that any node $k \neq p$ with maximal node voltage $v_k$ is connected by a positive-admittance branch to another node at the same voltage $v_k$---and not to any node at higher or lower voltage.
If this node is not $p$, one may repeat the argument as many times as necessary, exploiting the nonzero-connectedness of $\nw$ to identify successive nodes in $\nw$ not previously visited, ultimately exhausting all node possibilities $k \neq p$.
One concludes there is a path connecting node $k$ and node $p$, on which all admittances are nonzero yet all node voltages take the same value $v_k$ (and hence the current flow must be zero).
Moreover, there is no other path of nonzero admittances linking $k$ to any other node at different voltage, including $q$ (since $I > 0$ hence $v_p > v_q$ by nonzero-connectivity).
This establishes first that $v_p$ is also maximal, and secondly that $p$ lies on every path of nonzero admittances between $k$ and $q$, as desired.

If there is no node $k \neq p$ with maximal voltage as above, then $v_p$ itself is maximal, and in this case $p$ is the unique node of maximal voltage.

An analogous argument to the above establishes that any node $k \neq q$ of minimal voltage in $\nw$ is connected to $q$ by a path of nonzero admittances with all nodes at identical voltage, hence $v_q$ is also minimal and $q$ is on all paths of nonzero admittances from $k$ to $p$.
If there is no such node $k \neq q$, then $q$ is the unique node with minimal voltage $v_q$.
This completes the proof of Proposition \ref{prop:dcextrema}.

Proposition \ref{prop:dcextrema}, in turn, helps to establish an important property of driving-point impedances $Z_{jk}$ in a (nonzero-connected) DC network: that they satisfy all the properties of a \emph{metric} on the set of network nodes.
As such, $Z_{jk}$ is termed the \emph{resistance distance} between nodes $j$ and $k$ in theoretical chemistry literature \cite{kr:rd,k:rdsr}.
$Z_{jk}$ is indeed strictly positive for all node pairs $(j,k)$ as a consequence of nonzero-connectedness, and the basic properties of cofactors of the network admittance matrix $Y$ ensure that $Z_{jk} = Z_{kj}$ and $Z_{jj} = 0$ for all $j$ and $k$.
This leaves the \emph{triangle inequality} $Z_{pq} + Z_{qr} \geq Z_{pr}$ for triples $(p,q,r)$ as the only metric property remaining to be demonstrated, and this is now shown to be a consequence of Proposition \ref{prop:dcextrema}.

Consider three alternative connections of the same source $I > 0$ to the same network $\nw$: (a) from node $q$ to node $p$; (b) from node $r$ to node $q$; and (c) from node $r$ to node $p$.
Given any node index $k$, let $v_k^a$, $v_k^b$ and $v_k^c$ denote, respectively, the voltage at node $k$ in an arbitrary chosen solution in each case.
(These voltages are only determined up to a constant specific to each case, but as only voltage differences are used below, this will not be important.)
By definition of $Z_{jk}$, one has
\begin{equation}
v_p^a - v_q^a = Z_{pq} I, \qquad v_q^b - v_r^b = Z_{qr} I, \qquad v_p^c - v_r^c = Z_{pr} I.
\label{eq:dcvpqr}
\end{equation}
Observe now that case (c) can be obtained as the sum of cases (a) and (b), since the effect of the source in (c) is exactly that of the sources in (a) and (b) taken in combination.
Therefore by linearity,
\begin{equation}
v_p^c - v_r^c = (v_p^a - v_r^a) + (v_p^b - v_r^b).
\label{eq:dcvpr}
\end{equation}
Now one may invoke Proposition \ref{prop:dcextrema} to bound the voltage differences on the right of (\ref{eq:dcvpr}).
Since there is only one source in each case, one has
\begin{equation}
v_p^a - v_r^a \leq v_p^a - v_q^a \qquad \text{and} \qquad v_p^b - v_r^b \leq v_q^b - v_r^b
\label{eq:dcvbnd}
\end{equation}
with equality in both cases if any only if all paths of nonzero admittances between nodes $p$ and $r$ go via node $q$.
Using the bounds (\ref{eq:dcvbnd}) in (\ref{eq:dcvpr}), substituting (\ref{eq:dcvpqr}) and cancelling the common (positive) factor $I$ leads to $Z_{pr} \leq Z_{pq} + Z_{qr}$, as required.
One has the following:
\begin{proposition}
\label{prop:dcmetric}
Let $\nw$ be a nonzero-connected network whose admittances are nonnegative real numbers.
Then the driving-point impedances $Z_{jk}$, where $j$ and $k$ are any two node indices, define a metric on the set of nodes in $\nw$, satisfying the usual properties of a metric as follows:
\begin{enumerate}
\item
$Z_{jj} = 0$ and $Z_{jk} > 0$ for all $j$ and all $k \neq j$.
\item
$Z_{jk} = Z_{kj}$ for all $j$ and $k$.
\item
$Z_{pq} + Z_{qr} \geq Z_{pr}$ for all $p$, $q$ and $r$, with equality if any only if every path of nonzero admittances in $\nw$ between nodes $p$ and $r$ includes node $q$ ($q$ is a `bridge node' relative to $p$ and $r$).
\end{enumerate}
\end{proposition}
The behaviour of this impedance metric is particularly simple when $\nw$ is a tree: in this case, by Proposition \ref{prop:ztree}, $Z_{jk}$ for any node indices $j$ and $k$ is found just by summing the (real, positive) impedances $z_{\alpha}$ along the unique path connecting nodes $j$ and $k$.
The triangle inequality becomes an equality, $Z_{pq} + Z_{qr} = Z_{pr}$, if and only if node $q$ is on the unique path between nodes $p$ and $r$.

Using the triangle inequality in formula (\ref{eq:transimp}) for transfer impedances, after setting $k = q$ in that formula, leads to the following immediate corollary:
\begin{proposition}
\label{prop:dctransimp}
Let $\nw$ be a nonzero-connected network whose admittances are nonnegative real numbers.
Then all transfer impedances $\tz{pq}{rq}$ with identical second indices are nonnegative.
Further, $\tz{pq}{rq} = 0$ if and only if $q$ is a bridge node relative to $p$ and $r$; consequently, two node pairs with a node in common are conjugate if and only if the common node is a bridge relative to the others.
\end{proposition}
An interesting consequence of this result is that in the formula (\ref{eq:gsoln}), $\hat{v} = T \hat{\imath}$ for a grounded solution of a DC network, the matrix $T$ has all nonnegative entries.
Further, the zero entries in $T$ (if any) are exactly those for which the ground node itself is a bridge relative to the corresponding indexed nodes, and in this case it will immediately follow that $T$ separates (perhaps with a suitable relabelling) into a form like (\ref{eq:sep}) based on the components into which $\nw$ is separated by removal of the ground node.
Carrying this through, it is seen that $T$ can always be rearranged into a block diagonal form, in which each block has all \emph{positive} entries.

\subsection{Monotonicity of driving-point impedances in DC networks}
\label{sec:dcmono}

Since in a DC network both the branch admittances and the driving-point impedances are real numbers, with the latter constituting a metric---a well-defined `electrical distance' between any two nodes---it is natural to consider how these driving-point impedances vary as the individual branch admittances take greater or lesser values.
Accordingly let $j$, $k$ be any two nodes in the connected network $\nw$ and $\alpha$ any branch, and consider the derivative of $Z_{jk}(\nw)$ by (\ref{eq:kzjkc}) with respect to the branch admittance $y_{\alpha}$.

One may deal first with the special case where branch $\alpha$ is connected between nodes $j$ and $k$.
In this case, branch $\alpha$ vanishes in the contracted network $\nw^{[jk]} = \nw \circ \alpha$, and
\begin{equation}
\frac{\partial}{\partial y_{\alpha}} Z_{jk}(\nw)
   = \kappa(\nw \circ \alpha) \cdot \frac{\partial}{\partial y_{\alpha}} \paren{\frac{1}{\kappa(\nw)}}
   = - \paren{\frac{\kappa(\nw \circ \alpha)}{\kappa(\nw)}}^2 = - Z_{jk}(\nw)^2 < 0.
\label{eq:zderiv0}
\end{equation}
Note that the quantity on the right hand side will never be zero, since by hypothesis $\nw$, hence also $\nw \circ \alpha$, are connected with positive real admittances after disregarding zero-admittance branches.
(In case $\nw$ has only two nodes, one may again take $\kappa(\nw \circ \alpha) = 1$.)

Suppose now that branch $\alpha$ is not connected between nodes $j$ and $k$.
Then by (\ref{eq:kderiv}) and the quotient rule,
\begin{equation}
\frac{\partial}{\partial y_{\alpha}} Z_{jk}(\nw) = \frac{\kappa(\nw) \cdot \kappa(\nw^{[jk]} \circ \alpha)
   - \kappa(\nw^{[jk]}) \cdot \kappa(\nw \circ \alpha)}{\kappa(\nw)^2}.
\label{eq:zderiv}
\end{equation}
Here one may invoke Proposition \ref{prop:kcon2}, letting $p$ and $q$ denote the endpoint nodes of branch $\alpha$ so that $\nw \circ \alpha = \nw^{[pq]}$, $\nw^{[jk]} \circ \alpha = \nw^{[jk][pq]}$, and pairs $(j,k)$ and $(p,q)$ are distinct by assumption:
\begin{equation}
\frac{\partial}{\partial y_{\alpha}} Z_{jk}(\nw) = \frac{\kappa(\nw) \cdot \kappa(\nw^{[jk][pq]})
   - \kappa(\nw^{[jk]}) \cdot \kappa(\nw^{[pq]})}{\kappa(\nw)^2}
   = - \paren{\frac{C_{pq,jk}(Y)}{c(Y)}}^2 = - \tz{pq}{jk}(Y)^2.
\label{eq:zderivpq}
\end{equation}
Formula (\ref{eq:zderivpq}), in conjunction with (\ref{eq:zderiv0}), leads to the following general result, credited to Rayleigh:
\begin{proposition}[Rayleigh's Monotonicity Theorem]
\label{prop:rayleigh}
Let $\nw$ be a nonzero-connected network whose admittances are nonnegative real numbers, and let $Y$ be the corresponding admittance matrix.
Then for any node pair $(j,k)$, the driving-point impedance $Z_{jk}(\nw)$ is positive, and is monotonically nonincreasing as a function of any individual network admittance.
More specifically, if the admittance $y_{\alpha}$ is connected between nodes $p$ and $q$ then
\begin{equation}
\frac{\partial}{\partial y_{\alpha}} Z_{jk}(\nw) = - \tz{pq}{jk}^2(Y) \leq 0 \qquad
\paren{(p,q) \text{ endpoints of } y_{\alpha}}
\label{eq:rayleigh}
\end{equation}
where $\tz{pq}{jk}$ is the transfer impedance defined by (\ref{eq:transz}).
If $(j,k)$ and $(p,q)$ are conjugate pairs in $\nw$, then $C_{pq,jk}(Y) = \tz{pq}{jk}(Y) = 0$ and $Z_{jk}$ is insensitive to the value of $y_{\alpha}$.
In all other cases, including where the unordered pairs $(j,k)$ and $(p,q)$ coincide, $Z_{jk}$ strictly decreases with $y_{\alpha}$ at a rate given by $\tz{pq}{jk}^2$ at any value $y_{\alpha} > 0$ (extending to $y_{\alpha} = 0$ if $\nw - \alpha$ is connected).
\end{proposition}
Note that if $Z_{jk}$ is decreasing (respectively, constant) with respect to an admittance $y_{\alpha}$, this implies it is increasing (respectively, constant) with respect to the corresponding impedance $z_{\alpha} = 1 / y_{\alpha}$.

In practice, conjugacy implies that under the action of a source connected between nodes $p$ and $q$ (and no other sources), the nodes $j$ and $k$ are always driven to the same potential---and vice versa.
The canonical example is the Wheatstone bridge configuration, based on four admittances $y_{\alpha}$, $y_{\beta}$, $y_{\gamma}$, $y_{\delta}$ arranged in a square.
A source is presumed connected between two diagonally opposite nodes (labelled 1 and 2 below) and voltage sensed across the `bridge' between the other two diagonally opposite nodes (labelled 3 and 4).
For completeness, one may also assume a `test' admittance $y_{\tau}$ connected across the bridge nodes 3 and 4, and a source admittance $y_{\sigma}$ between nodes 1 and 2.
These may not be present in a given instance, in which case the admittance is set to zero.

The admittance matrix for the full Wheatstone bridge configuration is
\begin{equation}
Y = \mat{y_{\alpha} + y_{\gamma} + y_{\sigma} & -y_{\sigma} & -y_{\alpha} & -y_{\gamma} \\
   -y_{\sigma} & y_{\beta} + y_{\delta} + y_{\sigma} & -y_{\beta} & -y_{\delta} \\
   -y_{\alpha} & -y_{\beta} & y_{\alpha} + y_{\beta} + y_{\tau} & -y_{\tau} \\
   -y_{\gamma} & -y_{\delta} & -y_{\tau} & y_{\gamma} + y_{\delta} + y_{\tau}}.
\label{eq:whsy}
\end{equation}
The Kirchhoff characteristic can be found by computing any first cofactor of $Y$, but is more easily written down by inspection of the underlying network graph, which is the net of a tetrahedron.
The spanning trees are obtained by enumerating all $\binom{6}{3} = 20$ combinations of three branches and excluding those corresponding to the 4 triangular faces, for a total of 16 products of three admittances.
It is convenient to write this in the following form, which accounts easily for the removal of $y_{\sigma}$ or $y_{\tau}$:
\begin{multline}
c(Y) = \kappa(\nw) = y_{\alpha} y_{\beta} y_{\gamma} + y_{\alpha} y_{\beta} y_{\delta}
   + y_{\alpha} y_{\gamma} y_{\delta} + y_{\beta} y_{\gamma} y_{\delta} \\
   + y_{\sigma} \paren{y_{\alpha} y_{\gamma} + y_{\alpha} y_{\delta}
      + y_{\beta} y_{\gamma} + y_{\beta} y_{\delta}}
   + y_{\tau} \paren{y_{\alpha} y_{\beta} + y_{\alpha} y_{\delta}
      + y_{\beta} y_{\gamma} + y_{\gamma} y_{\delta}} \\
   + y_{\sigma} y_{\tau} \paren{y_{\alpha} + y_{\beta} + y_{\gamma} + y_{\delta}}.
\label{eq:whsk}
\end{multline}
One may also easily obtain expressions for the driving-point impedance $Z_{12}$ at the source, as well as the impedance $Z_{34}$ across the bridge; these are
\begin{align}
Z_{12} &= \frac{1}{\kappa(\nw)} \bbrak{\paren{y_{\alpha} + y_{\beta}} \paren{y_{\gamma} + y_{\delta}}
   + y_{\tau} \paren{y_{\alpha} + y_{\beta} + y_{\gamma} + y_{\delta}}},
\label{eq:whsz12} \\
Z_{34} &= \frac{1}{\kappa(\nw)} \bbrak{\paren{y_{\alpha} + y_{\gamma}} \paren{y_{\beta} + y_{\delta}}
   + y_{\sigma} \paren{y_{\alpha} + y_{\beta} + y_{\gamma} + y_{\delta}}}.
\label{eq:whsz34}
\end{align}
(Note that as $Y$ in this case is a $4 \times 4$ matrix, one must take care to avoid the easy mistake of constructing the relevant $2 \times 2$ minors by including rather than excluding the nominated indices.)

The key Wheatstone bridge property stems from the transfer impedance $\tz{12}{34}$.
From the matrix $Y$ in (\ref{eq:whsy}) it is not hard to see that the relevant cofactor is
\begin{equation}
C_{12,34}(Y) = y_{\alpha} y_{\delta} - y_{\beta} y_{\gamma}.
\label{eq:whsc}
\end{equation}
Setting this to zero leads directly to the classical `balance' criterion for the bridge, $y_{\alpha} / y_{\beta} = y_{\gamma} / y_{\delta}$.
Under this condition, $(1,2)$ and $(3,4)$ are a conjugate pair, from which it follows that:
\begin{itemise}
\item
Under the action of the source between nodes 1 and 2, nodes 3 and 4 are always at the same potential (by definition of $\tz{12}{34}$);
\item
Under the action of the same source, the current in $y_{\tau}$ is zero (by Ohm's law);
\item
The driving-point impedance $Z_{12}$ at the source is insensitive to $y_{\tau}$; and
\item
The driving-point impedance $Z_{34}$ at the bridge is insensitive to $y_{\sigma}$.
\end{itemise}
The insensitivity may be confirmed by direct calculation of partial derivatives, but is also revealed by analysis of the expressions for $\kappa(\nw)$, $Z_{12}$ and $Z_{34}$.
Let the bridge ratio be $y_{\alpha} / y_{\beta} = y_{\gamma} / y_{\delta} = 1 / k$, so that $y_{\beta} = k y_{\alpha}$ and $y_{\delta} = k y_{\gamma}$.
Crucially, the Kirchhoff characteristic then factors into two terms, one depending on $y_{\tau}$ but not $y_{\sigma}$, the other depending on $y_{\sigma}$ but not $y_{\tau}$:
\begin{equation}
\begin{split}
\kappa(\nw) &= k \paren{1 + k} y_{\alpha} y_{\gamma} \paren{y_{\alpha} + y_{\gamma}}
      + \paren{1+k}^2 y_{\sigma} y_{\alpha} y_{\gamma} + k y_{\tau} \paren{y_{\alpha} + y_{\gamma}}^2
      + \paren{1 + k} y_{\sigma} y_{\tau} \paren{y_{\alpha} + y_{\gamma}} \\
   &= \bparen{\paren{1 + k} y_{\alpha} y_{\gamma} + y_{\tau} \paren{y_{\alpha} + y_{\gamma}}}
      \bparen{k \paren{y_{\alpha} + y_{\gamma}} + \paren{1 + k} y_{\sigma}}.
\end{split}
\label{eq:whskb}
\end{equation}
In the formula (\ref{eq:whsz12}) for $Z_{12}$, the terms involving $y_{\tau}$ then cancel, leaving an expression involving admittances other than $y_{\tau}$ only:
\begin{equation}
Z_{12} = \frac{1 + k}{\kappa(\nw)} \bbrak{\paren{1 + k} y_{\alpha} y_{\gamma}
      + y_{\tau} \paren{y_{\alpha} + y_{\gamma}}}
   = \paren{\frac{k}{1 + k} \paren{y_{\alpha} + y_{\gamma}} + y_{\sigma}}^{-1}.
\label{eq:whsz12b}
\end{equation}
Likewise, in formula (\ref{eq:whsz34}) for $Z_{34}$ the terms involving $y_{\sigma}$ cancel:
\begin{equation}
Z_{34} = \frac{y_{\alpha} + y_{\gamma}}{\kappa(\nw)} \bbrak{k \paren{y_{\alpha} + y_{\gamma}}
      + \paren{1 + k} y_{\sigma}}
   = \paren{\paren{1 + k} \frac{y_{\alpha} y_{\gamma}}{y_{\alpha} + y_{\gamma}} + y_{\tau}}^{-1}.
\label{eq:whsz34b}
\end{equation}
These expressions for $Z_{12}$ and $Z_{34}$ have an intuitive form.
$Z_{12}$ is the impedance of a parallel combination of $y_{\sigma}$ with a scaled parallel combination of $y_{\alpha}$ and $y_{\gamma}$---as would be the case if $y_{\tau}$ were omitted.
Likewise $Z_{34}$ is the impedance of a parallel combination of $y_{\tau}$ with a scaled series combination of $y_{\alpha}$ and $y_{\gamma}$, as when the source and $y_{\sigma}$ are not present.

The general sensitivity formula (\ref{eq:rayleigh}) may also be confirmed by direct calculation.
For a simple example, let $y_{\sigma} = y_{\tau} = 0$ in the above network, so that
\begin{equation}
Z_{12} = \frac{\paren{y_{\alpha} + y_{\beta}} \paren{y_{\gamma} + y_{\delta}}}
   {y_{\alpha} y_{\beta} y_{\gamma} + y_{\alpha} y_{\beta} y_{\delta}
      + y_{\alpha} y_{\gamma} y_{\delta} + y_{\beta} y_{\gamma} y_{\delta}}.
\label{eq:whsz120}
\end{equation}
Consider the sensitivity to the admittance $y_{\alpha}$ (connected between nodes 1 and 3).
By the quotient rule, one has
\begin{align}
\frac{\partial Z_{12}}{\partial y_{\alpha}} &= \frac{\paren{y_{\gamma} + y_{\delta}}
      \paren{y_{\alpha} y_{\beta} y_{\gamma} + y_{\alpha} y_{\beta} y_{\delta}
         + y_{\alpha} y_{\gamma} y_{\delta} + y_{\beta} y_{\gamma} y_{\delta}}
   - \paren{y_{\alpha} + y_{\beta}} \paren{y_{\gamma} + y_{\delta}}
      \paren{y_{\beta} y_{\gamma} + y_{\beta} y_{\delta} + y_{\gamma} y_{\delta}}}
   {\paren{y_{\alpha} y_{\beta} y_{\gamma} + y_{\alpha} y_{\beta} y_{\delta}
      + y_{\alpha} y_{\gamma} y_{\delta} + y_{\beta} y_{\gamma} y_{\delta}}^2} \nonumber \\
   &= - \paren{\frac{y_{\beta} \paren{y_{\gamma} + y_{\delta}}}
   {y_{\alpha} y_{\beta} y_{\gamma} + y_{\alpha} y_{\beta} y_{\delta}
      + y_{\alpha} y_{\gamma} y_{\delta} + y_{\beta} y_{\gamma} y_{\delta}}}^2.
\label{eq:whsdz120}
\end{align}
Direct inspection of $Y$ in (\ref{eq:whsy}) confirms that $C_{13,12}(Y) = y_{\beta} (y_{\gamma} + y_{\delta})$ when $y_{\tau} = 0$, confirming (\ref{eq:rayleigh}) in this instance.

Lastly, Proposition \ref{prop:kdelc2} may also be illustrated using the network example above.
Consider applying (\ref{eq:kdelc2}) to the branches $\sigma$ and $\tau$ of the Wheatstone bridge, with $(jk) = (12)$ and $(pq) = (34)$.
The networks $\nw - \sigma \tau$, $\nw \circ \sigma \tau$, $\nw \circ \sigma - \tau$ and $\nw \circ \tau - \sigma$ are all configurations of the four branches $\alpha$, $\beta$, $\gamma$ and $\delta$, with $\nw \circ \sigma \tau$ on two nodes, $\nw \circ \sigma - \tau$ and $\nw \circ \tau - \sigma$ on three nodes and $\nw - \sigma \tau$ a four-node square.
So by inspection,
\begin{equation}
\begin{split}
\kappa(\nw \circ \sigma - \tau) &= \paren{y_{\alpha} + y_{\beta}} \paren{y_{\gamma} + y_{\delta}}, \\
\kappa(\nw \circ \tau - \sigma) &= \paren{y_{\alpha} + y_{\gamma}} \paren{y_{\beta} + y_{\delta}}, \\
\kappa(\nw - \sigma \tau) &= y_{\alpha} y_{\beta} y_{\gamma} + y_{\alpha} y_{\beta} y_{\delta}
   + y_{\alpha} y_{\gamma} y_{\delta} + y_{\beta} y_{\gamma} y_{\delta}, \\
\kappa(\nw \circ \sigma \tau) &= y_{\alpha} + y_{\beta} + y_{\gamma} + y_{\delta},
\end{split}
\label{eq:whsdelc}
\end{equation}
and a laborious but straightforward calculation shows that
\begin{multline}
\kappa(\nw \circ \sigma - \tau) \cdot \kappa(\nw \circ \tau - \sigma)
      - \kappa(\nw - \sigma \tau) \cdot \kappa(\nw \circ \sigma \tau)
   = y_{\alpha}^2 y_{\delta}^2 + y_{\beta}^2 y_{\gamma}^2
      - 2 y_{\alpha} y_{\beta} y_{\gamma} y_{\delta} \\
   = \paren{y_{\alpha} y_{\delta} - y_{\beta} y_{\gamma}}^2 = C_{12,34}(Y)^2.
\end{multline}

\section{Properties of Passive AC Networks}
\label{sec:acprop}

\subsection{Admittance, impedance and reactance functions}
\label{sec:acfuncs}

When one moves away from DC steady state operation of networks, the important practical consideration is that the voltage-current behaviour of elements is not merely static (depending on instantaneous magnitudes alone) but also dynamic (depending on the time evolution of these magnitudes).
In place of the pure Ohmic relation $i_{\alpha} = y_{\alpha} v_{\alpha} = y_{\alpha} \paren{v_{\alpha+} - v_{\alpha-}}$, a capacitor (with capacitance $C$) or an inductor (with inductance $L$) obey the respective \emph{differential} relations
\begin{equation}
i_{\alpha} = C \frac{dv_{\alpha}}{dt}, \qquad v_{\alpha} = L \frac{di_{\alpha}}{dt}.
\label{eq:lc}
\end{equation}
In either case, the `admittance' connecting the current and the voltage difference across the element is a differential operator.
Following Heaviside, one might informally write $p = d/dt$ in the above equations to obtain the equivalent admittance as an operator, a function of $p$.
More rigorously, one may take Laplace transforms of the voltages and currents, so that in place of time functions they become functions of a complex-valued variable $s$, and the operation of taking time derivatives becomes formally equivalent to multiplication by $s$.
This formal substitution $d/dt \mapsto s$ leads to
\begin{equation}
i_{\alpha}(s) = s C v_{\alpha}(s), \qquad v_{\alpha}(s) = s L i_{\alpha}(s),
\label{eq:lcs}
\end{equation}
leading to the admittances $y_{\alpha} = s C$ and $y_{\alpha} = 1 / s L$ respectively.
The admittance matrix in this formalism includes entries containing $s$ and $s^{-1}$, leading to cofactors (hence transfer and driving-point impedances) that are rational functions of $s$.

An important special case is an AC network, where sources are assumed to be sinusoidal functions of time with a fixed frequency $\omega$.
This permits a different, more straightforward mathematical simplification (anticipated by Maxwell and Heaviside, but first fully elaborated by Steinmetz) based on the fact that the algebra of fixed-frequency sinusoidal functions---linear combinations of $\sin \omega t$ and $\cos \omega t$ for fixed $\omega$---is closed not only under addition, subtraction and scalar multiplication but also under time derivatives (and antiderivatives, if additive constants are disregarded).

This sinusoidal algebra is built on a two-dimensional vector space, a general element of which is a function $\sqrt{2} r \cos\paren{\omega t + \phi}$ parameterised by the (root mean square) amplitude $r$ and phase displacement $\phi$.
An equivalent representation is as $\sqrt{2} \paren{a \cos \omega t - b \sin \omega t}$, where $a = r \cos\phi$, $b = r \sin\phi$ and $a^2 + b^2 = r^2$.
The zero element of the vector space has $r = a = b = 0$ and $\phi$ is indeterminate.
(These conventions are not entirely unique, and those used by others in defining $r$ and $\phi$ may differ in trivial ways from the above, for example by reversing the sign of $\phi$.)

There is a natural isomorphism between this vector space and the complex numbers $\cplx$: one simply identifies $r$ and $\phi$ with the magnitude and phase respectively of the complex \emph{phasor} $z$ (a portmanteau of `phase vector').
The coefficients $a$ and $b$ are then recognised as the real and imaginary parts of $z = r \e^{j \phi} = a + jb$.
Addition of sinusoidal functions maps naturally to addition of complex numbers, and vice versa, as does multiplication by real scalars, and the zero element of $\cplx$ has zero magnitude and indeterminate phase.
Differentiation maps sinusoidal functions to sinusoidal functions, and in the isomorphism with $\cplx$ corresponds to multiplication by $j \omega$:
\begin{multline}
\frac{d}{dt} \paren{\sqrt{2} a \cos \omega t - \sqrt{2} b \sin \omega t}
   = \sqrt{2} (- \omega b) \cos \omega t - \sqrt{2} (\omega a) \sin \omega t \\ \mapsto \qquad
j \omega \paren{a + jb} = (- \omega b) + j (\omega a).
\label{eq:jw}
\end{multline}
Thus, as long as one deals exclusively with sinusoidal functions, represented by complex numbers with frequency implied by context, one may specialise the formal substitution $d/dt \mapsto s$ above to $d/dt \mapsto j\omega$, where $\omega$ is interpreted as frequency.
Likewise, the Laplace transform formula itself under the substitution $s = j\omega$ takes on the appearance of a Fourier transform, as appropriate for the analysis of sinusoidal signals, although care is required in making this mathematically precise.

(Note that some abuse of notation will also be evident in the sequel, where $j$ is used both as the imaginary unit in $\cplx$ and as a node index.
The meaning should always be clear from context.)

For AC circuits, then, the characteristic formulae (\ref{eq:lcs}) for capacitors and inductors become
\begin{equation}
i_{\alpha} = j \omega C v_{\alpha}, \qquad v_{\alpha} = j \omega L i_{\alpha},
\label{eq:lcw}
\end{equation}
with equivalent admittances $y_{\alpha} = j \omega C$ and $- j / \omega L$, respectively.
The most general form of a branch admittance in a passive, linear AC network is $y_{\alpha} = G + jB$ with $G$ and $B$ real numbers, which allows for arbitrary combinations of resistive, capacitive and inductive characteristics.
The component $B$ is denoted the \emph{susceptance}, and $G$ the (parallel) conductance.
From above, the susceptance of a capacitor and inductor are $\omega C$ and $-1 / \omega L$ respectively.

As with DC admittances, whenever $y_{\alpha} \neq 0$ one may also define the branch impedance $z_{\alpha} = 1 / y_{\alpha}$.
This has the general form $z_{\alpha} = R + jX$, where $R$ denotes the (series) resistance and the imaginary component $X$ the \emph{reactance}.
The components $G$, $B$, $R$ and $X$ satisfy the elementary identities
\begin{equation}
R = \frac{G}{G^2 + B^2}, \quad X = - \frac{B}{G^2 + B^2}, \qquad \text{and} \qquad
G = \frac{R}{R^2 + X^2}, \quad B = - \frac{X}{R^2 + X^2},
\label{eq:gbrx}
\end{equation}
and as a consequence, the following are true for any AC impedance or admittance quantity:
\begin{itemise}
\item
When $G = 0$ then $R = 0$, and vice versa (unless $R$ is undefined).
Similarly, $B = 0$ implies $X = 0$ and vice versa (unless $X$ is undefined).
\item
When $G > 0$ then $R > 0$, irrespective of the value of $B$.
Likewise, $R > 0$ implies $G > 0$ irrespective of the value of $X$.
\item
When $B \neq 0$ then $X \neq 0$ and $X$ has the opposite sign to $B$, irrespective of the value of $G$.
LIkewise, $X \neq 0$ implies $B \neq 0$ with the opposite sign, irrespective of the value of $R$.
\item
When $G = R = 0$, then $X = - 1 / B$ and $B = - 1 / X$, provided both are nonzero.
\end{itemise}
It is clear in particular that the reactance of a capacitor and an inductor are, respectively, $- 1 / \omega C$ and $\omega L$.
Observe that whether one is considering susceptance $B$ or reactance $X$, each quantity for a single element varies either linearly with frequency $\omega$ or else as $- 1 / \omega$---both of which are strictly increasing functions of $\omega$ (although the latter tends to infinity as $\omega \goesto 0$).
It will be seen that this strict increase with $\omega$ of both susceptance and reactance extends more generally to AC networks made up of capacitors and inductors.

It is readily apparent on the other hand that due to the possibility of combining capacitive and inductive elements in an AC circuit, with susceptances or reactances of opposite signs, there is no equivalent of Proposition \ref{prop:dcextrema} for voltage magnitudes; that is to say, voltages are not orientable in AC networks.
As an elementary example, consider a three-node LC network with $y_{12} = - j / \omega L$, $y_{23} = j \omega C$, and a source $I$ connected between nodes 1 and 3.
The voltage developed across the source is $v_1 - v_3 = Z_{13} I = j \paren{\omega L - 1 / \omega C} I$, while that across the capacitor alone is $v_1 - v_2 = j \omega L I$ and with $\omega L > 0$ and $\omega C > 0$, is greater in magnitude than the voltage across the source.

(At the same time, the presence of reactances with opposite signs is not essential here.
The following section demonstrates that orientability of AC voltages, and the triangle inequality for AC impedances, can fail in general even when all reactances have the same sign.)

To gain further insight on impedance functions in LC and RLC networks it is helpful to revert to the general case of functions of the Laplace variable $s$: results specific to AC networks then follow upon substituting $s = j\omega$.
In the more general formalism based on $s$, the assumed form for an admittance or impedance element needs further development.
It will be sufficient to assume the linear fractional (M\"{o}bius) form with nonnegative real coefficients
\begin{equation}
y_{\alpha}(s) = \frac{g + sc}{r + sl}, \qquad z_{\alpha}(s) = \frac{r + sl}{g + sc}, \qquad
g, c, r, l \geq 0, \quad g + c \neq 0, \quad r + l \neq 0.
\label{eq:gcrl}
\end{equation}
This form is sufficiently general to represent a resistor, capacitor or inductor standing alone, or either of the latter in parallel or series combination with a resistor.
Without sacrificing this generality, it will also be assumed \emph{either} that $lc = 0$ \emph{or} that $gl - rc \neq 0$ (or both): this sidesteps unimportant technical issues with pole-zero cancellation.

For an AC steady state at frequency $\omega$, the functions (\ref{eq:gcrl}) again revert to the phasor forms $y_{\alpha} = G + jB$, $z_{\alpha} = R + jX$, where
\begin{equation}
G = \frac{rg + \omega^2 lc}{r^2 + \omega^2 l^2}, \qquad
B = \omega \frac{rc - gl}{r^2 + \omega^2 l^2}, \qquad
R = \frac{rg + \omega^2 lc}{g^2 + \omega^2 c^2}, \qquad
X = \omega \frac{gl - rc}{g^2 + \omega^2 c^2}.
\label{eq:gcrlac}
\end{equation}
The function $y_{\alpha}(s)$ or $z_{\alpha}(s)$ is termed a \emph{reactance function} if $G = 0$ (equivalently $R = 0$) under the substitution $s = j\omega$.
Assuming $y_{\alpha}$ and $z_{\alpha}$ are not identically zero this implies that either $l, g = 0$ and $c, r > 0$ (capacitive case) or that $c, r = 0$ and $l, g > 0$ (inductive case).
In either case, one readily establishes that both $B$ and $X$ are strictly increasing functions of $\omega$.
However, this strict increase is no longer guaranteed if $y_{\alpha}$ is not a reactance function.
The case $c = 0$, $l, r > 0$ in (\ref{eq:gcrl}) already furnishes a counterexample, in which $B(\omega)$ decreases to a minimum at $\omega^* = r/l$ and is only an increasing function for $\omega > r/l$ (though it is negative for all $\omega > 0$).
Electrically, this represents the inductive susceptance of a $R$--$L$ series filter.

The total AC admittance $|y_{\alpha}|$ and impedance $|z_{\alpha}|$ corresponding to $y_{\alpha}(j\omega)$ are
\begin{equation}
|y_{\alpha}| = \sqrt{G^2 + B^2} = \sqrt{\frac{g^2 + \omega^2 c^2}{r^2 + \omega^2 l^2}}, \qquad
|z_{\alpha}| = \sqrt{R^2 + X^2} = |y_{\alpha}|^{-1}
   = \sqrt{\frac{g^2 + \omega^2 c^2}{r^2 + \omega^2 l^2}},
\label{eq:gcrlyz}
\end{equation}
and one readily verifies not only that $G$, $B$, $R$, $X$ in (\ref{eq:gcrlac}) satisfy (\ref{eq:gbrx}), but also that
\begin{equation}
\frac{\partial |y_{\alpha}|}{\partial \omega}
   = \frac{\omega}{|y_{\alpha}|} \frac{r^2 c^2 - g^2 l^2}{\paren{r^2 + \omega^2 l^2}^2}, \qquad
\frac{\partial |z_{\alpha}|}{\partial \omega}
   = \frac{\omega}{|z_{\alpha}|} \frac{g^2 l^2 - r^2 c^2}{\paren{g^2 + \omega^2 c^2}^2}
   = - \frac{1}{|y_{\alpha}|^2} \frac{\partial |y_{\alpha}|}{\partial \omega}.
\label{eq:gcrlyzw}
\end{equation}
This shows that the strictly increasing property of susceptance or reactance with $\omega$ does not extend to total admittance or impedance, even for reactance functions.
Instead, for the functions (\ref{eq:gcrl}), $|y_{\alpha}|$ is a strictly increasing function of $\omega$ when $y_{\alpha}$ is `more capacitive' (in the sense that $rc > gl$), and a decreasing function of $\omega$ when `more inductive' ($gl > rc$).
The case $rc = gl$, for which $lc = 0$ by assumption, reduces to that of a pure resistance ($l = c = 0$).

Returning to (\ref{eq:gcrl}) and looking beyond pure AC systems, the Laplace variable $s$ takes in general any value $s = \sigma + j\omega \in \cplx$.
In this more general sense, $s$ takes on the role of a derivative operator in an alternative two-dimensional space of functions represented by complex numbers.
These are analogous to the pure sinusoidal functions with frequency $\omega$, but now incorporate an exponential time-weighting with the specific rate constant $\sigma$:
\begin{equation}
r \e^{j \phi} \mapsto f_{\sigma,\omega}(r, \phi) = \sqrt{2} r \e^{\sigma t} \cos(\omega t + \phi)
\quad \text{or} \quad
a + jb \mapsto g_{\sigma,\omega}(a, b) = \sqrt{2} \e^{\sigma t} \paren{a \cos\omega t - b \sin\omega t}.
\label{eq:trigexp}
\end{equation}
One recovers the space of pure sinusoidal functions at frequency $\omega$ by setting $\sigma = 0$.
If instead one sets $\omega = 0$, the variable $s = \sigma$ is the derivative operator for the one-parameter family of pure exponential time functions $\sqrt{2} a \e^{\sigma t}$, represented by real numbers $a$ (which are mapped to real numbers under the action of the operator $s$).

Where $s$ includes a nonzero real part $\sigma$, this leads to more general versions of (\ref{eq:gcrlac}) for $G$, $B$, $R$ and $X$, quantifying the relationship between voltage and current when either takes the more general form (\ref{eq:trigexp}) as a function of time.
In a purely algebraic sense, $g$ in the formulae analogous to (\ref{eq:gcrlac}) is replaced with $g + \sigma c$ and $r$ with $r + \sigma l$ when $\sigma \neq 0$.
This leads to the formulae
\begin{equation}
G(\sigma, \omega) = \frac{r g + \sigma \paren{rc + gl} + \paren{\sigma^2 + \omega^2} lc}
   {r^2 + 2 \sigma r l + \paren{\sigma^2 + \omega^2} l^2}, \qquad
B(\sigma, \omega) = \omega \frac{r c - g l}{r^2 + 2 \sigma r l + \paren{\sigma^2 + \omega^2} l^2},
\label{eq:gbs}
\end{equation}
wth similar formulae for $R(\sigma,\omega)$ and $X(\sigma,\omega)$.
Owing to the assumed nonnegativity of the parameters $g$, $c$, $r$, $l$, one finds in this more general case, as in (\ref{eq:gcrlac}), that $G(\sigma,\omega) \geq 0$ and $R(\sigma,\omega) \geq 0$ if $\sigma \geq 0$.

\subsection{Complex power and (non-)orientability in AC networks}
\label{sec:acpower}

The non-orientability in general of AC networks was noted in the previous section.
Further examples are provided here indicating particular cases where the orientability of DC networks can be extended to the AC case, and also of where this fails.
Groundwork will also be laid toward a more limited form of orientability of phase angles, described in the next section, which can be ascribed to a special subset of AC networks that occur commonly in practice.

First of all, if an AC network is purely resistive it is straightforward to verify that Propositions \ref{prop:dcextrema}, \ref{prop:dcmetric} and \ref{prop:dctransimp} apply just as in a DC network.
With the admittance matrix $Y$ real-valued, a source current $I \in \reals$ leads to orientable network voltages $v_k \in \reals$ satisfying Proposition \ref{prop:dcextrema}, with these real values now representing the amplitudes of AC quantities that are aligned in phase.
Indeed, even in the general case where $I = |I| \e^{j \phi} \in \cplx$ is arbitrary, the equivalence of the network equations
\begin{equation}
Y v = \paren{e_p - e_q} |I| \e^{j \phi} \qquad \text{and} \qquad
Y \paren{v \e^{-j \phi}} = \paren{e_p - e_q} |I|
\label{eq:acreal}
\end{equation}
(where $p$ and $q$ denote the source terminals) implies that for every node index $k$, the quantity $v_k \e^{-j \phi}$ is real-valued, and that these quantities together form an orientable set satisfying Proposition \ref{prop:dcextrema}, with maximum at node $p$ and minimum at node $q$.
The actual voltages $v_k$ have these real orientable values as amplitudes while aligning in phase with the source current $I$, all having phase angle $\phi$.
This phase alignment also carries over to the voltage differences and (by Ohm's law) to the current flows within the network.
Accordingly, in a resistive AC network as in a DC network, the direction of positive (active) power flow is from points of higher voltage amplitude to those of lower amplitude.

There is one important caveat of course: in DC networks it is natural for voltages to take either positive or negative magnitudes, but in AC networks (as in complex analysis) it is more convenient to regard all magnitudes as positive, and to re-interpret a `negative' magnitude as a corresponding unsigned \emph{amplitude} with a phase shift by $\pi$.
To avoid technical ambiguities raised by negative magnitudes, it is natural in AC network analysis to identify the negative source terminal $q$ as a ground node and to work with solutions grounded at $q$.
Since by Proposition \ref{prop:dcextrema} the voltage $v_q = 0$ is minimal among the values $v_k \e^{-j \phi}$, this is sufficient to ensure all other such values are nonnegative, thus yielding proper unambiguous amplitudes for AC voltages all having phase angle $\phi$.

(This is the case for a single source.
When multiple sources are present in an AC network, they will not in general all have their negative terminal at the same node $q$; however, it has been seen previously that by linearity one may notionally replace any arbitrary current source with two sources each having a terminal at $q$, one identical to the original and the other rotated by $\pi$.
When considering the effect of multiple sources, it is of course the full complex-valued voltages and currents that combine and not the real magnitudes, unless the source currents happen to all align in phase or differ in phase by $\pi$.
The reasoning leading to the proof of Proposition \ref{prop:dcmetric}, for AC networks with orientable (signed) voltage magnitudes as above, relies only on sources aligned in phase in this manner.
With multiple AC sources of arbitrary phase, one cannot rely on voltages being orientable.)

Observe now that the refactorisation as in (\ref{eq:acreal}) above applies not just when $Y$ is real-valued, but also if it happens that all branch admittances $y_{\alpha}$ themselves have identical phase angles $\theta$, meaning a common term $\e^{j \theta}$ can be factored out of $Y$ and absorbed, along with $\e^{j \phi}$, into the voltages $v$.
The following proposition, built on the argument above, describes the general (non-resistive) case.
\begin{proposition}
\label{prop:acmetric}
Let $\nw$ be a nonzero-connected network whose admittances are complex numbers $y_{\alpha} = g_{\alpha} + j b_{\alpha}$ with $g_{\alpha} \geq 0$ such that
\begin{enumerate}
\item
all branch susceptances $b_{\alpha}$ are nonzero and of the same sign; and
\item
the ratio $g_{\alpha} / b_{\alpha} (= \rho)$ is the same for all branches $\alpha$ (possibly zero).
\end{enumerate}
For any two nodes $(j,k)$ let the (complex-valued) driving-point impedance between $j$ and $k$ be denoted $Z_{jk} = R_{jk} + jX_{jk}$.
Then a metric on the set of nodes, in the sense of Proposition \ref{prop:dcmetric}, is defined by the impedance magnitudes $|Z_{jk}|$ or alternatively by the unsigned reactances $|X_{jk}|$ for node pairs $(j,k)$.
Further, if a single AC source with current $|I| \e^{j \phi}$ is connected from node $q$ to node $p$, then all the voltages in $\nw$ are aligned in phase with phase angle $\phi - \theta$, where $\theta = \arccot\rho \in [-\pi/2, \pi/2]$ (and for $\rho = 0$ the sign of $\theta$ is taken identical to that of the $b_{\alpha}$).
\end{proposition}
That the values $|Z_{jk}|$ define a metric follows immediately by factoring $\e^{j \theta}$ out of the elements of the admittance matrix and applying Proposition \ref{prop:dcmetric} to the equivalent resistive network.
That the $X_{jk}$ likewise define a metric is assured by the condition on the $b_{\alpha}$ values, which ensures that $X_{jk}$ and $|Z_{jk}|$ differ only by the multiplicative constant $\sin\theta \neq 0$.
While the conditions of Proposition \ref{prop:acmetric} are stated in terms of admittance components, they have straightforward equivalent statements in terms of the reactances $x_{\alpha}$ and ratios $r_{\alpha} / x_{\alpha} = -\rho$ of the branches having nonzero admittance.

Proposition \ref{prop:acmetric} is unfortunately not very useful in practice.
In general, neither branch admittances nor voltages in an AC network will align in phase, and the phase angles $\delta_k = \arg{v_k}$ will vary widely.
This leads to a breakdown of both the orientability principle and the metric property in general meshed AC networks.
This will be shown below after developing further the idea of power flow in AC networks and its connection to complex voltages.

Voltage phase angles $\delta_k$ in AC networks have a useful relation to power flows, stemming from a basic identity for complex numbers: if $u, v \in \cplx$, then the product $z = u v^*$ has magnitude $|z| = |u| |v|$---the product of the original magnitudes---and argument $\arg z = \delta_{uv} = \arg u - \arg v$, the difference in phase angles between $u$ and $v$ (equivalently, the amount by which $v$ `lags' $u$ in phase).
In particular, if $u$ and $v$ have the same phase angle then $z$ is real.
More generally, if $v$ is analysed as the vector sum of two components, one `in phase' with $u$ and having the same argument, and the other `in quadrature' with an argument 90 degrees out of phase, then $\re{z}$ can be identified as the product of $|u|$ with the `in phase' component of $v$, while $\im{z}$ is the product of $|u|$ with the `quadrature' component.
(Equivalently, $\re{u v^*}$ defines an inner product on the complex plane, vanishing precisely when $u$ and $v$ are orthogonal as Euclidean vectors.)

A consequence of the above is that quantities of the form $u v^*$ are invariant under a global `reference frame shift' that advances or retards the phases of all network voltages and currents by the same angle---as suggested by (\ref{eq:acreal}) above.
It is a basic physical fact that power flows in networks are invariant under this kind of shift in reference frame.
Indeed if the principle is applied to a node voltage $v_k$ and a branch current $i_{\alpha}$ entering or leaving this node, the product $s(v_k, i_{\alpha}) = v_k i_{\alpha}^*$ is denoted the \emph{complex power} associated with the current $i_{\alpha}$ at the voltage $v_k$.
It has the form $s = p + j q$, where the \emph{active power} $p$ is derived from the component of $i_{\alpha}$ in phase with $v_k$, and the \emph{reactive power} $q$ from the component of $i_{\alpha}$ in quadrature.
In electrical terms $p$ is the mean rate of `real' (instantaneous) energy flow on branch $\alpha$ where it meets node $k$.
Note it is important to identify the specific endpoint $k$ where power is measured.

It also follows that Tellegen's theorem, identity (\ref{eq:tellegen}), is more conveniently restated for AC networks in terms of complex-power products $v_{\alpha}^* i_{\alpha}$ or $v_{\alpha} i_{\alpha}^*$, rather than $v_{\alpha} i_{\alpha}$.
Fortunately, this requires only a trivial modification of the argument from Kirchhoff's laws to the identity (\ref{eq:tellegen}) using matrix $\gfull$ and vectors $\vfull$ and $\ifull$.
As $\gfull$ is real-valued, it satisfies $\gfull^* = \gfull^T$, and since $\vfull = \gfull^T v = \gfull^* v$, one has
\begin{equation}
\sum_{\alpha} v_{\alpha}^* i_{\alpha} = \vfull^* \ifull = v^* \gfull \ifull = 0 \qquad \text{and} \qquad
\sum_{\alpha} v_{\alpha} i_{\alpha}^* = \ifull^* \vfull = \ifull^* \gfull^* v = (\gfull \ifull)^* v = 0.
\label{eq:tellegenc}
\end{equation}
Tellegen's theorem in full generality is thus the following:
\begin{proposition}[Tellegen]
\label{prop:tellegen}
Suppose that in a network $\nw$ with arbitrary sources, the voltages and currents are all either real-valued or complex-valued.
Then one has
\begin{equation}
\sum_{\alpha} v_{\alpha} i_{\alpha} = \sum_{\alpha} v_{\alpha}^* i_{\alpha}
   = \sum_{\alpha} v_{\alpha} i_{\alpha}^* = 0,
\label{eq:tellegenrc}
\end{equation}
where the sum is taken over all branches $\alpha$ of $\nw$ \emph{and} all source branches, $v_{\alpha}$ is the voltage difference across the network or source branch $\alpha$, and $i_{\alpha}$ is the branch or source current.
\end{proposition}

As noted above in connection with negative magnitudes, there is an important technical point to bear in mind when interpreting power flow quantities in AC networks.
By Proposition \ref{prop:unique} and Proposition \ref{prop:grounded}, the equations $Y v = i$ for any network define not a unique solution but rather a family of solutions, differing by the addition of a common value to all network voltages $v_k$.
In AC networks, this equates to a uniform translation of all voltages in the complex plane.
Any such `level shift' in voltages preserves the value of all voltage \emph{differences} $v_j - v_k$, hence also the value of the power quantities $v_{\alpha} i_{\alpha}^*$ appearing in Tellegen's theorem, which represent the complex power $p_{\alpha} + j q_{\alpha}$ generated or consumed in a single branch $\alpha$.
However, such a shift obviously does not preserve the value of a voltage $v_k$ \emph{in isolation}, nor the magnitude $|v_k|$ or phase angle $\delta_k$, nor even the difference of two such magnitudes.

This becomes important when associating the power \emph{flow} $v_k i_{\alpha}^* = p + j q$ with a current $i_{\alpha}$ leaving or entering a node at voltage $v_k$ as above.
Defining $p + j q$ as a `flow' implicitly assumes a network solution grounded at some node $g$, with $p + j q$ being the power notionally consumed were branch $\alpha$ connected from node $k$ directly to node $g$, or notionally generated were branch $\alpha$ connected from node $g$ to node $k$, as appropriate.
So a power flow of the form $v_k i_{\alpha}^*$ with a \emph{single} node voltage $v_k$ is only ever defined relative to a \emph{particular} grounded network solution, after identifying the grounded node $g$.
Likewise, magnitudes or phase angles of individual node voltages are also implicitly defined relative to this particular solution grounded at node $g$.

(The same holds true in DC networks or any network in general; DC networks however provide for a closer identification between power flow and current flow as each has a clear and consistent direction.
In AC networks by contrast the direction of power flow can only be inferred from a branch current $i_{\alpha}$ by describing its relationship to an explicitly quantified endpoint node voltage.)

Because a single node voltage is involved in defining power flows, the following is obtained as a direct corollary to Kirchhoff's current law, generalising it from currents to power flows.
\begin{proposition}
\label{prop:kclpq}
Let $k$ be any node in an AC network $\nw$ with nonzero voltage $v_k$.
Then the algebraic sum of complex power flows on branches incident with node $k$ is zero.
That is, one has $\sum_{\alpha} p_{\alpha} = \sum_{\alpha} q_{\alpha} = 0$, where the sum is over all branches $\alpha$ incident with node $k$, the power flows are measured at the node $k$ end of each branch, and power entering $k$ given the opposite sign to power leaving $k$.
\end{proposition}

Returning again to complex power in single branches---that is, to quantities $s_{\alpha} = v_{\alpha} i_{\alpha}^*$ where $v_{\alpha} = v_j - v_k$ is the voltage `drop' on a branch from $j$ to $k$---it is evident that, for a branch with admittance $y_{\alpha}$:
\begin{equation}
s_{\alpha} = v_{\alpha} i_{\alpha}^* = v_{\alpha} \paren{y_{\alpha} v_{\alpha}}^*
   = y_{\alpha}^* |v_{\alpha}|^2
\label{eq:sbrch}
\end{equation}
so that $s_{\alpha}$ has a phase angle equal and opposite to that of the admittance $y_{\alpha}$ (and equal to that of the branch impedance $z_{\alpha} = y_{\alpha}^{-1}$).
In particular, the reactive power consumption $q$ on a branch is zero for resistive branches, positive for inductive branches and negative for capacitive branches.
The identity
\begin{equation}
s_{\alpha} = v_j i_{\alpha}^* - v_k i_{\alpha}^*
\label{eq:sbrchjk}
\end{equation}
where $\alpha$ is directed from node $j$ to node $k$, equates the branch's own complex power to the amount by which the power flow decreases from the `sending' to the `receiving' end of the branch.

One might readily imagine that identity (\ref{eq:sbrchjk}) restores a form of the `power flows downhill' principle if, for all branches $\alpha$, the real and imaginary parts of $s_{\alpha}$ can separately be assumed to have a common sign.
When for example the network contains resistors and inductors only, both $p_{\alpha}$ and $q_{\alpha}$ are nonnegative for all branches, so power flow differences as in (\ref{eq:sbrchjk}) are also uniformly nonnegative.

In these situations some stronger conclusions regarding the driving-point impedances $Z_{jk}$ are indeed possible.
Let a single source with arbitrary (nonzero) current $I$ be directed from node $k$ to node $j$.
By Proposition \ref{prop:transz} the voltage developed across the source is $V_{jk} = Z_{jk} \cdot I$.
Now apply Tellegen's theorem (\ref{eq:tellegenrc}), noting that the current $I$ (directed from the `negative' to the `positive' terminal) has the opposite sign convention to the usual branch current $i_{\alpha}$ (directed from `positive' to `negative'), leading to the complex-power identity
\begin{equation}
V_{jk} I^* = \sum_{\alpha} v_{\alpha} i_{\alpha}^* = \sum_{\alpha} s_{\alpha}
\label{eq:tellegens}
\end{equation}
where the sum on the right is over all network branches.
Then substituting $V_{jk} = Z_{jk} I$ and $s_{\alpha}$ as in (\ref{eq:sbrch}) leads to
\begin{equation}
V_{jk} I^* = Z_{jk} |I|^2 = \sum_{\alpha} y_{\alpha}^* |v_{\alpha}|^2,
\label{eq:tellegenz}
\end{equation}
or with $I$ nonzero,
\begin{equation}
Z_{jk} = \sum_{\alpha} \left|\frac{v_{\alpha}}{I}\right|^2 y_{\alpha}^*
\label{eq:zjky}
\end{equation}
which expresses $Z_{jk}$ as a sum of conjugated branch admittances $y_{\alpha}^*$ weighted by nonnegative real numbers.
Since \emph{cones} in the complex plane (sectors $\{z : \arg z \in I\}$ for intervals $\ivl$ of width at most $\pi$) are closed under addition and positive scaling, one has the following:
\begin{proposition}
\label{prop:zcone}
Let $\ivl$ denote a closed, open, or semi-open interval of phase angles, of length $|\ivl| \leq \pi$, and $-\ivl$ its negation.
(For example, if $\ivl = (\theta_1, \theta_2]$ then $-\ivl = [-\theta_2, -\theta_1)$.)
Let $\nw$ be a network with admittance matrix $Y$, all of whose admittances $y_{\alpha}$ have $\arg y_{\alpha} \in \ivl$.
Then all driving-point impedances $Z_{jk}$ in $\nw$ have phase angles confined to the negated interval:
\begin{equation}
\arg Z_{jk}(Y) \in -\ivl, \qquad j, k \in \nw.
\label{eq:zcone}
\end{equation}
Equivalently, the impedances $Z_{jk}(Y)$ are confined in phase to the same sector as branch impedances $z_{\alpha} = y_{\alpha}^{-1}$.
\end{proposition}
Suppose for example that all admittances in $\nw$ are strictly inductive: $y_{\alpha} = g_{\alpha} - j b_{\alpha}$ with $g_{\alpha} \geq 0$ and $b_{\alpha} > 0$ for all branches $\alpha$.
Then $\arg y_{\alpha} \in [-\pi, 0)$ for all $\alpha$ and Proposition \ref{prop:zcone} asserts that $\arg Z_{jk} \in (0, \pi]$ for all pairs $(j,k)$.
If a source $I$ is connected with its positive terminal at node $p$ and its negative terminal at some node $q$ (taken as ground for the time being), then the voltage $v_p = Z_{jk} I$ will lead the current $I$ in phase and there will be a strictly positive reactive power flow $q = \im{v_p I^*}$ into the network at node $p$.

Intuitively, the idea suggests itself that in such a network, the positive flows of reactive power induced by the source would lead to a component of voltage being orientable, by virtue of the identity (\ref{eq:sbrchjk}).
This works, for example, when $\nw$ is a tree.
Then there is a unique path $\path_{pq}$ from node $p$ to node $q$, and the current is $I$ in all branches on $\path_{pq}$.
Suppose the source current $I$ has phase angle $\phi$, then (\ref{eq:sbrchjk}) on any branch in $\path_{pq}$ gives
\begin{equation}
y_{\alpha}^* |v_{\alpha}|^2 = \paren{v_{\alpha+} - v_{\alpha-}} |I| \e^{-j \phi}, \qquad \text{or} \qquad
v_{\alpha+} - v_{\alpha-} = \frac{|v_{\alpha}|^2}{|I|} \e^{j \phi} \paren{g_{\alpha} + j b_{\alpha}}
\label{eq:vtreeq}
\end{equation}
and hence
\begin{equation}
\im{\e^{-j \phi} v_{\alpha+} - \e^{-j \phi} v_{\alpha-}} > 0 \qquad \text{for all } \alpha \in \path_{pq}.
\label{eq:vtree}
\end{equation}
This suffices to establish that $\im{\e^{-j\phi} v}$ is an orientable component of voltage across the network $\nw$---in the sense of Proposition \ref{prop:dcextrema} and its proof---relative to sources with current phase angle $\phi$.
Again admitting the possibility that $\im{\e^{-j\phi} v}$ may take negative values for the purpose of argument, the reasoning for DC networks may now be followed to establish that the reactances $X_{jk} = \im{\e^{-j\phi} (Z_{jk} I)} / |I|$ define a metric on the nodes of $\nw$.

Notwithstanding the above, there is a simpler method for establishing metric properties of impedances in tree networks, the above example being a special case of the following:
\begin{proposition}
\label{prop:treemetric}
Let $\nw$ be an AC network that is a tree, and suppose $\re{\e^{j\theta} y_{\alpha}} > 0$ for all branch admittances $y_{\alpha} \in \nw$ and some $\theta \in \reals$.
(For example, a strictly inductive network satisfies this property with $\theta = \pi/2$.)
Then the quantities $d_{jk} = \re{\e^{-j\theta} Z_{jk}}$, with $(j,k)$ nodes in $\nw$ and $Z_{jk}$ the driving point impedance, define a metric on the nodes of $\nw$.
The triangle inequality $d_{pq} + d_{qr} \geq d_{pr}$ is an equality if and only if node $q$ is on the unique path $\path_{pr}$ between nodes $p$ and $r$.
\end{proposition}
To establish Proposition \ref{prop:treemetric} it suffices to observe that by Proposition \ref{prop:ztree} and inspection of the paths $\path_{pq}$, $\path_{qr}$ and $\path_{pr}$ in a tree network (and noting the $y_{\alpha}$ are nonzero by hypothesis), one has
\begin{equation}
Z_{pq} + Z_{qr} - Z_{pr} = 2 \sum_{\alpha \in \path_{pq \cap qr}} y_{\alpha}^{-1},
\qquad \text{hence} \qquad
d_{pq} + d_{qr} - d_{pr} = 2 \sum_{\alpha \in \path_{pq \cap qr}} \re{\paren{\e^{j\theta} y_{\alpha}}^{-1}},
\label{eq:treez}
\end{equation}
where $\path_{pq \cap qr}$ is shorthand for the intersection $\path_{pq} \cap \path_{qr}$.
If this intersection is empty, the right hand side is zero; this occurs precisely when $q$ is on $\path_{pr}$.
Otherwise, the fact that $\re{z} > 0$ implies $\re{z^{-1}} > 0$ for all $z \in \cplx$ establishes that the right hand side is strictly positive.

Proposition \ref{prop:treemetric} unfortunately lacks any straightforward generalisation to more general (meshed) networks $\nw$.
Observe, however, that the argument leading up to (\ref{eq:vtree}) primarily depends not on $\nw$ being a tree, so much as the branch currents on $\path_{pq}$ all having the same phase angle $\phi$.
Relaxing this assumption on the currents allows the same reasoning to be applied in more general networks.
Returning then to the formulae (\ref{eq:sbrch}) and (\ref{eq:sbrchjk}), the general formula that results for the voltage difference across a branch (analogous to (\ref{eq:vtreeq}) for trees) is
\begin{equation}
v_{\alpha+} - v_{\alpha-} = |v_{\alpha}|^2 \paren{y_{\alpha} i_{\alpha}^{-1}}^* \qquad (i_{\alpha} \neq 0).
\label{eq:vbrchq}
\end{equation}
The key to generalising (\ref{eq:vtree})---thus obtaining an orientable component of voltage in the more general network---is to establish that for all branches $\alpha$ with nonzero current (hence nonzero voltage difference), $v_{\alpha+} - v_{\alpha-}$ lies within a cone in $\cplx$ (in the sense of Proposition \ref{prop:zcone}) defined by an \emph{open} interval $(\theta - \pi/2, \theta + \pi/2)$ of phase angles, for some suitable value of $\theta$.
It will then follow that
\begin{equation}
\re{\e^{-j\theta} \paren{v_{\alpha+} - v_{\alpha-}}} > 0 \qquad
   \text{for all } \alpha \in \nw \text{ with } i_{\alpha} \neq 0,
\label{eq:vcone}
\end{equation}
and consequently that $\re{\e^{-j\theta} v}$ is an orientable component of voltage across $\nw$.

The previous argument for trees, leading to (\ref{eq:vtree}), amounts to doing this for a cone with $\theta = \phi + \pi/2$, relying on all currents having phase angle $\phi$ and all admittances having $\arg y_{\alpha} \in [-\pi/2, 0)$.
Indeed, observe that on this assumption one has by (\ref{eq:vbrchq})
\begin{equation}
\arg y_{\alpha} i_{\alpha}^{-1} \in \left[-\pi/2 - \phi ,\: -\phi\right), \qquad \text{hence} \qquad
\arg (v_{\alpha+} - v_{\alpha-}) \in \left(\phi ,\: \phi + \pi/2\right].
\label{eq:vbrchqtree}
\end{equation}
Notice that while the interval of phase angles in (\ref{eq:vbrchqtree}) is only semi-open, its width is $\pi/2$, well under the maximum width $\pi$ for a cone.
This provides substantial headroom for the phase angles of $i_{\alpha}$ to deviate from $\phi$ while preserving the key `cone' property for the voltage drops.
The most general result along these lines is the following:
\begin{proposition}
\label{prop:conemetric}
Let $\ivl = [-\phi_1, \phi_2)$ be a semi-open interval covering zero, with $\phi_1$, $\phi_2 \geq 0$ and $\phi_1 + \phi_2 < \pi$.
(The trivial case $\phi_1 = \phi_2 = 0$, $\ivl = \{0\}$ is also admitted.)
An AC network $\nw$ with admittances $y_{\alpha} \in \cplx$ is called \emph{dendromorphic} with \emph{interval of spread} $\ivl$ if, whenever a single current source $I$ is connected between an arbitrary node pair $(j,k) \in \nw$, the branches can be oriented such that all branch currents are either zero or satisfy $\arg i_{\alpha} \in [\phi - \phi_1, \phi + \phi_2)$ where $\phi$ is the phase angle of $I$.
(By phase invariance, it is sufficient to take $\phi = 0$, $\arg i_{\alpha} \in \ivl$ without loss of generality.)

Suppose $\nw$ is dendromorphic with interval of spread $\ivl$, and that all branch admittances $y_{\alpha}$ satisfy $\arg y_{\alpha} \in [\theta_1, \theta_2)$ where $\theta_2 - \theta_1 < \pi$ and $\theta_2 - \theta_1 \leq \pi - \phi_1 - \phi_2$.
Then for any choice of $\theta$ such that
\begin{alignat}{2}
\phi_2 - \theta_1 - \frac{\pi}{2} \leq \theta &\leq \frac{\pi}{2} - \phi_1 - \theta_2 & \qquad &
   \ivl \text{ nontrivial},
\label{eq:iconetheta} \\
- \theta_1 - \frac{\pi}{2} < \theta &\leq \frac{\pi}{2} - \theta_2 & &
   \ivl = \{0\},
\label{eq:iconetheta0}
\end{alignat}
the quantities
\begin{equation}
d_{jk} = \re{\e^{-j \theta} Z_{jk}} \qquad j, k \in \nw
\label{eq:iconedjk}
\end{equation}
where $Z_{jk}$ are the driving-point impedances, define a metric on the nodes of $\nw$.
The triangle inequality $d_{pq} + d_{qr} \geq d_{pr}$ is an equality if and only if $q$ is a `bridge node' relative to nodes $p$ and $r$, in the sense of Proposition \ref{prop:dcmetric}.

The above remains true when the interval of spread is $\ivl = (-\phi_1, \phi_2]$ and the admittances satisfy $\arg y_{\alpha} \in (\theta_1, \theta_2]$, provided (\ref{eq:iconetheta0}) is replaced with the condition $-\theta_1 - \pi/2 \leq \theta < \pi/2 - \theta_2$.
It also remains true when $\ivl = [-\phi_2, \phi_2]$ and $\arg y_{\alpha} \in (\theta_1, \theta_2)$, with (\ref{eq:iconetheta}) applying without restriction.
\end{proposition}
To verify Proposition \ref{prop:conemetric} it suffices to consider a source with $\phi = 0$.
One may then check that, when $\ivl$ is nontrivial and under the conditions given, one has
\begin{equation}
\arg \paren{y_{\alpha} i_{\alpha}^{-1}}^* \in \paren{-\phi_1 - \theta_2, \phi_2 - \theta_1} \qquad
   \text{for all } \alpha \in \nw
\label{eq:vbrchicone}
\end{equation}
and for any $\theta$ chosen as in (\ref{eq:iconetheta}) the latter interval will be covered by the (open) interval $(\theta - \pi/2, \theta + \pi/2)$.
The quantity $\re{\e^{-j\theta} v}$ will then define an orientable component of voltage that strictly decreases on all branches of $\nw$ with nonzero current as per (\ref{eq:vcone}).
When $\ivl$ is trivial, the only difference is that the interval in (\ref{eq:vbrchicone}) is $(-\theta_2, -\theta_1]$ and is only semi-open, so the choice $\theta = -\theta_1 - \pi/2$ is not available (as the resulting measure may not be strictly positive).

The strict decrease of the quantity $\re{\e^{-j\theta} v}$ on all branches with nonzero current, oriented as required, then allows the reasoning to proceed as for a DC network, showing that the quantities $\re{\e^{-j\theta} v_j}$ and $\re{\e^{-j\theta} v_k}$ are maximal and minimal, respectively, for a source connected from node $k$ to $j$.
To verify that $d_{jk}$ as in (\ref{eq:iconedjk}) is the desired metric, observe that for such a source connection one has
\begin{equation}
\re{\e^{-j\theta} v_j} - \re{\e^{-j\theta} v_k} = \re{\e^{-j\theta} Z_{jk} I} = d_{jk} I
\label{eq:iconedjki}
\end{equation}
(keeping in mind $I$ is assumed to have phase angle zero).
The positivity of the left hand side when $j \neq k$ thus establishes the positivity of $d_{jk}$, and an argument based on linearity as for Proposition \ref{prop:dcmetric} establishes the triangle inequality (with the same topological condition for equality).

Proposition \ref{prop:conemetric} also encapsulates, as special cases, both Proposition \ref{prop:acmetric} for networks with common admittance phase angles, and Proposition \ref{prop:treemetric} for networks that are trees.
All such networks are dendromorphic with a trivial interval of spread, in the latter case due to the tree structure, and in the former case because all currents automatically align in phase (as do all node voltages, with a fixed phase offset between any node voltage and any branch current).

Proposition \ref{prop:acmetric} asserts that $\arg y_{\alpha} = \theta_Y$ for all $\alpha$, where $\theta_Y$ is a constant, hence lies in $[\theta_Y, \theta_Y + \epsilon)$ for any $\epsilon > 0$.
Proposition \ref{prop:conemetric} now asserts that $\re{\e^{-j\theta} Z_{jk}}$ is a metric for any $\theta$ with $|\theta - \theta_Y| < \pi/2$.
This includes $\theta_Y$ itself, ensuring that $|Z_{jk}|$ is a metric, and also either $\pi/2$ or $-\pi/2$ depending on the sign of $\theta_Y \neq 0$, ensuring the unsigned reactance $|X_{jk}|$ is also a metric.
But the case $\theta_Y = 0$ is also now included, establishing the metric principle already noted for resistive networks based on driving-point resistances $R_{jk}$.

Proposition \ref{prop:treemetric} for trees supposes that $\re{\e^{j\theta} y_{\alpha}} > 0$ for some $\theta$, equivalent to assuming $\arg y_{\alpha} \in (-\pi/2 - \theta, \pi/2 - \theta)$.
With a trivial interval of spread, and with an open interval for $\arg y_{\alpha}$, this same value of $\theta$ (but only this value) also defines a metric under Proposition \ref{prop:conemetric}. 

The following is another useful special case of Proposition \ref{prop:conemetric}, which establishes the metric property for certain dendromorphic networks with nontrivial intervals of spread and suitable restrictions on the admittances---including the strictly inductive and strictly capacitive cases.
\begin{proposition}
\label{prop:dendrolc}
Let $\nw$ be an AC network that is dendromorphic with interval of spread $\ivl = [-\phi_1, \phi_2)$ (respectively, $(-\phi_1, \phi_2]$) such that $\phi_1 + \phi_2 \leq \pi/2$, and let all branch admittances $y_{\alpha} \in \nw$ satisfy $\arg y_{\alpha} \in [\theta_1, \theta_2)$ (respectively, $\arg y_{\alpha} \in (\theta_1, \theta_2]$) where $\theta_2 - \theta_1 \leq \pi/2$.
Then for $\theta$ given by (\ref{eq:iconetheta}), the quantities $d_{jk} = \re{\e^{-j\theta} Z_{jk}}$ define a metric on the nodes of $\nw$.

In particular, if the admittances of $\nw$ are strictly inductive ($\arg y_{\alpha} \in [-\pi/2, 0)$) then a metric is defined with $\theta \in [\phi_2, \pi/2 - \phi_1]$ where $\ivl = [-\phi_1, \phi_2)$, while if the admittances are strictly capacitive ($\arg y_{\alpha} \in (0, \pi/2]$) then a metric is defined with $\theta \in [-\pi/2 + \phi_2, -\phi_1]$ where $\ivl = (-\phi_1, \phi_2]$.
\end{proposition}
Observe in particular that where the interval of spread $\ivl$ is symmetric around zero (that is, $\phi_1 = \phi_2 \leq \pi/4$), then for (say) a strictly inductive network, $\theta = \pi/4$ will always define a metric.
The metric in this case is, up to an unimportant scale factor, given by $d_{jk} = R_{jk} + X_{jk}$ where $Z_{jk} = R_{jk} + jX_{jk}$ is the driving-point impedance.

In a certain sense, Propositions \ref{prop:conemetric} and \ref{prop:dendrolc} are as far as one may go in arguing orientability in general AC networks.
There exist meshed networks that are not dendromorphic and in which there is no readily identifiable metric definable from admittances that satisfies the triangle inequality---despite the admittances all being strictly inductive.
Such networks can also violate intuitions regarding orientability in quite radical ways.
An example is now given based on the Wheatstone bridge network already invoked in Section \ref{sec:dcmono}.

As in the previous example, let a four-node network be constructed from admittances $y_{\alpha}$, $y_{\beta}$, $y_{\gamma}$, $y_{\delta}$ arranged in a square, with a source $I$ connected between nodes 1 and 2 and a further impedance $y_{\tau}$ connected across nodes 3 and 4.
(For concreteness the source current $I$ may be taken to be 1A, with phase angle zero.)
There is a guiding assumption that all admittances are strictly inductive, with $\re{y} \geq 0$ and $\im{y} < 0$.
The admittance matrix $Y$ is as in (\ref{eq:whsy}), but with $y_{\sigma}$ omitted for simplicity.
The Kirchhoff characteristic, simplified from (\ref{eq:whsk}), is
\begin{equation}
\kappa(\nw) = y_{\alpha} y_{\beta} y_{\gamma} + y_{\alpha} y_{\beta} y_{\delta}
   + y_{\alpha} y_{\gamma} y_{\delta} + y_{\beta} y_{\gamma} y_{\delta}
   + y_{\tau} \paren{y_{\alpha} y_{\beta} + y_{\alpha} y_{\delta}
      + y_{\beta} y_{\gamma} + y_{\gamma} y_{\delta}}
\label{eq:whsk0}
\end{equation}
and the driving-point impedance at the source terminals was found in (\ref{eq:whsz12}):
\begin{equation}
Z_{12} = \frac{v_1}{I}
   = \frac{1}{\kappa(\nw)} \bbrak{\paren{y_{\alpha} + y_{\beta}} \paren{y_{\gamma} + y_{\delta}}
   + y_{\tau} \paren{y_{\alpha} + y_{\beta} + y_{\gamma} + y_{\delta}}}.
\label{eq:whsv1}
\end{equation}
The internal node voltages $v_3$ and $v_4$ relative to the ground node $g = 2$ (or rather, their ratio to the source current $I$) are found directly from the transfer impedances
\begin{align}
\frac{v_3}{I} &= \tz{12}{32} = \frac{1}{\kappa(\nw)}
   \bbrak{y_{\alpha} \paren{y_{\gamma} + y_{\delta}} + y_{\tau} \paren{y_{\alpha} + y_{\gamma}}},
\label{eq:whsv3} \\
\frac{v_4}{I} &= \tz{12}{42} = \frac{1}{\kappa(\nw)}
   \bbrak{y_{\gamma} \paren{y_{\alpha} + y_{\beta}} + y_{\tau} \paren{y_{\alpha} + y_{\gamma}}}.
\label{eq:whsv4}
\end{align}
The formulae for the branch currents $i_{\beta}$ and $i_{\delta}$ follow directly from those for $v_3$ and $v_4$:
\begin{align}
\frac{i_{\beta}}{I} = \frac{y_{\beta} v_3}{I} &= \frac{1}{\kappa(\nw)}
   \bbrak{y_{\alpha} y_{\beta} \paren{y_{\gamma} + y_{\delta}}
      + y_{\tau} y_{\beta} \paren{y_{\alpha} + y_{\gamma}}},
\label{eq:whsib} \\
\frac{i_{\delta}}{I} = \frac{y_{\delta} v_4}{I} &= \frac{1}{\kappa(\nw)}
   \bbrak{\paren{y_{\alpha} + y_{\beta}} y_{\gamma} y_{\delta}
      + y_{\tau} y_{\delta} \paren{y_{\alpha} + y_{\gamma}}},
\label{eq:whsid}
\end{align}
and one may check from these formulae that $(i_{\beta} + i_{\delta}) / I = 1$, as given by Kirchhoff's current law, is an identity.
The remaining branch currents can be derived from the voltage formulae as
\begin{align}
\frac{i_{\alpha}}{I} = \frac{y_{\alpha} (v_1 - v_3)}{I} &= \frac{1}{\kappa(\nw)}
   \bbrak{y_{\alpha} y_{\beta} \paren{y_{\gamma} + y_{\delta}}
      + y_{\tau} y_{\alpha} \paren{y_{\beta} + y_{\delta}}},
\label{eq:whsia} \\
\frac{i_{\gamma}}{I} = \frac{y_{\gamma} (v_1 - v_4)}{I} &= \frac{1}{\kappa(\nw)}
   \bbrak{\paren{y_{\alpha} + y_{\beta}} y_{\gamma} y_{\delta}
      + y_{\tau} y_{\gamma} \paren{y_{\beta} + y_{\delta}}}, \qquad \text{and}
\label{eq:whsic} \\
\frac{i_{\tau}}{I} = \frac{y_{\tau} (v_3 - v_4)}{I} &= \frac{1}{\kappa(\nw)} \cdot
   y_{\tau} \paren{y_{\alpha} y_{\delta} - y_{\beta} y_{\gamma}}.
\label{eq:whsit}
\end{align}
One observes again in (\ref{eq:whsit}) the conjugacy condition $y_{\alpha} y_{\delta} = y_{\beta} y_{\gamma}$ that zeros the current in branch $\tau$, irrespective of the value of $y_{\tau}$ or the other branch currents.
The Kirchhoff current law identities $(i_{\alpha} + i_{\gamma}) / I = 1$ and $i_{\tau} = i_{\alpha} - i_{\beta} = i_{\delta} - i_{\gamma}$ are also readily verified.

It will also be useful to have expressions for the `triangle' of driving-point impedances $Z_{13}$, $Z_{14}$ and $Z_{34}$.
The latter was found back at (\ref{eq:whsz34}), and simplifies with $y_{\sigma} = 0$ to:
\begin{equation}
Z_{34} = \frac{1}{\kappa(\nw)} \paren{y_{\alpha} + y_{\gamma}} \paren{y_{\beta} + y_{\delta}}.
\label{eq:whsz340}
\end{equation}
For the others one calculates directly from cofactors of (\ref{eq:whsy}):
\begin{equation}
Z_{13} = \frac{1}{\kappa(\nw)} \bbrak{\paren{y_{\beta} + y_{\delta}} \paren{y_{\gamma} + y_{\tau}} + y_{\beta} y_{\delta}}, \qquad
Z_{14} = \frac{1}{\kappa(\nw)} \bbrak{\paren{y_{\beta} + y_{\delta}} \paren{y_{\alpha} + y_{\tau}} + y_{\beta} y_{\delta}}.
\label{eq:whsz1314}
\end{equation}
These formulae have an obvious similarity, and one readily computes
\begin{equation}
Z_{13} - Z_{14} = \frac{1}{\kappa(\nw)} \paren{y_{\beta} + y_{\delta}} \paren{y_{\gamma} - y_{\alpha}},
\qquad \text{and} \qquad
Z_{13} + Z_{34} - Z_{14} = \frac{2 y_{\gamma}}{\kappa(\nw)} \paren{y_{\beta} + y_{\delta}}.
\label{eq:whsz134}
\end{equation}
Observe that a violation of the triangle inequality occurs for a candidate metric $d_{jk} = \re{\e^{-j\theta} Z_{jk}}$ if $\e^{-j\theta} \paren{Z_{13} + Z_{34} - Z_{14}}$ falls in the left half plane.
Factoring $\kappa(\nw)$ from (\ref{eq:whsk0}) as
\begin{equation}
\kappa(\nw) = y_{\gamma} \paren{y_{\beta} + y_{\delta}} \paren{y_{\alpha} + y_{\tau}}
   + y_{\beta} y_{\delta} \paren{y_{\alpha} + y_{\gamma}}
   + y_{\tau} y_{\alpha} \paren{y_{\beta} + y_{\delta}}
\label{eq:whsk134}
\end{equation}
one may combine (\ref{eq:whsz134}) and (\ref{eq:whsk134}) to obtain
\begin{equation}
Z_{13} + Z_{34} - Z_{14} = 2 \paren{y_{\alpha}
   + \paren{\frac{1}{y_{\beta}^{-1} + y_{\delta}^{-1}} + y_{\tau}}
      \paren{1 + \frac{y_{\alpha}}{y_{\gamma}}}}^{-1}.
\label{eq:whstri}
\end{equation}
One now notices that $1 / (y_{\beta}^{-1} + y_{\delta}^{-1})$ is the admittance of the series combination of $(y_{\beta}, y_{\delta})$ and is strictly inductive if the original admittances are; adding $y_{\tau}$ to this likewise yields a strictly inductive admittance, which can however be placed more or less arbitrarily in the quadrant $\arg y \in [-\pi, 0)$.
Meanwhile, the ratio $y_{\alpha} / y_{\gamma}$ can be placed independently in the right half plane.

It is therefore not difficult to see that if both $y_{\alpha}$ and $y_{\tau} + 1 / (y_{\beta}^{-1} + y_{\delta}^{-1})$ are close to $-\pi/2$ in phase, while $y_{\gamma}$ is close to the real axis and $|y_{\alpha} / y_{\gamma}|$ sufficiently large, the quantity (\ref{eq:whstri}) can be placed well within the opposing quadrant $\re{z} < 0$, $\im{z} \geq 0$.
If for example $y_{\gamma}$ is almost purely resistive and the other admittances pure inductors:
\begin{equation}
y_{\alpha} = -j b, \quad y_{\beta} = y_{\delta} = -j1, \quad y_{\tau} = -j (b - 1/2), \quad
y_{\gamma} = 1 - j\epsilon, \quad (b \gg 1, \epsilon \ll 1)
\label{eq:whstrival}
\end{equation}
then one readily calculates
\begin{equation}
\begin{split}
Z_{13} + Z_{34} - Z_{14}
   &= 2 \paren{-jb - jb \paren{1 - \frac{b (j - \epsilon)}{1 + \epsilon^2}}}^{-1}
   = - 2 (1 + \epsilon^2) \bbrak{b^2 + j b \paren{2 + \epsilon b + 2 \epsilon^2}}^{-1} \\
   &\approx - \frac{2}{b^2} \brak{1 - j \paren{\frac{2}{b} + \epsilon}}.
\end{split}
\label{eq:whstri1}
\end{equation}
This places the `triangle defect' close to the negative real axis in $\cplx$, with $\arg\paren{Z_{13} + Z_{34} - Z_{14}} \approx \pi - \arctan\paren{2/b + \epsilon}$.
Meanwhile, by exchanging $y_{\alpha}$ and $y_{\gamma}$ in (\ref{eq:whsz134}) one also calculates that
\begin{equation}
Z_{14} + Z_{34} - Z_{13} = 2 \paren{y_{\gamma}
   + \paren{\frac{1}{y_{\beta}^{-1} + y_{\delta}^{-1}} + y_{\tau}}
      \paren{1 + \frac{y_{\gamma}}{y_{\alpha}}}}^{-1},
\label{eq:whstrix}
\end{equation}
and with admittance values as in (\ref{eq:whstrival}),
\begin{equation}
Z_{14} + Z_{34} - Z_{13} = 2 \paren{1 - j \epsilon - jb \paren{1 + \frac{j + \epsilon}{b}}}^{-1}
   = \paren{1 - j \paren{\frac{b}{2} + \epsilon}}^{-1} \approx \frac{1 + j b / 2}{(b / 2)^2 + 1}.
\label{eq:whstrix1}
\end{equation}
This results in $\arg\paren{Z_{14} + Z_{34} - Z_{13}} \approx \arctan(b/2)$, which for $b \gg 1$ is close to $\pi/2$.

Notice now that for a quantity of the form $d_{jk} = \e^{-j\theta} Z_{jk}$ to satisfy the triangle inequality and therefore serve as a metric over $\nw$, based on (\ref{eq:whstri1}) and (\ref{eq:whstrix1}) $\theta$ would need to satisfy as a minimum
\begin{equation}
\frac{\pi}{2} - \arctan\paren{\frac{2}{b} + \epsilon} \leq \theta \leq \frac{\pi}{2} + \arctan\frac{b}{2}.
\label{eq:whstheta}
\end{equation}
For $b$ sufficiently large this excludes most values of $\theta$ less than $\pi/2$---in particular it excludes $\theta = \pi/4$ ($d_{jk} \sim R_{jk} + X_{jk}$), which Proposition \ref{prop:dendrolc} otherwise suggests as the most promising candidate for a general metric valid over a large class of strictly inductive networks.

For a specific numerical example, take $b = 10$ and $\epsilon = 1/10$; then one may calculate
\begin{equation}
\begin{split}
\kappa(\nw) = -21 + \frac{121}{10} j + \paren{-\frac{19}{2} j} \paren{-\frac{101}{5} - 2j} &= -40 + 204j,
   \qquad \frac{1}{\kappa(\nw)} = -\frac{10 + 51j}{10804}, \\
Z_{13} = -\frac{10 + 51j}{10804} \brak{(-2j) \paren{1 - \frac{48}{5} j} - 1}
   &= \frac{500 + 5251j}{54020}, \\
Z_{14} = -\frac{10 + 51j}{10804} \brak{(-2j) \paren{-\frac{39}{2} j} - 1}
   &= \frac{100 + 510j}{2701}, \\
Z_{34} = -\frac{10 + 51j}{10804} \brak{(-2j) \paren{1 - \frac{101}{10} j}}
   &= \frac{500 + 5251j}{54020} \:\paren{= Z_{13}}.
\end{split}
\label{eq:whsnum}
\end{equation}
Accordingly one sees that
\begin{equation}
Z_{13} + Z_{34} = \frac{500 + 5251j}{27010}, \qquad Z_{14} = \frac{1000 + 5100j}{27010}
\label{eq:whsnumz}
\end{equation}
and one may verify that while the triangle inequality $d_{13} + d_{34} \geq d_{14}$ is satisfied for $d_{jk} = X_{jk}$ in this example, it is violated for both $d_{jk} = R_{jk}$ and for $d_{jk} = R_{jk} + X_{jk}$.

But now suppose that in place of $\nw$ one considers a network $\nw'$ configured identically other than having its admittances replaced with new values
\begin{equation}
y_{\alpha}' = b - j\epsilon, \qquad y_{\beta}' = y_{\delta}' = 1 - j\epsilon, \qquad
y_{\gamma}' = -j, \qquad y_{\tau} = (b - 1/2) - j\epsilon,
\label{eq:whstrivalx}
\end{equation}
obtained by reflecting the values (\ref{eq:whstrival}) through the line $\arg y = -\pi/4$ (converting conductances to susceptances and vice versa) and perturbing by $-j\epsilon$ where necessary to keep the values in the region $\im{y} < 0$.
Calculating explicitly once again with $b = 10$ and $\epsilon = 1/10$ (and defining $\zeta = 1 - (1/10) j$ for convenience):
\begin{equation}
\begin{split}
\kappa(\nw) = \paren{\zeta^2 - 2 j \zeta} \paren{9 + \zeta} - j \zeta^2
   + \paren{17 + 2 \zeta} \paren{\zeta \paren{9 + \zeta} - j \zeta} & \\
   = \zeta \brak{\paren{9 - j + \zeta} \paren{17 - 2j + 3 \zeta} + 2}
   = \zeta \paren{199.47 - 45j}, & \qquad \frac{\zeta}{\kappa} = \frac{664900 + 150000j}{139377603}, \\
Z_{13} = \frac{1}{\kappa} \brak{2 \zeta \paren{\frac{17}{2} - j + \zeta} + \zeta^2}
   = \frac{\zeta}{\kappa} \paren{17 - 2j + 3\zeta}
   &= \frac{13643000 + 1470730j}{139377603}, \\
Z_{14} = \frac{1}{\kappa} \brak{2 \zeta \paren{\frac{35}{2} + 2 \zeta} + \zeta^2}
   = \frac{\zeta}{\kappa} \paren{35 + 5 \zeta}
   &= \frac{26671000 + 5667550j}{139377603}, \\
Z_{34} = \frac{1}{\kappa} \paren{9 - j + \zeta} \cdot 2\zeta = \frac{\zeta}{\kappa} \paren{20 - 2.2j}
   &= \frac{13628000 + 1537220j}{139377603},
\end{split}
\label{eq:whsxnum}
\end{equation}
and hence
\begin{equation}
Z_{13} + Z_{34} = \frac{27271000 + 3007950j}{139377603}, \qquad
Z_{14} = \frac{26671000 + 5667550j}{139377603}.
\label{eq:whsxnumz}
\end{equation}
It is now clear that while $d_{jk} = R_{jk}$ will satisfy the triangle inequality $d_{13} + d_{34} \geq d_{14}$, this inequality now fails for both $d_{jk} = X_{jk}$ and for $d_{jk} = R_{jk} + X_{jk}$.

These examples together demonstrate by construction the following result, which answers in the negative the question of finding a metric $d_{jk}$ as a component of the driving-point impedances $Z_{jk}$ that is valid across all strictly inductive networks without further restrictions.
\begin{proposition}
\label{prop:noacmetric}
Fix a value $\theta \in \reals$ and for an AC network $\nw$ let $d_{jk} = \e^{-j\theta} Z_{jk}$ for all node pairs $(j,k) \in \nw$ be a candidate for a distance metric.
Then there is a choice of strictly inductive branch admittances for $\nw$, having $\re{y_{\alpha}} \geq 0$ and $\im{y_{\alpha}} < 0$ for all branches $\alpha \in \nw$, such that the triangle inequality $d_{pq} + d_{qr} \geq d_{pr}$ is violated for at least one triple of nodes $(p,q,r) \in \nw$.
Thus, there is no metric of this form universally valid over the class of strictly inductive networks.
\end{proposition}

\subsection{Phase angles, inductive flows and semi-orientability}
\label{sec:acphase}

This section further develops the relationship between power flows and AC voltages, and in particular indicates how AC voltage phase angles can be weakly orientable in practical cases, based on a limited version of the `power flows downhill' principle.

Suppose the AC network $\nw$ contains no zero-admittance branches, so that every branch $\alpha$ can be described by its complex impedance $z_{\alpha} = y_{\alpha}^{-1} = r_{\alpha} + j x_{\alpha}$, and the fundamental Ohm's law relation can be written
\begin{equation}
v_{\alpha+} - v_{\alpha-} = z_{\alpha} i_{\alpha},
\label{eq:acohm}
\end{equation}
where the ordered pair $(\alpha+, \alpha-)$ denotes the endpoint node indices.
If one multiplies through by the conjugate voltage $v_{\alpha-}^*$, this relation becomes
\begin{equation}
v_{\alpha+} v_{\alpha-}^* - |v_{\alpha-}|^2 = z_{\alpha} v_{\alpha-}^* i_{\alpha}
   = z_{\alpha} s_{\alpha-}^* = \paren{r_{\alpha} + j x_{\alpha}} \paren{p_{\alpha-} - j q_{\alpha-}}
\label{eq:acohm1}
\end{equation}
where the notations $s_{\alpha-}$, $p_{\alpha-}$, $q_{\alpha-}$ specifically denote the power flows measured at the second endpoint $\alpha-$ of the branch $\alpha$.
Equation (\ref{eq:acohm1}) can then be restated as
\begin{equation}
v_{\alpha+} v_{\alpha-}^* = |v_{\alpha+}| |v_{\alpha-}| \e^{j \paren{\delta_{\alpha+} - \delta_{\alpha-}}}
   = |v_{\alpha-}|^2 + r_{\alpha} p_{\alpha-} + x_{\alpha} q_{\alpha-}
      + j \paren{x_{\alpha} p_{\alpha-} - r_{\alpha} q_{\alpha-}}
\label{eq:acdelta1}
\end{equation}
where the connection with the voltage phase angles $\delta_{\alpha+}$, $\delta_{\alpha-}$ starts to become more apparent.
If instead one multiplies (\ref{eq:acohm}) through by $v_{\alpha+}^*$, one obtains
\begin{equation}
|v_{\alpha+}|^2 - v_{\alpha+}^* v_{\alpha-} = z_{\alpha} v_{\alpha+}^* i_{\alpha}
   = z_{\alpha} s_{\alpha+}^* = \paren{r_{\alpha} + j x_{\alpha}} \paren{p_{\alpha+} - j q_{\alpha+}}
\label{eq:acohm2}
\end{equation}
and
\begin{equation}
v_{\alpha+}^* v_{\alpha-} = |v_{\alpha+}| |v_{\alpha-}| \e^{-j \paren{\delta_{\alpha+} - \delta_{\alpha-}}}
   = |v_{\alpha+}|^2 - r_{\alpha} p_{\alpha+} - x_{\alpha} q_{\alpha+}
      - j \paren{x_{\alpha} p_{\alpha+} - r_{\alpha} q_{\alpha+}},
\label{eq:acdelta2}
\end{equation}
where $s_{\alpha+}$, $p_{\alpha+}$, $q_{\alpha+}$ denote the power flows measured at the first endpoint $\alpha+$ on the branch $\alpha$.

Equations (\ref{eq:acdelta1}) and (\ref{eq:acdelta2}) provide two identities for the complex quantity $v_{\alpha+} v_{\alpha-}^*$ and its conjugate.
It follows that when considering the right hand sides of these identities, one may equate the real parts, and also equate the imaginary parts with a change of sign.
Equating the real parts gives the identity
\begin{equation}
|v_{\alpha+}|^2 - |v_{\alpha-}|^2
   = r_{\alpha} \paren{p_{\alpha+} + p_{\alpha-}} + x_{\alpha} \paren{q_{\alpha+} + q_{\alpha-}}
\label{eq:acv2}
\end{equation}
relating the voltage magnitudes at the endpoints.
This can be rendered more useful in practice by utilising the additional identities
\begin{equation}
p_{\alpha+} - p_{\alpha-} = r_{\alpha} |i_{\alpha}|^2, \quad
q_{\alpha+} - q_{\alpha-} = x_{\alpha} |i_{\alpha}|^2, \quad \text{and} \quad
p_{\alpha\pm}^2 + q_{\alpha\pm}^2 = |v_{\alpha\pm}|^2 |i_{\alpha}|^2.
\label{eq:aci2}
\end{equation}
The first two of these follow directly from (\ref{eq:acohm}) by multiplying through by $i_{\alpha}^*$, noting that the left hand side becomes $s_{\alpha+} - s_{\alpha-}$, and resolving into real and imaginary components.
The third (a separate identity for each choice of sign) is a corollary of the defining formulae for $s_{\alpha+}$ and $s_{\alpha-}$ when magnitudes are taken.
Observe that the first two identities give
\begin{equation}
\paren{r_{\alpha}^2 + x_{\alpha}^2} |i_{\alpha}|^2
   = r_{\alpha} \paren{p_{\alpha+} - p_{\alpha-}} + x_{\alpha} \paren{q_{\alpha+} - q_{\alpha-}},
\label{eq:aci2rx}
\end{equation}
and by combining with (\ref{eq:acv2}) this leads to (by addition)
\begin{equation}
|v_{\alpha+}|^2 - |v_{\alpha-}|^2 = 2 \paren{r_{\alpha} p_{\alpha+} + x_{\alpha} q_{\alpha+}}
   - \paren{r_{\alpha}^2 + x_{\alpha}^2} |i_{\alpha}|^2
\label{eq:acv2a}
\end{equation}
or alternatively to (by subtraction)
\begin{equation}
\begin{split}
|v_{\alpha+}|^2 - |v_{\alpha-}|^2 &= 2 \paren{r_{\alpha} p_{\alpha-} + x_{\alpha} q_{\alpha-}}
      + \paren{r_{\alpha}^2 + x_{\alpha}^2} |i_{\alpha}|^2 \\
  &= 2 \paren{r_{\alpha} p_{\alpha-} + x_{\alpha} q_{\alpha-}}
      + \frac{\paren{r_{\alpha}^2 + x_{\alpha}^2} \paren{p_{\alpha-}^2 + q_{\alpha-}^2}}{|v_{\alpha-}|^2}.
\end{split}
\label{eq:acv2b}
\end{equation}
Either (\ref{eq:acv2a}) or (\ref{eq:acv2b}) provides a criterion for when the voltage magnitudes are higher at one end of the branch than another; for example, when $r_{\alpha} > 0$ and $x_{\alpha} > 0$, (\ref{eq:acv2b}) shows that $|v_{\alpha+}| > |v_{\alpha-}|$ whenever $p_{\alpha-}$ and $q_{\alpha-}$ are nonnegative with at least one positive.
But (\ref{eq:acv2b}) also---after multiplying through by $|v_{\alpha-}|^2$ to clear fractions---gives a useful equation that can be solved for an unknown `near end' voltage $|v_{\alpha-}|$ given a known `far end' voltage $|v_{\alpha+}|$ and `near end' power flows $p_{\alpha-}$ and $q_{\alpha-}$ on branch $\alpha$.
Together with a more straightforward application of the identities (\ref{eq:aci2}), this gives
\begin{proposition}
\label{prop:acv}
Let $\alpha$ be a branch in a network $\nw$ with AC impedance $r_{\alpha} + j x_{\alpha}$, directed from node $j$ with known voltage magnitude $|v_j|$ to node $k$ with unknown voltage $v_k$.
If the complex power flow $s_j = p_j + j q_j$ from node $j$ on branch $\alpha$ is known, then the voltage magnitude $|v_k|$ is given by
\begin{equation}
|v_k|^2 = \frac{p_k^2 + q_k^2}{|i_{\alpha}|^2} \quad \text{where} \quad
|i_{\alpha}|^2 = \frac{p_j^2 + q_j^2}{|v_j|^2}, \quad
p_k = p_j - r_{\alpha} |i_{\alpha}|^2, \quad q_k = q_j - x_{\alpha} |i_{\alpha}|^2.
\label{eq:acvsj}
\end{equation}
If on the other hand the complex power flow $s_k = p_k + j q_k$ toward node $k$ on branch $\alpha$ is known, then the following provides a quadratic equation for $|v_k|^2$ where $|v_k|$ is the unknown voltage magnitude:
\begin{equation}
|v_k|^4 - \brak{|v_j|^2 - 2 \paren{r_{\alpha} p_k + x_{\alpha} q_k}} |v_k|^2
   + \paren{r_{\alpha}^2 + x_{\alpha}^2} \paren{p_k^2 + q_k^2} = 0.
\label{eq:acvsk}
\end{equation}
If instead the voltage magnitude $|v_k|$ is known and $v_j$ unknown, the above formulae apply after interchanging the indices $j$ and $k$ and changing the signs of $p_j$, $q_j$, $p_k$, $q_k$.
The modified formula (\ref{eq:acvsj}) applies when the flow $s_k = p_k + j q_k$ into node $k$ is known, and the modified formula (\ref{eq:acvsk}) when $s_j = p_j + j q_j$ from node $j$ is known.
\end{proposition}
Proposition \ref{prop:acv} relates the voltage \emph{magnitudes} at the two ends of an AC circult branch with known power flow, using formula (\ref{eq:acv2}) which was found by equating the real parts of the two formulae (\ref{eq:acdelta1}) and (\ref{eq:acdelta2}).
If instead one equates the imaginary parts of these two formulae (requiring a change of sign in the process), making these explicitly equal to the imaginary part of $v_{\alpha+} v_{\alpha-}^*$, this yields an important result for phase angles.
\begin{proposition}
\label{prop:acdelta}
Let $\alpha$ be a branch in a network $\nw$ with AC impedance $r_{\alpha} + j x_{\alpha}$, directed from node $j$ to node $k$.
Let the complex power flow on the branch $\alpha$ be $s_j = p_j + j q_j$ out of node $j$ and $s_k = p_k + j q_k$ into node $k$.
Then the following is an identity:
\begin{equation}
\mu_{\alpha} = x_{\alpha} p_j - r_{\alpha} q_j = x_{\alpha} p_k - r_{\alpha} q_k
   = |v_j| |v_k| \sin\paren{\delta_j - \delta_k},
\label{eq:acdelta}
\end{equation}
where $\delta_j = \arg v_j$ and $\delta_k = \arg v_k$.
In particular, $\mu_{\alpha} = 0$ whenever $v_j = 0$ or $v_k = 0$; in all other cases, the voltage at node $k$ lags, leads, or matches that at node $j$ in phase, according as $\mu_{\alpha}$ is positive, negative or zero, respectively.
Further, if $v_j, v_k \neq 0$ and the power flow satisfies either of
\begin{equation}
r_{\alpha} p_j + x_{\alpha} q_j \leq |v_j|^2 \qquad \text{or} \qquad
r_{\alpha} (-p_k) + x_{\alpha} (-q_k) \leq |v_k|^2,
\label{eq:acdeltau}
\end{equation}
then $\delta_j - \delta_k$ is uniquely determined in the range $[-\pi/2, \pi/2]$ by (\ref{eq:acdelta}).
\end{proposition}
Note the `end-invariance' implied by the identity (\ref{eq:acdelta}): one may calculate the coefficient $\mu_{\alpha}$ as $x_{\alpha} p - r_{\alpha} q$, irrespective of the end at which the power flows $p$, $q$ are measured---provided only that the branch orientation is respected.
That the differences in $\mu_{\alpha}$ between the two ends cancel is assured by the identities (\ref{eq:aci2}).
Condition (\ref{eq:acdeltau}) is equivalent to $\re{v_j v_k^*} \geq 0$ by (\ref{eq:acdelta1}) and (\ref{eq:acdelta2}).

Proposition \ref{prop:acdelta} generalises an important identity relating active power and phase angles for a \emph{lossless} branch having $r_{\alpha} = 0$.
In that case one may write as a corollary to (\ref{eq:acdelta})
\begin{equation}
p_j = p_k = p_{\alpha} = \frac{|v_j| |v_k|}{x_{\alpha}} \sin\paren{\delta_j - \delta_k}.
\label{eq:acdelta0}
\end{equation}
It follows that when $r_{\alpha} = 0$ and $x_{\alpha} > 0$, and $v_j$ and $v_k$ are both nonzero, there is a direct correspondence between the sign of $p_{\alpha}$ and that of the voltage phase lag from nodes $j$ to $k$.
This makes the node voltages `semi-orientable' based on the power flow, with the orientability confined to phase angles rather than magnitudes and excluding any node at zero voltage.

Realistic networks have $r_{\alpha} > 0$ and $p_j \neq p_k$ in general, but the semi-orientability in the lossless case can be preserved under more general conditions, as the following proposition describes.
\begin{proposition}
\label{prop:induc}
Let $\alpha$ be a branch in a network $\nw$ with impedance $r_{\alpha} + j x_{\alpha}$, directed from node $j$ to node $k$.
Let the complex power flow on the branch $\alpha$ be $s_j = p_j + j q_j$ out of node $j$ and $s_k = p_k + j q_k$ into node $k$, and let $\mu_{\alpha}$ be defined as in (\ref{eq:acdelta}).
The branch $\alpha$ is said to be \emph{inductively loaded} if
\begin{equation}
\mu_{\alpha} p_j \geq 0 \qquad \text{and} \qquad \mu_{\alpha} p_k \geq 0,
\label{eq:induc}
\end{equation}
with equality occurring in either case only if an endpoint node voltage is zero (the branch is a \emph{shunt}), or if $r_{\alpha} = 0$ (the branch is lossless).
If a branch $\alpha$ is inductively loaded and not a shunt, then:
\begin{enumerate}
\item
There is a well-defined voltage phase lag $\delta_j - \delta_k$ whose sign (positive, negative or zero) matches the sign of both $p_j$ and $p_k$.
\item
If (\ref{eq:acdeltau}) holds, then $\mu_{\alpha}$ uniquely determines the voltage phase lag in $[-\pi/2,\pi/2]$ by (\ref{eq:acdelta}).
\item
Zero active power flow $p = 0$ at either end of the branch implies $p = 0$ at both ends.
\item
If the active power flow $p$ is nonzero at both ends, its direction is the same at both ends.
\end{enumerate}
A sufficient condition for a branch $\alpha$ to be inductively loaded is that $r_{\alpha} = 0$ and $x_{\alpha} > 0$, irrespective of its power flow.
The branch is then said to be \emph{(purely) inductive}.
With $r_{\alpha}, x_{\alpha} > 0$, it suffices that at each end of $\alpha$, either $p = q = 0$, or $p$ and $q$ have opposite signs, or $|q| < \paren{x_{\alpha} / r_{\alpha}} |p|$.
\end{proposition}
That a purely inductive branch is automatically inductively loaded ensures that Proposition \ref{prop:induc} provides the desired generalisation of the lossless branch formula (\ref{eq:acdelta0}).

The other properties of inductively loaded branches in Proposition \ref{prop:induc} are easily checked.
In particular, a branch with zero active power flow at either end (say $j$) cannot be inductively loaded unless it is lossless or a shunt---both of which imply that $\mu_{\alpha} = 0$ also---and in the former case will evidently have zero active power flow at both ends.
A branch with zero current flow ($p_j = p_k = q_j = q_k = 0$) may or may not meet the criteria of Proposition \ref{prop:induc} but satisfies the remaining conditions \emph{a fortiori} (and by virtue of Ohm's law has $v_j = v_k$, with zero phase lag), while a lossless branch is inductively loaded only if $x_{\alpha} > 0$ (the pure inductive case) or its active power is zero at both ends.
Either way, an inductively loaded branch with $p = 0$ at one end necessarily has $p = 0$ at both ends unless it is a shunt.

For a shunt, the power flow at the end opposite a zero-voltage (grounded) end is identical to the power generated or consumed in the shunt itself.
Any shunt in a network automatically has $\mu_{\alpha} = 0$, though lacks a well-defined phase angle difference between its endpoints.
All other inductively loaded branches have a well-defined phase lag between their endpoints that matches the sense of the active power flow according to Proposition \ref{prop:induc}.

It is possible for a (non-shunt) branch with $r_{\alpha} > 0$ to fail the inductive loading condition.
The `edge case' is where $\mu_{\alpha} = 0$ with nonzero endpoint voltages: this occurs when $p_j / q_j = p_k / q_k = r_{\alpha} / x_{\alpha}$.
The voltages $v_j$ and $v_k$ then have no relative phase shift, yet there is a nonzero active power flow.
It nonetheless follows in this case that the complex power flows $s_j = v_j i_{\alpha}^*$ and $s_k = v_k i_{\alpha}^*$ are also in phase, ensuring that the active (and reactive) power flows, if nonzero, are still of the same sign.

Other cases breaching Proposition \ref{prop:induc} include those where $x_{\alpha} < 0$ and $p \neq 0$, and those having relatively large resistance $r_{\alpha}$ and substantial reactive power flow $q$.
However, in practice non-shunt branches typically have $x_{\alpha} > 0$, $r_{\alpha}$ of similar or lesser magnitude to $x_{\alpha}$, and power flows satisfying $|q| < (x_{\alpha} / r_{\alpha}) |p|$, which generally suffices to ensure $\mu_{\alpha} p > 0$ and establish inductive loading.

The following proposition is the AC equivalent to Proposition \ref{prop:dcextrema} for DC networks, and applies when the sources in the network are configured in a practically relevant way and all branches in the network solution are inductively loaded, forming what will be termed an \emph{inductive flow}.
\begin{proposition}
\label{prop:acflow}
Let $\nw$ be a connected AC network with all branch admittances nonzero, such that all branches have an impedance $r_{\alpha} + j x_{\alpha}$ with $r_{\alpha} \geq 0$.
Let one node $g$ be distinguished as a grounding point, having the property that the network $\nw - g$, having $g$ and its incident branches removed, remains connected.
Suppose that one or more sources are connected from the common grounding point $g$ to distinct nodes in $\nw$, and that in the resulting network solution grounded at $g$, all branches in $\nw$ are inductively loaded.
The solution is then termed an \emph{inductive flow} grounded at $g$.
For any such inductive flow, there exists an assignment of real numbers $\delta_k$ to all network nodes $k \in \nw$ having $v_k \neq 0$, with the following properties:
\begin{enumerate}
\item
$\delta_k$ is congruent modulo $2\pi$ to the phase angle $\arg v_k$.
\item
If all branch resistances $r_{\alpha}$ are zero and all source currents agree in phase angle (modulo $\pi$) then all angles $\delta_k$ take the same value, which differs by $\pi/2$ from all source currents (modulo $2\pi$), and there is no active power flow in $\nw$.
In all other cases, there is at least one source acting as a \emph{generator} of active power $p$ into $\nw$.
\item
The maximal value among the $\delta_k$ occurs at the terminal of a source which---except in the zero-$p$ case immediately above---acts as a generator.
Any other node with the same maximal value $\delta_k$ is connected to such a source terminal by a path in $\nw$ on which the active power flow $p$ is zero.
\item
Whenever two nodes $j$ and $k$ (having $v_j, v_k \neq 0$) are connected by a branch $\alpha \in \nw$, the quantity $\delta_j - \delta_k$ has the same sign as the power flow $p_{\alpha}$ on the branch from $j$ to $k$ (as measured at either endpoint).
Further, if the power flow satisfies one of the conditions (\ref{eq:acdeltau}), $\delta_j - \delta_k$ is uniquely determined from $\mu_{\alpha}$, $|v_j|$ and $|v_k|$ by (\ref{eq:acdelta}) with magnitude at most $\pi/2$.
\item
The assignment of $\delta_k$ is unique up to a shift $\delta_k' = \delta_k + 2 \pi m$ for all $k$, with $m$ an integer.
\end{enumerate}
\end{proposition}
Much of Proposition \ref{prop:acflow} follows immediately upon applying Proposition \ref{prop:induc} to the full set of branches in $\nw$---the exceptions being the maximum principle for the $\delta_k$, and the exceptional case when all $\delta_k$ are equal.
It is convenient to first dispose of this exceptional case (which is of course also a special case of Proposition \ref{prop:acmetric}).
Observe first that if any $r_{\alpha}$ is positive (and the others nonnegative), Tellegen's theorem asserts that the total power consumed within $\nw$ is positive and will be balanced by net positive generated power from the sources.
Since active power $p$ is a real number it directly follows that there is at least one source having positive generated power.
Now suppose all resistances are in fact zero, and observe that the admittance matrix of $\nw$ then has the form $Y = jB$ where $B$ is a real matrix.
By Proposition \ref{prop:acmetric}, or the definition of a transfer impedance from cofactors, it follows that all transfer impedances $\tz{qg}{kg}$ in $\nw$---those that determine the contribution to the node voltage $v_k$ from the current $I_q$ in a source at node $q$---are also of the form $jX$ with $X$ real.
Then by linearity, it follows that if there is only one source, or the source currents agree in phase mod $\pi$, the nonzero voltages in $\nw$ must also agree mod $\pi$ (all being exactly $\pi/2$ out of phase with any source current).
The connectivity of $\nw - g$ then ensures that all phase angles are in fact equal; Proposition \ref{prop:induc} consequently asserts that the active power flow in each branch is zero, and Proposition \ref{prop:kclpq} that all sources (being connected at distinct terminals) neither generate nor consume active power.

For the converse, suppose that all branch resistances in $\nw$ are zero and that no source acts as a generator.
By Tellegen's theorem the total active power $p$ contributed from all sources is zero, and since by assumption none of these contributions is positive, all source contributions $p$ must in fact be zero.
Assuming there are nonzero voltages $v_k \in \nw$ (hence at least one source current is nonzero), all these voltages must have the same phase $\delta_k$ given there is no active power flow and $\nw - g$ is connected.
(That all $\delta_k$ must be equal when there is no generator will also follow from the proof of the maximum principle, below.)
But then any source with nonzero current must have a current exactly $\pi/2$ out of phase with its terminal voltage in order that its active power $p$ is zero.
With all $\delta_k$ being equal, this implies all source currents agree in phase angle modulo $\pi$.

(Note it is possible for active power to be transferred in a network with no resistances, despite such a network having zero internal power consumption.
But this requires the existence of multiple sources, at least one of which acts as a generator and at least one as a sink of active power.
These sources will necessarily have currents that do not agree in phase modulo $\pi$.
And just as in networks with positive resistances, the phase angles $\delta_k$ in such networks will strictly decrease in the direction of power transfer.)

The maximum principle, meanwhile, is demonstrated in almost the same way as in the proof of Proposition \ref{prop:dcextrema}.
Since each branch is inductively loaded, any branch that is not a shunt and has nonzero $p$ flow can be given a well-defined orientation agreeing with the direction of $p$ at both ends.
By Proposition \ref{prop:induc} this orientation (from node $j$ to $k$ say) is necessarily that for which $\delta_j - \delta_k > 0$.
If on the other hand the branch is a shunt, an end with nonzero $p$ is necessarily the unique end with nonzero voltage, and the branch can likewise be given a unique orientation agreeing with the direction of $p$.
But this direction is necessarily from the end at nonzero voltage to that at zero voltage, as the contrary would imply the shunt is a net generator of power, contradicting the assumption that $r_{\alpha} \geq 0$ for all branches.
This has established the following:
\begin{itemise}
\item
the only non-oriented branches in $\nw$ are those with zero power flow $p$ at both ends, all of which are shunts or have $\delta_j = \delta_k$ by Proposition \ref{prop:induc};
\item
for all oriented branches, the orientation agrees with the direction of their nonzero power flow $p$ and also (unless the branch is a shunt) with the sign of $\delta_j - \delta_k$; and
\item
nodes at zero voltage cannot have outgoing branches within $\nw$, the only ones being sources.
\end{itemise}
The proof of the maximum principle can now be carried through in essentially the same way as that of Proposition \ref{prop:dcextrema}, but with the role of current played by the active power flows $p$, and that of voltage by the phase angles $\delta_k$.
Since Kirchhoff's current law generalises to the power flows $p$ at each node with nonzero voltage by Proposition \ref{prop:kclpq}, it follows that if $k$ is a node in $\nw$ that has an outgoing directed branch (and hence has $v_k \neq 0$) and is \emph{not} the terminal of a source with positive generated power $p$, then $k$ must also have an incoming directed branch within $\nw$.
Let $k$ be a node whose angle $\delta_k$ is maximal under the hypothetical assignment.
If $k$ is the terminal of a generator, there is no more to say.
Otherwise, $k$ has no incoming branches (which imply another node has an angle greater than $\delta_k$) hence has no outgoing branches either.
By connectedness of $\nw$, $k$ is connected to other nodes, at least one of which (say $j$) has $v_j \neq 0$ since an isolated shunt cannot maintain nonzero voltage with angle $\delta_k$ by itself.
Carrying through as before, this implies the existence of a path within $\nw$ from any node $k$ whose angle $\delta_k$ is maximal to the terminal node of a source, with all branches along this path having zero active power flow and all nodes having the same angle $\delta_k$, as required.
But again by connectedness of $\nw$, if this fails to exhaust all nodes with nonzero voltage, there must exist some connection from a node with angle $\delta_k$ to a node $j$ with $\delta_j < \delta_k$, implying the existence of an outgoing branch from a node at $\delta_k$ (which is necessarily the source terminal) and that the source has positive generated power to balance the outgoing branch.
This establishes that the maximal $\delta_k$ must occur at the terminal of a source, and this source must have positive generated power \emph{except} in the case where all $\delta_k$ are identical.

Proposition \ref{prop:acflow} is now fully established.
It may readily be applied to practical AC networks where active power $p$ and reactive power $q$ are generated or consumed within active sources, sinks and shunts that are assumed to be connected between an `active' terminal in the network and a common ground point.
This also applies to balanced three-phase networks, with each `active' terminal representing a set of three individual phase busbars with balanced voltages and currents in suitable phase sequence, and the phase angle being that in one designated representative phase.

The principle of inferring power flows in AC networks from voltage phase angles, that `lag' in the direction of net flow, is sometimes described as \emph{DC load flow}.
This is not because it mimics the behaviour of DC networks (except by analogy with the orientability principle) but because it was the principle behind a type of analogue computer (the \emph{DC network analyser}) used for power system analysis in the 20th century.
This approximated (\ref{eq:acdelta0}) in the per-unit system as $p_{\alpha} \approx \paren{\delta_j - \delta_k} / x_{\alpha}$, which can be likened to Ohm's law for DC circuits with $\delta_k$ as voltage and $p_{\alpha}$ as current.

\subsection{The positive-real property and its consequences}
\label{sec:acposreal}

Further important properties of admittance and impedance functions become apparent when considering the general dynamic response in terms of the Laplace variable $s$.

It was seen with the aid of formulae (\ref{eq:gbs}) that when a general AC branch admittance or impedance is expressed as a function $y_{\alpha}(s)$, respectively $z_{\alpha}(s)$, of the complex variable $s = \sigma + j\omega$, both have the property that their real part $G = \re{y_{\alpha}}$ or $R = \re{z_{\alpha}}$ is always nonnegative provided $\sigma$ is nonnegative.
This key \emph{positive-real} property is valid not only for individual branch admittances and impedances, but also for the driving-point impedances $Z_{jk}$ of the network and their reciprocals.

\begin{proposition}[Cauer--Brune]
\label{prop:posreal}
Let a network be composed of admittances of the form (\ref{eq:gcrl}), where $s$ is a complex variable, and let $Y(s)$ be the corresponding admittance matrix.
Then for any nodes $j$, $k$ the driving-point impedance $Z_{jk}(Y(s))$ and its reciprocal $1 / Z_{jk}(Y(s))$ are both rational functions of $s$ with real coefficients that satisfy the \emph{positive-real} property.
A rational function $f(s)$ with $s \in \cplx$ is positive-real (also known as a \emph{Brune function}) if at all points $s$ other than poles:
\begin{itemise}
\item $f(s)$ is real whenever $s$ is real; and
\item $\re{f(s)} \geq 0$ whenever $\re{s} \geq 0$.
\end{itemise}
\end{proposition}
It follows immediately from the definition of the positive-real property in Proposition \ref{prop:posreal} that the sum of any two positive-real functions is positive-real (since $\re{z_1 + z_2} = \re{z_1} + \re{z_2}$), and also that the \emph{composition} of two positive-real functions is positive-real.
It is therefore sufficient to prove Proposition \ref{prop:posreal} just for impedances $Z_{jk}$, as the property for admittances $1 / Z_{jk}$ then follows from the fact that $f(z) = 1/z$ is itself a positive-real function.

That $Z_{jk}$ and $1 / Z_{jk}$ for any pair $(j,k)$ are rational functions of $s$ with real coefficients (and hence real for real $s$), is a straightforward consequence of their construction, either as cofactors of $Y(s)$ or as Kirchhoff characteristics built up from the admittances $y_{\alpha}(s)$, in the field of rational functions over $\reals$.
The remaining positive-real property---$\re{Z_{jk}(s)} \geq 0$ whenever $\re{s} \geq 0$---will be found to be a direct consequence of Proposition \ref{prop:zcone}.

One readily sees that individual branch admittances $y_{\alpha}(s)$ having the general form (\ref{eq:gcrl}) automatically have the positive-real property---that is to say, setting $s = \sigma + j\omega$ leads to $y_{\alpha}(s) = G(\sigma,\omega) + jB(\sigma,\omega)$ with $G \geq 0$ whenever $\sigma \geq 0$ and $B = 0$ whenever $\omega = 0$.
Here $G$ and $B$ are given explicitly by (\ref{eq:gbs}).

Now observe that the positive-real property amounts to asserting that $y_{\alpha}$ for all branches $\alpha$ is confined to a particular range of phase angles: either to $\arg y_{\alpha}(s) = 0$ (when $s$ is real), or to $\arg y_{\alpha}(s) \in [-\pi, \pi]$ (when $\re{s} \geq 0$).
Under these circumstances, Proposition \ref{prop:zcone} immediately asserts that $Z_{jk}(Y(s))$ is also confined in phase to the same range: to $\arg Z_{jk} = 0$ or to $\arg Z_{jk} \in [-\pi, \pi]$ respectively.
This establishes that $Z_{jk}(Y(s))$ is positive-real.
(One may in fact infer this directly from the identity (\ref{eq:zjky}) and addititivity of positive-real functions, without resorting to the invariance of cones in $\cplx$.)

The positive-real property was first elaborated by Wilhelm Cauer, who demonstrated that the impedance functions of passive electrical networks (those made up of resistors, capacitors and inductors) necessarily have this property.
His student Otto Brune demonstrated the converse---that any rational positive-real function $Z(s)$ with real coefficients could be realised as the impedance function (or the admittance function) of a passive network---thereby paving the way for the general theory of passive network synthesis.

In a complex analytic sense, to say $f(s)$ is positive-real asserts that both $f(s)$ and $1/f(s)$ map the real line into itself, and the half plane $\re{s} \geq 0$---alternatively, the right hemisphere of the Riemann sphere---into itself.
With the help of continuity arguments (for example, that an analytic function from the region $D \cup \partial D$ into itself, where $D$ is a domain and $\partial D$ its boundary, will map points in $D$ to points in $D$ rather than $\partial D$) and the maximum principles for analytic and harmonic functions, a number of other important consequences follow, gathered into the following proposition.
\begin{proposition}
\label{prop:poles}
Let a network have branch admittances of the form (\ref{eq:gcrl}) with $s \in \cplx$.
Then for any node pair $(j,k)$ the driving-point impedance function $Z_{jk}(s)$ is a rational function $B(s) / A(s)$ that has the positive-real property and the following properties as a consequence:
\begin{enumerate}
\item
$Z_{jk}(s)$ and $1 / Z_{jk}(s)$ are analytic (and finite) on the open half plane $\re{s} > 0$.
\item
$A(s)$ and $B(s)$ are both polynomials in $s$ with positive real coefficients, all of whose roots $s_k$ satisfy $\re{s_k} \leq 0$ and, if not real, occur as complex conjugate pairs.
(Thus, $Z_{jk}(s)$ and $1 / Z_{jk}(s)$ are at least weakly Hurwitz stable and minimum phase.)
\item
The degrees of polynomials $A(s)$ and $B(s)$ differ by at most 1.
(Equivalently, both $Z_{jk}(s)$ and $1 / Z_{jk}(s)$ are $O(s)$, $O(1)$ or $O(s^{-1})$ as $|s| \goesto \infty$, in particular in the half plane $\re{s} \geq 0$.
In terms of the extended complex plane, $Z_{jk}(s)$ has at most one pole or zero at infinity.)
\item
Any roots of $A(s)$ or $B(s)$ on the imaginary axis are simple.
\item
Roots of $A(s)$ (respectively $B(s)$) on the imaginary axis are poles of $Z_{jk}(s)$ (respectively $1 / Z_{jk}(s)$) with real strictly positive residues.
\item
The resistance $R_{jk}(s) = \re{Z_{jk}(s)}$ and reactance $X_{jk}(s) = \im{Z_{jk}(s)}$ are harmonic functions in the region $\re{s} \geq 0$, whose local maxima and minima occur on the imaginary axis.
\item
$Z_{jk}(s)$ and $1 / Z_{jk}(s)$ both satisfy a contractive angle condition: $|\arg Z_{jk}(s)| \leq |\arg s|$ when $|\arg s| \leq \pi/2$.
\end{enumerate}
\end{proposition}
To see why a positive-real function $Z(s)$ has no poles or zeros with $\re{s} > 0$---hence both $Z$ and $1/Z$ are analytic on the right half plane---it is best to appeal to the extended complex plane, or rather its compact representation as the Riemann sphere.
Let $H$ denote the right hemisphere $\re{s} > 0$ on the Riemann sphere, $\partial H$ its boundary (the union of the imaginary axis with the point $\infty$) and $\bar{H} = H \cup \partial H$ the closed right hemisphere.
By definition, $Z$ is an analytic function that maps $\bar{H}$ into itself, and by continuity, points within the interior $H$ will not be mapped to points on $\partial H$, including the points 0 and $\infty$.
This suffices to show that $Z$ (and likewise $1/Z$) has no poles or zeros in $H$.

Turning to the nature of the poles and zeros themselves, note that Proposition \ref{prop:poles}'s guarantee of simple poles with positive residues for $Z$ and $1/Z$ applies only to poles on the imaginary axis.
One readily verifies, for example, that if $a$, $b$ are real parameters with $a > b > 0$ but otherwise arbitrary, the functions
\begin{equation}
f_{a,b}(s) = \frac{a}{s + a} - \frac{b}{s + b} = \frac{(a - b) s}{s^2 + (a + b) s + a b}
\qquad \text{and} \qquad
g_a(s) = \frac{s}{(s + a)^2}
\label{eq:posrealnegrep}
\end{equation}
are positive-real and can be realised as driving-point impedances of RLC networks yet have, respectively, a negative residue at the pole $s = -b$, and a double pole at $s = -a$.

Suppose however that the positive-real function $Z(s)$ has a pole of order $N$ on the imaginary axis, at $s^* = j\omega$ for $\omega \in \reals$.
Then $Z$ can be represented by a Laurent series in a neighbourhood of $s^*$, within some positive radius of convergence $r$:
\begin{equation}
Z(s) = P(s - j\omega) + \sum_{k=1}^N \frac{b_k}{\paren{s - j\omega}^k}, \qquad 0 < |s - j\omega| < r,
\label{eq:posreallaurent}
\end{equation}
where $P(z)$ is a polynomial and $b_k \in \cplx$.
Now suppose $s$ is a point in $H$ near $s^*$, so can be represented as $s = j\omega + \epsilon\,\e^{j \theta}$ where $0 < \epsilon < r$ and $|\theta| \leq \pi/2$, and the Laurent series becomes
\begin{equation}
Z\!\paren{j\omega + \epsilon\,\e^{j \theta}}
   = P\!\paren{\epsilon\,\e^{j \theta}} + \sum_{k=1}^N b_k \e^{-j k \theta} \paren{\frac{1}{\epsilon}}^k.
\label{eq:posreallaureps}
\end{equation}
For $\epsilon$ sufficiently small, the term with $k = N$ will dominate the others in both real and absolute value; thus setting $b_N = |b_N| \e^{j \beta}$ one has
\begin{equation}
\re{Z} = |b_N| \cos\paren{\beta - N \theta} \paren{\frac{1}{\epsilon}}^N
   + O\!\brak{\paren{\frac{1}{\epsilon}}^{N-1}}, \qquad \epsilon \goesto 0.
\label{eq:posreallaureal}
\end{equation}
Since $Z$ is positive-real, this expression must remain positive for all $\theta \in [-\pi/2, \pi/2]$, which by inspection can only occur if $N = 1$ and $\beta = 0$.
This confirms that any pole on the imaginary axis must be simple, and its residue $b_1$ must be real and positive.

The above suffices to establish all properties of positive-real functions stated in Proposition \ref{prop:poles}, aside from the contractive angle condition $|\arg Z(s)| \leq |\arg s|$ when $|\arg s| \leq \pi/2$ and $Z(s)$ is positive-real.
Observe that the condition $|\arg s| \leq \pi/2$ is equivalent to saying $\re{s} \geq 0$, so the basic positive-real condition already ensures $|\arg Z(s)| \leq \pi/2$ under these conditions.
For the stronger result, consider the M\"{o}bius transformation $z = (s - a) / (s + a)$, where $a$ is any real, positive constant.
This maps the right half-plane $\re{s} > 0$ to the interior of the unit circle in the $z$-plane, the imaginary axis to the unit circle itself, and the point $s = a$ to the origin.
Also, lines with $\arg s = \theta = \text{constant}$ in the right half $s$-plane are mapped to circular arcs in the unit circle in the $z$-plane, connecting the points $z = -1$ (image of $s = 0$) and $z = 1$ (image of $s = \infty$) via the point $z = j \tan(\theta / 2)$ on the imaginary $z$-axis.
(For $\theta = 0$ this is the line segment $[-1,1]$ while for $\theta = \pm\pi/2$ it is the unit circle itself.)
Importantly, this latter property is independent of the value of $a$.

Now apply this transformation both to $s$ and to a new function $W = (Z - Z(a)) / (Z + Z(a))$, where $a$ is a positive constant, and $Z(s)$ being positive real guarantees that $Z(a)$ is also real and strictly positive.
This implicitly defines a function $W_a(z)$ that is, in a complex-analytic sense, \emph{conjugate} to the original function $Z(s)$---albeit depending on the choice of $a$ as well as on $Z(s)$.
One may verify that $W_a(z)$, for any choice of $a$, is analytic in the unit circle $|z| < 1$ and satisfies $W_a(0) = 0$.
Further, the condition on $Z(s)$ to be positive-real translates to the condition that $|W_a(z)| \leq 1$ whenever $|z| < 1$.
Accordingly one has a parametric family of functions $W_a(z)$ defined from $Z(s)$ that, for any $a > 0$, satisfies all the conditions for Schwarz's lemma (the maximum modulus principle applied to $W_a(z) / z$), and therefore has the property that $|W_a(z)| \leq |z|$ whenever $|z| < 1$.

Now fix a value $s$ with $|\arg s| = \theta < \pi/2$ (hence $\re{s} > 0$) and let $a = |s| > 0$.
Under the M\"{o}bius transformation with this value of $a$, the image of $s$ in the $z$-plane is one of the two points $z_{\pm} = \pm j \tan(\theta/2)$.
Further, under the function $W_a$ for this choice of $a$, one has $|W_a(z)| \leq |z_+| = \tan(\theta/2)$, ensuring in particular that the two points $W_a(z_{\pm})$ lie within the (closed) region bounded by the two circular arcs with endpoints $1$, $-1$ and passing through the two points $z_{\pm}$.
This suffices to confirm that the relevant preimage $Z(s)$ of $W_a(z_{\pm})$ satisfies $|\arg Z(s)| \leq \theta$.
Since $s$ is arbitrary, this establishes the angle condition whenever $|\arg s| < \pi/2$; and given $|\arg Z(s)| \leq \pi/2$ due to $Z(s)$ being positive real, this establishes the condition for all $s$ with $|\arg s| \leq \pi/2$ as required.

It may be noted that the angle condition $|\arg Z(s)| \leq |\arg s|$ for $|\arg s| \leq \pi/2$ is also \emph{sufficient} for $Z(s)$ to be positive-real---and, therefore, can serve as an equivalent (if somewhat arcane) definition of positive-real functions.
Indeed, taking $s$ to be real (hence $\arg s = 0$) immediately forces $Z(s)$ to be real (indeed positive), while taking $\re{s} \geq 0$ is equivalent to asserting $|\arg s| \leq \pi/2$ which forces $|\arg Z(s)| \leq |\arg s| \leq \pi/2$, so that $\re{Z(s)} \geq 0$ also.

\subsection{Reactance functions and strictly positive-real functions}
\label{sec:acspecial}

Within the class of positive-real impedance functions that can appear as driving-point impedances $Z_{jk}(s)$, the reactance functions form an important subclass.
These represent lossless networks built from capacitors and inductors only.
Analytically, reactance functions can be characterised as those rational functions $f(s)$ which in addition to being positive-real, also map the \emph{imaginary} axis into itself---that is, satisfy $\re{f(j\omega)} = 0$ for all $\omega \in \reals$.
(Reactance functions will however generally still have nonzero real parts for general $s = \sigma + j\omega$ with $\sigma \neq 0$.)
As a consequence of this property, reactance functions $f(s) = B(s) / A(s)$ necessarily have all their poles and zeros (roots of $A(s)$ and $B(s)$ respectively) on the imaginary axis, and these are all simple with real positive residues as a consequence of the positive-real property.

The monotonicity of reactances in lossless networks is a corollary of Proposition \ref{prop:poles}.
Indeed, this proposition implies that when $Z_{jk}(s) = B(s) / A(s)$ is a reactance function, with $n$ the number of poles (degree of $A(s)$, possibly zero), it will necessarily have a partial fraction expansion of the form
\begin{equation}
Z_{jk}(s) = \frac{B(s)}{A(s)} = a_0 s + j \omega_0 + \sum_{k=1}^n \frac{a_k}{s + j \omega_k}, \qquad
a_0 \in \reals_{\geq 0}, \quad a_k \in \reals_{>0}, \quad \omega_0, \omega_k \in \reals.
\label{eq:partfrac}
\end{equation}
On the imaginary axis this evaluates to $Z_{jk}(j \omega) = jX_{jk}(\omega)$, where $X_{jk}(\omega)$ is the real function
\begin{equation}
X_{jk}(\omega) = a_0 \omega + \omega_0 - \sum_{k=1}^n \frac{a_k}{\omega + \omega_k}, \qquad
a_0 \geq 0, \quad a_k > 0
\label{eq:partfracx}
\end{equation}
with derivative
\begin{equation}
\frac{d}{d\omega} X_{jk}(\omega) = a_0 + \sum_{k=1}^n \frac{a_k}{(\omega + \omega_k)^2}, \qquad
a_0 \geq 0, \quad a_k > 0.
\label{eq:partfracxd}
\end{equation}
This derivative exists for all $\omega$ bar the poles, $\omega = -\omega_k$.
Provided it exists, its value is always strictly positive: either $Z_{jk}(s)$ has at least one pole, in which case positivity is guaranteed by the positive residue $a_1$, or else $Z_{jk}(s)$ is a polynomial of degree 1 (representing an inductor), with strictly positive coefficient $a_0$.
Accordingly one has the following, originally due to Foster:
\begin{proposition}[Foster's reactance theorem]
\label{prop:reactance}
Let $Z_{jk}(s)$ be a driving-point impedance of a network as in Proposition \ref{prop:poles} and suppose further that $Z_{jk}(s)$ is a reactance function: that is, $Z_{jk}(j \omega) = j X(\omega)$ where $X(\omega)$ is real-valued.
Then $X(\omega)$ has a finite number (possibly zero) of simple poles $\omega_k$, corresponding to poles $s_k = j\omega_k$ of $Z_{jk}(s)$.
If $\omega_k$ is a pole then so is $-\omega_k$.
For all $\omega \in \reals$ excluding poles, $X(\omega)$ is continuous and strictly increasing, and consequently has a unique zero $\bar{\omega}_k$ in between any two adjacent poles $\omega_k$ and $\omega_{k+1}$.
Poles and zeros of $X(\omega)$ thus interleave: $\omega_1 < \bar{\omega}_1 < \omega_2 < \bar{\omega}_2 < \ldots < \omega_n$.
Lastly, all the above properties of $X(\omega)$ are shared by the susceptance $B(\omega) = - 1 / X(\omega)$, whose poles are the zeros $\bar{\omega}_k$ of $X(\omega)$ and correspond to poles of $1 / Z_{jk}(s)$.
\end{proposition}

In contrast to the reactance functions are the \emph{strictly positive-real} or \emph{strictly passive} impedance functions---those which have \emph{no} poles or zeros on the imaginary axis (hence have all poles and zeros in the open half-plane $\re{s} < 0$), other than possibly a pole or zero at infinity.
This is the subclass of functions $f(s)$ for which the condition $\re{f(s)} \geq 0$ (other than at poles) in Proposition \ref{prop:posreal} can be replaced with the stricter condition $\re{f(s)} \in (0, \infty)$ on the (non-extended) half-plane $\re{s} \geq 0$.
(Note it was previously established that when $f(s)$ is positive-real, $\re{f(s)}$ is finite with no zeros in the region $\re{s} > 0$.)

As seen previously for positive-real functions, the stricter definition entails that the sum of two strictly positive-real functions is strictly positive-real, as is the composition of two strictly positive-real functions.
Indeed one can say more, as the following (building on previous results) sums up:
\begin{proposition}
\label{prop:posrealfg}
Let $f(s)$ and $g(s)$ be rational functions satisfying the positive-real property in Proposition \ref{prop:posreal}.
Then
\begin{itemise}
\item
$\re{f(s)} > 0$ whenever $\re{s} > 0$, and indeed this condition is equivalent to the condition $\re{f(s)} \geq 0$ whenever $\re{s} \geq 0$ of Proposition \ref{prop:posreal}.
The same is true of $g(s)$.
\item
The functions $f(s) + g(s)$ and $g(f(s))$ are positive-real.
\item
If $f(s)$ is strictly positive-real then the function $g(f(s))$ is strictly positive-real.
\item
If $f(s) = B(s) / A(s)$ is strictly positive-real with $\deg B = \deg A$, then the function $f(g(s))$ is strictly positive-real.
\item
If $f(s)$ is strictly positive-real and $g(s)$ has no poles on the imaginary axis, then the functions $f(s) + g(s)$ and $f(g(s))$ are strictly positive-real.
\end{itemise}
\end{proposition}
Most of Proposition \ref{prop:posrealfg} restates or follows from previous results.
If $f(s)$ is positive-real, then by Proposition \ref{prop:poles} $\re{f(s)}$ cannot take a zero value in the region $\re{s} > 0$---only on its boundary, the imaginary axis.
It immediately follows that $\re{f(s)} > 0$ whenever $\re{s} > 0$.
Conversely, suppose $\re{f(s)}$ is strictly positive on the entire open half-plane $\re{s} > 0$: then to show that $\re{f(s)} \geq 0$ (excluding poles) when $\re{s} \geq 0$, it suffices to show that all points on the imaginary axis, other than poles, satisfy $\re{f(s)} \geq 0$.
But this follows immediately from the fact that $f(s)$ is a rational function, hence $\re{f(s)}$ is continuous.

If $f(s)$ is strictly positive-real, then so are $f(s) + g(s)$ and $f(g(s))$ for any positive-real $g$ essentially by definition, \emph{except} at points $\re{s} \geq 0$ where $g(s)$ has a pole.
As previously seen, such a pole is necessarily on the imaginary axis, but if present, may cause the composite function to take the value 0 or $\infty$ and thus fail to satisfy the strict positive-real property.
By working in the extended complex plane, one may `remove' the pole in the composite $f(g(s))$ when $f(s) = B(s) / A(s)$ with $\deg B = \deg A$, since then as $g(s) \goesto \infty$ one has $f(g(s)) \goesto b_n / a_n$ where $b_n$, $a_n$ are the leading coefficients of $B$ and $A$ respectively, and are necessarily positive.
But otherwise, one needs the additional condition ruling out poles of $g$ when forming the composites $f(s) + g(s)$ and $f(g(s))$.

The composiite $g(f(s))$, on the other hand, is strictly positive-real for any positive-real $g$ without restriction, since for $\re{s} \geq 0$ the strictly positive-real $f$ generates an argument $s' = f(s)$ for $g$ such that $\re{s'} \in (0,\infty)$.
The positive-real function $g$ has no poles with $\re{s} > 0$ and by the argument above satisfies $\re{g(s')} > 0$, ensuring the composite $g(f(s))$ is strictly positive-real.

Using the above results one obtains necessary and sufficient conditions for a network impedance function $Z_{jk}(s)$ to be strictly positive-real when branch admittances take the form (\ref{eq:gcrl}).
First, observe that an admittance or impedance of the form (\ref{eq:gcrl}) is strictly positive-real if and only if $g > 0$ and $r > 0$.
Now, let network $\nw$ be driven by a single source connected between nodes $j$ and $k$ with current $I \neq 0$ and consider a general term $y_{\alpha} |v_{\alpha}|^2$ in the formula (\ref{eq:tellegenz}).
This is a positive-real function of $s$, that is strictly positive-real provided $y_{\alpha}(s)$ is strictly positive-real and $v_{\alpha}(s)$ is finite and nonzero whenever $\re{s} \geq 0$.
Letting $p,q$ denote the endpoint nodes of branch $\alpha$, one sees that $v_{\alpha}(s) = \tz{pq}{jk}(s) \cdot I$ is finite and nonzero precisely when $\tz{pq}{jk}$ is: that is, when the node pairs $(j,k)$ and $(p,q)$ are not conjugate in $\nw$ and $s$ is not a pole or zero of the transfer impedance $\tz{pq}{jk}(s)$.
Provided any one term satisfies these criteria, and none of the $\tz{pq}{jk}$ have poles for $\re{s} \geq 0$ (which are necessarily on the imaginary axis), the sum of these terms will be strictly positive-real.
This leads to the following statement:
\begin{proposition}
\label{prop:sposreal}
Let $\nw$ be a network whose branch admittances are of the form (\ref{eq:gcrl}) and let $j$ and $k$ be any two nodes.
The driving-point impedance function $Z_{jk}(s)$ is strictly positive-real if and only if the admittance matrix $Y(s)$ satisfies the following criteria:
\begin{itemise}
\item
For all branches $\alpha$ the transfer impedance $\tz{pq}{jk}(s)$ has no poles on the imaginary axis, where $(p,q)$ are the node indices of the endpoints of branch $\alpha$.
\item
There exists at least one branch $\alpha$, with endpoint nodes $(p,q)$, such that the pairs $(p,q)$ and $(j,k)$ are not conjugate to one another, the transfer impedance $\tz{pq}{jk}$ has no zeros on the imaginary axis, and the branch admittance $y_{\alpha}(s)$ is of the form (\ref{eq:gcrl}) with $g > 0$ and $r > 0$.
\end{itemise}
\end{proposition}
In order to apply Proposition \ref{prop:sposreal}, it would be useful to have some general conditions under which a collection of network transfer impedances has, between them, no poles on the imaginary axis.
The following, by appealing to Kirchhoff characteristics, provides a sufficient condition.
\begin{proposition}
\label{prop:cpoly}
Let $\nw$ be a network all of whose branch admittances are of the form (\ref{eq:gcrl}).
For any branch $\alpha$ in $\nw$, denote the branch admittance as $y_{\alpha}(s) = p_{\alpha}(s) / q_{\alpha}(s)$ where $p_{\alpha}$ and $q_{\alpha}$ are linear or constant polynomials, and define the polynomial
\begin{equation}
C(s) = \sum_{\tr \subseteq \nw}
   \prod_{\alpha \in \tr} p_{\alpha}(s) \prod_{\beta \in \nw \setminus \tr} q_{\beta}(s),
\label{eq:cpoly}
\end{equation}
where the sum is over all spanning trees $\tr$ in $\nw$ as in (\ref{eq:kirchhoff}), and the second product is precisely over those branches $\beta$ not contained in $\tr$ (also known as the cotree).
Then $C(s)$ has no roots in the region $\re{s} > 0$ in the complex plane.
If in addition $C(s)$ has no roots with $\re{s} = 0$, then no transfer impedance $\tz{pq}{jk}(s)$ of $\nw$ has a pole on the imaginary axis.
\end{proposition}
Proposition \ref{prop:cpoly} follows from the definition of $\tz{pq}{jk}$ as a ratio of cofactors $C_{pq,jk}(Y) / c(Y)$ of the admittance matrix $Y$, and the identity of $c(Y)$ with the Kirchhoff characteristic $\kappa(\nw)$.
Any cofactor of $Y$, by construction, is expressible as a sum-of-products form on the admittances $y_{\alpha}$, and accordingly a rational function of $s$.
It can further be seen by construction that the ratio $C_{pq,jk}(Y) / c(Y)$ can itself be `regularised' as a rational function---that is, clear fractions in numerator and denominator---by multiplying above and below by the product
\begin{equation}
Q(s) = \prod_{\alpha \in \nw} q_{\alpha}(s).
\label{eq:qpoly}
\end{equation}
The polynomial $C(s)$ in (\ref{eq:cpoly}) is simply the product $Q(s) \kappa(\nw)$.
From all the above it follows that the poles of the rational function $\tz{pq}{jk}(s)$ will be a subset of the roots of $C(s)$---perhaps a proper subset if any pole-zero cancellation occurs.
Thus, as long as no roots of $C(s)$ are on the imaginary axis, neither will be the poles of $\tz{pq}{jk}$.

That no roots of $C(s)$ have positive real part is a corollary of the following formula for a driving-point impedance $Z_{jk}(s)$.
\begin{proposition}
\label{prop:zpoly}
Let $\nw$ be a network all of whose branch admittances are of the form (\ref{eq:gcrl}) and let $p_{\alpha}(s)$, $q_{\alpha}(s)$ and $C(s)$ be defined as in Proposition \ref{prop:cpoly}.
Then for any node pair $j$ and $k$, the driving-point impedance has the expression
\begin{equation}
Z_{jk}(s) = \frac{1}{C(s)} \sum_{\tr_j, \tr_k \subseteq \nw} \prod_{\alpha \in \tr_j} p_{\alpha}(s)
   \prod_{\beta \in \tr_k} p_{\beta}(s) \prod_{\gamma \in \nw \setminus \tr_j \cup \tr_k} q_{\gamma}(s)
\label{eq:zpoly}
\end{equation}
where the sum is over pairs of disjoint trees $\tr_j$, $\tr_k$ with $j \in \tr_j$, $k \in \tr_k$ and $\tr_j \cup \tr_k$ covering all nodes of $\nw$ as in Proposition \ref{prop:kcon}.
Pole-zero cancellation may occur in the rational function defined by (\ref{eq:zpoly}).
If nodes $j$ and $k$ are themselves connected by a branch $\delta$, then $q_{\delta}$ appears in every term of the sum (\ref{eq:zpoly}) and consequently $q_{\delta}(s)$ divides $Z_{jk}(s)$.
\end{proposition}
Since $Z_{jk}(s)$ is known to be a positive-real function, the right hand side of (\ref{eq:zpoly}) is necessarily a rational function whose poles and zeros satisfy $\re{s} \leq 0$.
This does not quite necessarily imply that all roots of $C(s)$ satisfy $\re{s} \leq 0$ as asserted in Proposition \ref{prop:cpoly}, due to the possibility of pole-zero cancellation.
However, one may bridge this technical gap by appealing to the fact that the roots of a polynomial are continuous functions of its coefficients; hence if pole-zero cancellation does occur in (\ref{eq:zpoly}), one may perturb admittance parameters by a small amount $\epsilon$ without violating the constraints in (\ref{eq:gcrl}), resulting in a family of perturbed impedance functions $Z_{jk}(s,\epsilon)$ all of which remain positive-real and where pole-zero cancellation does not occur except at $\epsilon = 0$.
The corresponding polynomials $C(s,\epsilon)$ for $\epsilon > 0$ necessarily have all roots satisfying $\re{s} \leq 0$, and by continuity this also extends to $C(s,0) = C(s)$.

\section{Sources in Networks}
\label{sec:sources}

This section provides more explicit consideration of sources in linear networks described by admittance matrices.
To this point, the only sources explicitly considered have been independent current sources connected between pairs of nodes, as these readily fit into the framework via the right-hand side of the basic network equation $Y v = i$.
While most scenarios can be adapted to this framework, it is helpful to consider other types of source more explicitly; in particular independent voltage sources.
Although network solutions that fix a \emph{single} node voltage have been considered, fixing a constant difference between \emph{two} node voltages is a more delicate proposition.

\subsection{Norton and Th\'{e}venin sources, and equivalence of subnetworks}
\label{sec:thevnort}

The most straightforward extension of the framework based on $Y$ to explicitly consider voltage sources is provided by the classical `Th\'{e}venin--Norton' equivalence.
As in previous discussions, let a source with fixed current $I$ be connected between network nodes $p$ and $q$, with the current sensed positive from $q$ to $p$.
Suppose further that a nonzero admittance $y$ also appears between the same nodes $p$ and $q$.
Then, writing out rows $p$ and $q$ of the network equations $Y v = i$ explicitly, one has
\begin{equation}
\begin{split}
\text{(other branch currents)} + y \paren{v_p - v_q} &= +I + \text{(other source currents)}, \\
\text{(other branch currents)} + y \paren{v_q - v_p} &= -I + \text{(other source currents)}.
\end{split}
\label{eq:norton}
\end{equation}
Observe that both these equations express the quantity
\begin{equation}
N_{pq}(I, y) = y \paren{v_p - v_q} - I
\label{eq:nortpq}
\end{equation}
in terms of other quantities independent of the source current $I$ and admittance $y$.
The parallel combination of the source $I$ and the admittance $y$ between two designated nodes is known as a \emph{Norton source}.
It is characterised by the expression (\ref{eq:nortpq}), which describes the net current extracted from the network at node $p$ and injected at node $q$ due to the presence of the Norton source.
Setting $N_{pq}$ formally to zero has the effect of replacing the source (and its admittance) with an open circuit, while setting $I = 0$ leaves just the admittance $y$.

Provided again that $y$ is nonzero, the expression (\ref{eq:nortpq}) can also be written as
\begin{equation}
N_{pq}(I, y) = y \paren{v_p - v_q - y^{-1} I}.
\label{eq:nortpq1}
\end{equation}
This form effectively builds the source $I$ into the voltage drop seen by the admittance $y$.
More explicitly, define an alternative source current as
\begin{equation}
T_{pq}(V, y) = y \paren{v_p - \brak{v_q + V}}.
\label{eq:thevpq}
\end{equation}
From first principles, this is the current in an admittance element $y$ that sees at one terminal the voltage $v_p$ at node $p$ and at the other the voltage $v_q + V$: that at node $q$ offset by a fixed voltage $V$.
In other words, it expresses the current flow in a compound element comprising the fixed voltage source $V$ (sensed as a rise from $q$ to $p$) in series with the admittance $y$.
This series combination is known as a \emph{Th\'{e}venin source}.
As with the Norton source, setting $T_{pq}$ formally to zero replaces the source with an open circuit, and setting $V = 0$ leaves the admittance $y$ connected between $p$ and $q$.

Comparing (\ref{eq:thevpq}) with (\ref{eq:nortpq1}) and (\ref{eq:nortpq}), it is immediately evident that
\begin{equation}
N_{pq}(I, y) = T_{pq}(y^{-1} I, y) \qquad \text{and} \qquad T_{pq}(V, y) = N_{pq}(y V, y).
\label{eq:thevnort}
\end{equation}
This expresses in pure algebraic form the celebrated \emph{Th\'{e}venin--Norton equivalence}, as follows:
\begin{proposition}[Th\'{e}venin--Norton]
\label{prop:thevnort}
Suppose a network contains a Norton source with current $I$ and parallel admittance $y \neq 0$ connected between some pair of nodes $p$ and $q$.
Then this Norton source may be replaced with a Th\'{e}venin source connected between nodes $p$ and $q$, having source voltage $V = y^{-1} I$ and series admittance $y$, without affecting the behaviour of the network.
Similarly, if a network contains a Th\'{e}venin source with voltage $V$ and series admittance $y$, it may be replaced with a Norton source having current $I = yV$ and parallel admittance $y$ between the same pair of nodes, without affecting the behaviour of the network.
\end{proposition}
That the source substitution as in (\ref{eq:thevnort}) does not affect the network behaviour is guaranteed by the form of the network equations in (\ref{eq:norton}), given the terms and the equations not shown are all independent of the source quantities.
Indeed, as the basic principle does not even depend on the form of these unknown terms, the balance of the network does not need to be linear or passive (provided no other network quantities have a hidden dependence on the internal source quantities).

As is well known, Norton and Th\'{e}venin sources can also stand in as `one-port equivalents' for entire linear subnetworks of an electrical network.
This will be demonstrated explicitly on the assumption that any voltage sources in the network are of Th\'{e}venin type, so that the network equations may be rendered into the form $Y v = i$ using Th\'{e}venin--Norton conversions as necessary.
For present purposes, sources will be included as branches of the underlying graph of the network, so that a general network branch is either an admittance represented in the matrix $Y$, or a current source represented in the vector $i$.

A \emph{one-port decomposition} of network $\nw$ with \emph{port nodes} $p$ and $q$ is a partition of the nodes in $\nw \setminus \{p, q\}$ into two sets $A$ and $B$, which are separated from each other by the nodes $p$ and $q$.
This separation may be formally defined in two equivalent ways:
\begin{itemise}
\item
Nodes $p$ and $q$ together form a \emph{node cutset} in the graph of $\nw$, meaning their removal along with the incident branches leaves behind connected components (`islands') such that no connected component contains nodes from both $A$ and $B$.
\item
Node sets $A$ and $B$ are exclusively path-connected in the graph of $\nw$ via nodes $p$ and $q$.
That is, if $P$ is any path in $\nw$ not including either of nodes $p$ or $q$, then $P$ necessarily includes only nodes in $A$, or only nodes in $B$.
Consequently, any circuit in $\nw$ (a closed path on distinct nodes, equivalently a degree--2 regular subnetwork of $\nw$) that includes nodes in both $A$ and $B$ necessarily includes both nodes $p$ and $q$.
\end{itemise}
Note that the formal definitions above technically admit trivial cases where sets $A$ or $B$ are empty; such trivial cases exist for any choice of $p$ and $q$ in any network with at least two nodes.
(If $\nw$ has exactly two nodes, the definition is satisfied with $A$ and $B$ both empty.)
These trivial cases are logically admissible and consistent, but generally of little practical interest.

A one-port decomposition is the simplest way in which a network can be notionally `split' into subnetworks that are capable of influencing each other's behaviour.
This mutual influence is not possible with disconnected subnetworks, nor is it possible with subnetworks separated by a \emph{single} common node.
(To see the latter, consider a solution of the network grounded at the common node; the structure of the grounded network equations will then separate into a block diagonal form with the subnetwork solutions being independent of one another.
As a consequence of Proposition \ref{prop:grounded}, the subnetwork solutions will remain independent in the general case.)

With any one-port decomposition as above---denoted as $\nw = A \cup \{p,q\} \cup B$ with a slight abuse of notation---one may associate \emph{open-circuit subnetworks} $\nw_A$ and $\nw_B$ in a natural way as follows:
\begin{itemise}
\item
Subnetwork $\nw_A$ is derived from $\nw$ by deleting all nodes in $B$ and incident branches.
\item
Subnetwork $\nw_B$ is derived from $\nw$ by deleting all nodes in $A$ and incident branches.
\end{itemise}
It follows from the nature of the decomposition that forming network $\nw_A$ has no effect on the branches of the `A network' (the subnetwork on nodes $A \cup \{p,q\}$).
If the set $B$ is empty, then $\nw_A$ and $\nw$ are identical, while if $A$ is empty, then $\nw_A$ consists only of nodes $p$ and $q$ and any branch or source connecting them in $\nw$.
Similar observations apply to the network $\nw_B$.

The one-port decomposition and open-circuit subnetworks can be related directly to a block structure on the admittance matrix $Y$.
To describe this, it will be convenient to relabel the nodes of $\nw$ so that $q = p + 1$, the nodes in $A$ are numbered 1 through $p - 1$, and nodes in $B$ are numbered $q + 1$ through $n$, without loss of generality.
Then the basic network equation $Y v = i$ can be decomposed into the following form:
\begin{equation}
\mat{\qquad\vblock Y_{AA} \qquad & -Y_{Ap} & -Y_{Aq} & 0 \\
   -Y_{Ap}^T & Y_{pp} & -y_{pq} & -Y_{Bp}^T \\
   -Y_{Aq}^T & -y_{pq} & Y_{qq} & -Y_{Bq}^T \\
   0 & -Y_{Bp} & -Y_{Bq} & \qquad\vblock Y_{BB} \qquad}
\mat{\vblock v_A \\ v_p \\ v_q \\ \vblock v_B}
   = \mat{\vblock i_{AA} + i_{Ap} + i_{Aq} \\ -\Sigma i_{Ap} - \Sigma i_{Bp} + I_{pq} \\
      -\Sigma i_{Aq} - \Sigma i_{Bq} - I_{pq} \\ \vblock i_{BB} + i_{Bp} + i_{Bq}}.
\label{eq:oneport}
\end{equation}
To explain further the nomenclature in the schema (\ref{eq:oneport}):
\begin{itemise}
\item
$Y_{AA}$ and $Y_{BB}$ are submatrices (indeed principal minors) of $Y$ corresponding to nodes in $A$ and $B$ respectively.
Likewise $v_A$ and $v_B$ are the corresponding node voltage vectors.
\item
$Y_{Ap}$ is a column vector made up of admittances connecting each node in $A$ with the port node $p$ (or zero if there is no connection to $p$ from the relevant node).
Column vectors $Y_{Aq}$, $Y_{Bp}$ and $Y_{Bq}$ are defined similarly.
\item
Vector $i_{AA}$ represents current sources connected between node pairs in $A$, and therefore sums to zero.
Likewise $i_{BB}$ represents current sources connected between node pairs in $B$.
\item
Vector $i_{Ap}$ represents current sources connected between the port node $p$ and each node in $A$, with zeros where no source exists.
Vectors $i_{Aq}$, $i_{Bp}$ and $i_{Bq}$ are defined similarly.
(By assumption, there are no sources connected between a node in $A$ and one in $B$.)
The notation $\Sigma v$, where $v$ is any vector, is shorthand for $1^T v$, the algebraic sum of all elements in $v$.
\item
$y_{pq}$ is the branch admittance connecting nodes $p$ and $q$, or zero if no such admittance is present.
Similarly $I_{pq}$ if nonzero represents a current source connected between nodes $p$ and $q$ in $\nw$ and sensed positive from $q$ to $p$ (analogous to the current $I$ in (\ref{eq:norton})).
\item
With the $\Sigma$ notation above, the diagonal elements of $Y$ for the port nodes have the expressions
\begin{equation}
Y_{pp} = \Sigma Y_{Ap} + \Sigma Y_{Bp} + y_{pq} = \sum_{k=1}^n y_{kp}, \qquad
Y_{qq} = \Sigma Y_{Aq} + \Sigma Y_{Bq} + y_{pq} = \sum_{k=1}^n y_{kq},
\label{eq:oneporty}
\end{equation}
consistent with the definition of a diagonal element of $Y$.
\end{itemise}
Using the scheme (\ref{eq:oneport}) for the network equations, one may readily derive the corresponding equations for the networks $\nw_A$ and $\nw_B$.
Consider the network $\nw_B$: to form this, one deletes the rows of (\ref{eq:oneport}) corresponding to nodes in $A$, and the corresponding columns of $Y$; one also deletes the terms $\Sigma Y_{Ap}$ and $\Sigma Y_{Aq}$ from $Y_{pp}$ and $Y_{qq}$ respectively, and the terms $\Sigma i_{Ap}$ and $\Sigma i_{Aq}$ from the vector of currents.
As a technicality, the row and column indices will also change so that index $p$ in the equations becomes index 1, and index $q$ becomes index 2.
The below presentation of the resulting equations takes advantage of this by labelling the new port voltages $v_1$ and $v_2$ to distinguish them from the voltages $v_p$ and $v_q$ in the original network, since these will not necessarily take the same values.
Likewise, the notation $B'$ distinguishes node voltages $v_{B'}$ in $\nw_B$ from their values $v_B$ in (\ref{eq:oneport}).
\begin{equation}
\mat{\Sigma Y_{Bp} + y_{pq} & -y_{pq} & -Y_{Bp}^T \\
   -y_{pq} & \Sigma Y_{Bq} + y_{pq} & -Y_{Bq}^T \\
   -Y_{Bp} & -Y_{Bq} & \qquad\vblock Y_{BB} \qquad}
\mat{v_1 \\ v_2 \\ \vblock v_{B'}}
   = \mat{-\Sigma i_{Bp} + I_{pq} \\ -\Sigma i_{Bq} - I_{pq} \\ \vblock i_{BB} + i_{Bp} + i_{Bq}}.
\label{eq:oneportb}
\end{equation}
The equations for network $\nw_A$ are obtained analogously.
Note that (\ref{eq:oneportb}) adheres to the basic form $Y v = i$ where $Y$ and $i$ retain all the properties of an admittance matrix and source current vector respectively.
Assuming the original network $\nw$ to be connected, the same holds for $\nw_A$ and $\nw_B$ by virtue of the one-port decomposition, and accordingly the matrix $Y$ is of maximum rank in each case.
Indeed, taking the reduction one step further and contracting either $\nw_A$ or $\nw_B$ on nodes $p$ and $q$, it is evident each network remains connected and each contracted admittance matrix of maximum rank.
This implies in turn that the submatrices $Y_{AA}$ and $Y_{BB}$ are nonsingular.

Consider now a grounded solution of (\ref{eq:oneportb}) for network $\nw_B$, with ground at node 2 (which was $q$).
Removing row and column 2 of (\ref{eq:oneportb}) as redundant, the equations for this particular solution are
\begin{equation}
\mat{\Sigma Y_{Bp} + y_{pq} & -Y_{Bp}^T \\ -Y_{Bp} & \qquad\vblock Y_{BB} \qquad}
\mat{\vgnd{2}_1 \\ \vblock \vgnd{2}_{B'}}
   = \mat{-\Sigma i_{Bp} + I_{pq} \\ \vblock i_{BB} + i_{Bp} + i_{Bq}}.
\label{eq:oneportbgnd}
\end{equation}
Recall from earlier discussion that $\vgnd{2}_1$ (for example) is the solution for $v_1$ with $v_2 = 0$, and also the value of $v_1 - v_2$ in any arbitrary solution to (\ref{eq:oneportb}).
Note also that the matrix at left of (\ref{eq:oneportbgnd}) is no longer strictly an admittance matrix but instead a minor of one, and so the equations have a unique solution.
The value of $\vgnd{2}_1$ in this unique solution is known as the \emph{open-circuit voltage} $V_B^{oc}$ for network $\nw_B$ at the identified port.

Before solving explicitly for $V_B^{oc}$, consider first the special conditions under which this voltage would take a zero value---turning the port into a \emph{short circuit}.
Setting $\vgnd{2}_1 = 0$ and recalling that $Y_{BB}$ is nonsingular by assumption, one necessarily has that
\begin{equation}
\vgnd{2}_{B'} = Y_{BB}^{-1} \paren{i_{BB} + i_{Bp} + i_{Bq}} \qquad \text{and} \qquad
-Y_{Bp}^T \vgnd{2}_{B'} = -\Sigma i_{Bp} + I_{pq},
\label{eq:oneportbsc1}
\end{equation}
or after eliminating $\vgnd{2}_{B'}$,
\begin{equation}
I_{pq} = \Sigma i_{Bp} - Y_{Bp}^T Y_{BB}^{-1} \paren{i_{BB} + i_{Bp} + i_{Bq}}.
\label{eq:oneportbsc2}
\end{equation}
In words, formula (\ref{eq:oneportbsc2}) expresses the quantity of `bias' current that, when directed from node 2 to node 1 in the grounded network $\nw_B$, makes the voltage zero at both nodes.
It then follows that, if the source $I_{pq}$ is removed and replaced with an ideal short circuit between nodes 1 and 2, the same amount of current will flow in $\nw_B$ (since the uniqueness of (\ref{eq:oneportbsc2}) means any different amount of current would result in a voltage difference between these nodes).
When reversed in sign to denote the current flow from node 1 to node 2, this is the \emph{short-circuit current} $I_B^{sc}$ for network $\nw_B$ at its port, that is
\begin{equation}
I_B^{sc} = Y_{Bp}^T Y_{BB}^{-1} \paren{i_{BB} + i_{Bp} + i_{Bq}} - \Sigma i_{Bp}.
\label{eq:oneportbsc}
\end{equation}

Returning now to (\ref{eq:oneportbgnd}), its solution will bear on that for the original network $\nw$, as seen by comparing (\ref{eq:oneportbgnd}) directly with the grounded solution of equations (\ref{eq:oneport}) with ground at node $q$:
\begin{equation}
\mat{\qquad\vblock Y_{AA} \qquad & -Y_{Ap} & 0 \\
   -Y_{Ap}^T & Y_{pp} & -Y_{Bp}^T \\
   0 & -Y_{Bp} & \qquad\vblock Y_{BB} \qquad}
\mat{\vblock \vgnd{q}_A \\ \vgnd{q}_p \\ \vblock \vgnd{q}_B}
   = \mat{\vblock i_{AA} + i_{Ap} + i_{Aq} \\ -\Sigma i_{Ap} - \Sigma i_{Bp} + I_{pq} \\
      \vblock i_{BB} + i_{Bp} + i_{Bq}}.
\label{eq:oneportgnd}
\end{equation}
Observe in particular that the third group of equations in (\ref{eq:oneportgnd})---with right hand side $i_{BB} + i_{Bp} + i_{Bq}$---exactly match a corresponding group of equations in (\ref{eq:oneportbgnd}).
With the assumed connectivity of $\nw$ guaranteeing $Y_{BB}$ is nonsingular, these two groups of equations can be written explicitly as
\begin{equation}
\vgnd{2}_{B'} = Y_{BB}^{-1} \paren{i_{BB} + i_{Bp} + i_{Bq} + Y_{Bp} \vgnd{2}_1}, \qquad
\vgnd{q}_B = Y_{BB}^{-1} \paren{i_{BB} + i_{Bp} + i_{Bq} + Y_{Bp} \vgnd{q}_p}
\label{eq:oneportbsoln}
\end{equation}
with the difference between $\vgnd{2}_{B'}$ and $\vgnd{q}_B$ depending ultimately on the difference in port voltage:
\begin{equation}
\vgnd{q}_B = \vgnd{2}_{B'} + Y_{BB}^{-1} Y_{Bp} \paren{\vgnd{q}_p - \vgnd{2}_1}.
\label{eq:oneportbdiff}
\end{equation}
This formula serves as a certificate that the act of disconnecting subnetwork A, and replacing it with a grounded voltage source with the \emph{same} potential at node $p = 1$, reproduces the same behaviour within subnetwork B.
Owing to the properties of grounded solutions, this generalises to any situation where the voltage difference across the port nodes $p$ and $q$ is maintained at its prior value.
By symmetry, it also applies equally to the behaviour of A when subnetwork B is disconnected and replaced with an equivalent voltage source.

The promised unique solution of (\ref{eq:oneportbgnd}) for $\vgnd{2}_1$, equivalently the open-circuit voltage $V_B^{oc}$, is
\begin{equation}
V_B^{oc} = \vgnd{2}_1
   = \frac{Y_{Bp}^T Y_{BB}^{-1} \paren{i_{BB} + i_{Bp} + i_{Bq}} - \Sigma i_{Bp} + I_{pq}}
      {\Sigma Y_{Bp} + y_{pq} - Y_{Bp}^T Y_{BB}^{-1} Y_{Bp}}
   = \frac{I_{pq} + I_B^{sc}}{\Sigma Y_{Bp} + y_{pq} - Y_{Bp}^T Y_{BB}^{-1} Y_{Bp}}.
\label{eq:oneportvoc}
\end{equation}
The expression at right makes clear that a short circuit condition $V_B^{oc} = 0$ may be induced using a current source $I_{pq} = -I_B^{sc}$ to enforce current flow equal to the short-circuit current, reproducing the condition (\ref{eq:oneportbsc2}).
Any other source value $I_{pq}$ (including none at all if $I_B^{sc} \neq 0$) will lead to a nonzero open-circuit voltage $V_B^{oc}$.

With $\vgnd{2}_1$ known, the first equation of (\ref{eq:oneportbsoln}) furnishes the unique solution for the B network voltages $\vgnd{2}_{B'}$.
Observe that the quantities $\vgnd{2}_1$ and $\vgnd{2}_{B'}$ are now \emph{uniquely determined} by the admittances and source currents in the subnetwork induced by the nodes $\{p, q\} \cup B$---in particular they are independent of the source currents and admittances in subnetwork A.
This independence carries over to the second equation of (\ref{eq:oneportbsoln}), which expresses the voltages $\vgnd{q}_B$ for the \emph{full} network $\nw$ in the form
\begin{equation}
\vgnd{q}_B = Y_{BB}^{-1} \paren{i_{BB} + i_{Bp} + i_{Bq}} + Y_{BB}^{-1} Y_{Bp} \vgnd{q}_p
   = v_{B0} + k_{B0} \vgnd{q}_p
\label{eq:oneportbpq}
\end{equation}
where the vector $v_{B0}$ depends only on the source currents and admittances, and the vector $k_{B0}$ only on the admittances, in the network $\nw_B$.
One may in fact go further, and note that in this case neither $v_{B0}$ nor $k_{B0}$ depend on any source $I_{pq}$ or admittance $y_{pq}$ that may be present between nodes $p$ and $q$ themselves.

Equation (\ref{eq:oneportbpq}) may now be used to eliminate the voltages $\vgnd{q}_B$ from the equations (\ref{eq:oneportgnd}) for $\nw$.
This leads to the equivalent set of equations
\begin{equation}
\mat{\qquad\vblock Y_{AA} \qquad & -Y_{Ap} \\ -Y_{Ap}^T & Y_{pp} - Y_{Bp}^T k_{B0}}
\mat{\vblock \vgnd{q}_A \\ \vgnd{q}_p}
   = \mat{\vblock i_{AA} + i_{Ap} + i_{Aq} \\ -\Sigma i_{Ap} - \Sigma i_{Bp} + I_{pq} + Y_{Bp}^T v_{B0}}.
\label{eq:oneportgndp}
\end{equation}
Now, using the identity $Y_{BB} 1_m - Y_{Bp} - Y_{Bq} = 0$ (implied by the form of equations (\ref{eq:oneport}), with $m$ the number of nodes in B) one may write
\begin{multline}
Y_{pp} - Y_{Bp}^T k_{B0}
   = \Sigma Y_{Ap} + \Sigma Y_{Bp} + y_{pq} - Y_{Bp}^T Y_{BB}^{-1} \paren{Y_{BB} 1_m - Y_{Bq}} \\
   = \Sigma Y_{Ap} + \Sigma Y_{Bp} + y_{pq} - \Sigma Y_{Bp} + Y_{Bp}^T Y_{BB}^{-1} Y_{Bq}
   = \Sigma Y_{Ap} + y_{pq} + Y_{Bp}^T Y_{BB}^{-1} Y_{Bq}
\label{eq:oneportyp}
\end{multline}
leading to the representation of (\ref{eq:oneportgndp}) as
\begin{equation}
\mat{\qquad\vblock Y_{AA} \qquad & -Y_{Ap} \\ -Y_{Ap}^T & \Sigma Y_{Ap} + y_{pq}'}
\mat{\vblock \vgnd{q}_A \\ \vgnd{q}_p}
   = \mat{\vblock i_{AA} + i_{Ap} + i_{Aq} \\ -\Sigma i_{Ap} + I_{pq}'},
\label{eq:oneportgndp2}
\end{equation}
where
\begin{equation}
y_{pq}' = y_{pq} + Y_{Bp}^T Y_{BB}^{-1} Y_{Bq} \qquad \text{and} \qquad
I_{pq}' = I_{pq} + I_B^{sc}.
\label{eq:nortonyi}
\end{equation}
Crucially, one may observe that the form of (\ref{eq:oneportgndp2}) is identical to that for a grounded solution of the network $\nw_A$---a dual form to (\ref{eq:oneportbgnd})---only with the admittance $y_{pq}$ and source current $I_{pq}$ replaced with the augmented quantities $y_{pq}'$ and $I_{pq}'$ in (\ref{eq:nortonyi}).
These quantities in turn depend only on admittances and source currents connecting to nodes in the eliminated B subnetwork.

From the above it may be concluded that the solution for $\nw$ grounded at $q$---hence the more general solution also---is unchanged (from the point of view of all nodes and branches other than in subnetwork B) if one notionally disconnects network B entirely and substitutes a Norton source between the port nodes $p$ and $q$ with the values $I_{pq}'$ and $y_{pq}'$ indicated by (\ref{eq:nortonyi}).
In addition, one observes by inspection of (\ref{eq:oneportvoc}) that
\begin{equation}
V_B^{oc} = \frac{I_{pq}'}{y_{pq}'},
\label{eq:oneportvoca}
\end{equation}
that is, the ratio of the equivalent source current and admittance in (\ref{eq:oneportgndp2}) is identical to the open-circuit voltage of $\nw_B$.
To see that the denominator of (\ref{eq:oneportvoc}) is just $y_{pq}'$, observe that using the identity above for $Y_{Bp}$ one has
\begin{equation}
\Sigma Y_{Bp} - Y_{Bp}^T Y_{BB}^{-1} Y_{Bp}
   = \Sigma Y_{Bp} - Y_{Bp}^T Y_{BB}^{-1} \paren{Y_{BB} 1_m - Y_{Bq}}
   = Y_{Bp}^T Y_{BB}^{-1} Y_{Bq}.
\label{eq:oneportypq}
\end{equation}
The values added to $I_{pq}$ and $y_{pq}$ to form $I_{pq}'$ and $y_{pq}'$ respectively in (\ref{eq:nortonyi}) thus define a \emph{Norton source equivalent} to the subnetwork $\nw_B$ in $\nw$.
The equivalence takes the form of a Norton source as it takes the form of an additional source current and parallel admittance connected between the port nodes $p$ and $q$.
By a simple application of Proposition \ref{prop:thevnort} one may also define a Th\'{e}venin source equivalent, in which case (\ref{eq:oneportvoca}) defines the Th\'{e}venin source voltage (assumed to absorb the original $I_{pq}$ and $y_{pq}$ as well, if present).
Summing up these results one has:
\begin{proposition}
\label{prop:tnequiv}
Let $\nw$ be a connected network that remains connected when sources are suppressed, and let $A \cup \{p,q\} \cup B$ be any one-port decomposition of $\nw$ into node-sets $A$ and $B$ separated by port nodes $p$ and $q$.
Let the open-circuit subnetwork $\nw_B$ have open-circuit voltage $V_B^{oc} \neq 0$ and short-circuit current $I_B^{sc}$, and let $y = I_B^{sc} / V_B^{oc}$.
Then the subnetwork $\nw_B$ can be removed and replaced either with a Th\'{e}venin source $T_{pq}(V_B^{oc}, y)$ or with a Norton source $N_{pq}(I_B^{sc},y)$ between nodes $p$ and $q$ without any change to conditions in the remainder of $\nw$.
If the admittance matrix and current injections of $\nw$ are represented as in (\ref{eq:oneport}), then with $I_{pq} = y_{pq} = 0$ the explicit formulae for the equivalent source quantities are
\begin{equation}
I_B^{sc} = Y_{Bp}^T Y_{BB}^{-1} \paren{i_{BB} + i_{Bp} + i_{Bq}} - \Sigma i_{Bp}, \qquad
y = Y_{Bp}^T Y_{BB}^{-1} Y_{Bq},
\label{eq:tnequiv}
\end{equation}
with $V_B^{oc}$ given as $I_B^{sc} / y$.
If a current source $I_{pq}$ or admittance $y_{pq}$ does appear across the port nodes $p$ and $q$ in $\nw$, its value should be added to $I_B^{sc}$ or $y$ respectively.
In the exceptional case $V_B^{oc} = 0$, one has $I_B^{sc} = 0$ also (and vice versa) and the equivalent `source' degenerates to an admittance $y$ between $p$ and $q$ which is always nonzero as a consequence of the connectivity of $\nw$, and corresponds to the driving-point admittance $Y_{pq} = 1 / Z_{pq}$ in $\nw_B$.
The degenerate case occurs in particular whenever $\nw_B$ contains no sources.
\end{proposition}
Note that with the same one-port decomposition, Proposition \ref{prop:tnequiv} can equally well be applied to define an equivalent source for subnetwork $\nw_A$ simply by interchanging the sets $A$ and $B$.

\subsection{Isolated voltage sources}
\label{sec:isovsrc}

Isolated current or voltage sources in a network (with no associated parallel or series admittances) are exceptional cases, not amenable to Th\'{e}venin--Norton conversion.
An isolated current source $I$, while it has no Th\'{e}venin equivalent, does not present any particular difficulty for any of the theory above, since it can be introduced to the fundamental network equations $Y v = i$ by means of additive entries $+I$ and $-I$ in the vector $i$.
Isolated voltage sources are a different matter: they can only be introduced automatically in this way after conversion to an equivalent Norton source, and an isolated source $V$ with no series admittance does not admit a Norton equivalent.

Of course, if the impedances of a network are fixed in advance one can always replace an isolated voltage source of value $V$ between nodes $p$ and $q$ with an isolated current source of value $I = Y_{pq} V$, where $Y_{pq}$ is the driving-point admittance.
However, what one generally desires is a way of formulating the network equations in the presence of the source $V$, in such a way that the formulation survives arbitrary changes in the balance of the network---in much the same way the Th\'{e}venin or Norton equivalent in Proposition \ref{prop:tnequiv} survives arbitrary changes in subnetwork $\nw_A$.

Fortunately, a formulation of equivalent current sources can be derived by closer consideration of the network equations $Y v = i$.
Again, suppose an isolated voltage source with value $V$ is connected between nodes $p$ and $q$.
Since the source fixes the voltage between these nodes, it may be assumed to `absorb' any current source or admittance element connected directly in parallel; thus, without loss of generality one may assume the source vector $i$ omits any current source between nodes $p$ and $q$, and the admittance matrix $Y$ omits any admittance between nodes $p$ and $q$.

Now, consider the network equations for the special case of a solution grounded at node $q$.
Purely for convenience of presentation, let the $n$ network nodes be temporarily relabelled so that $p = n - 1$ and $q = n$ are the last two in sequence, and set $\nu = n - 2$.
Then the (full) equations take the form
\begin{equation}
\mat{&&& -y_{p1} & -y_{q1} \\ & Y_{(pq,pq)} & & \vdots & \vdots \\ &&& -y_{p\nu} & -y_{q\nu} \\
      -y_{p1} & \cdots & -y_{p\nu} & Y_{pp} & 0 \\ -y_{q1} & \cdots & -y_{q\nu} & 0 & Y_{qq}}
   \mat{\vgnd{q}_1 \\ \vdots \\ \vgnd{q}_{\nu} \\ V \\ 0}
   = \mat{i_1 \\ \vdots \\ i_{\nu} \\ i_p \\ i_q},
\qquad
\begin{aligned} Y_{pp} &=\sum_{k=1}^{\nu} y_{pk}, \\ Y_{qq} &=\sum_{k=1}^{\nu} y_{qk}, \end{aligned}
\label{eq:vgndq}
\end{equation}
where $y_{pk}$ are the admittances (where present, else zero) connecting node $p$ to other nodes $k$, $y_{qk}$ are defined similarly for node $q$, $V$ is the source voltage, and the current injections on the right hand side represent the sources other than $V$.
(In particular, $i_p$ is nonzero if and only if there are one or more current sources directed to or from node $p$ itself.)
Let $m_p$ denote the number of admittance branches incident with node $p$; then $m_p$ is also the number of nonzero values $y_{pk}$.

By separating out column $p$ of the equations (\ref{eq:vgndq}), and omitting the redundant column $q$ entirely, one obtains the equivalent set of equations
\begin{equation}
Y_{(pq,pq)} \mat{\vgnd{q}_1 \\ \vdots \\ \vgnd{q}_{\nu}}
   = \mat{i_1 + y_{p1} V \\ \vdots \\ i_{\nu} + y_{p\nu} V}
\label{eq:vgndqi}
\end{equation}
corresponding to rows 1 through $\nu$, and the two equations
\begin{equation}
Y_{pp} V - \sum_{k=1}^{\nu} y_{pk} \vgnd{q}_k = i_p, \qquad
- \sum_{k=1}^{\nu} y_{qk} \vgnd{q}_k = i_q,
\label{eq:vgndqipq}
\end{equation}
corresponding to rows $p$ and $q$.
Now, equations (\ref{eq:vgndqi}) may be recognised as those for a grounded network similar to the original one, but where node $p$ has been deleted, and the voltage source $V$ replaced with $m_p$ current sources.
Each of these current sources corresponds to an admittance in the original network connecting node $p$ to some other node $k$; it is directed from the ground node $q$ to node $k$ and takes the value $y_{pk} V$.
Meanwhile, the original admittance element $y_{pk}$ persists in the truncated admittance matrix $Y_{(pq,pq)}$, but only in the diagonal element $Y_{kk}$ of that matrix.
It follows that equation (\ref{eq:vgndqi}) retains its validity with $y_{pk}$ now assumed to be connected between node $k$ and the ground node $q$ which is omitted from the grounded network equations.

Finally, observe that equations (\ref{eq:vgndqi}) can be `completed' to a \emph{new} set of full network equations $Y' v = i'$ by adjoining a \emph{single} row and column, such that the completed matrix $Y'$ has vanishing row and column sums, and vector $i'$ likewise a vanishing sum, giving the equations
\begin{equation}
\mat{&&& -(y_{p1} + y_{q1}) \\ & Y_{(pq,pq)} & & \vdots \\ &&& -(y_{p\nu} + y_{q\nu}) \\ -(y_{p1} + y_{q1}) & \cdots & -(y_{p\nu} + y_{q\nu}) & Y_{pp} + Y_{qq}}
   \mat{\vgnd{q}_1 \\ \vdots \\ \vgnd{q}_{\nu} \\ 0}
   = \mat{i_1 + y_{p1} V \\ \vdots \\ i_{\nu} + y_{p\nu} V \\ i_p + i_q - Y_{pp} V}.
\label{eq:vgndqq}
\end{equation}
By construction, (\ref{eq:vgndqq}) is a valid set of network equations: it may be verified by inspection that the required row and column sums are zero, that rows 1 through $\nu$ of these equations are exactly equivalent to (\ref{eq:vgndqi}), and that the last row is equivalent to (\ref{eq:vgndqipq}) by addition of the latter.
Notice that the effect of adding $i_p$ to the last equation here is to move all current source connections at node $p$ to node $q$, while the term $Y_{pp} V$ simply accounts for the combined `negative ends' of the current sources $y_{pk} V$ introduced in place of the voltage source $V$.
Since the node voltages in (\ref{eq:vgndqq}) are identical to those for the equivalent nodes in the original grounded network (that is, all nodes other than $p$), this will remain true for any general solutions of the original and new networks.

The following proposition sums up the equivalence thus obtained.
\begin{proposition}
\label{prop:isovsrc}
Let the network $\nw$ contain a single voltage source with value $V$, having its positive terminal at some node $p$ and its negative terminal at node $q$, with all other branches of $\nw$ being admittances and current sources.
Then at all nodes other than $p$, network $\nw$ behaves identically to the network $\nw'$ with admittances and current sources only, obtained from $\nw$ as follows:
\begin{itemise}
\item
voltage source $V$ is deleted and network $\nw$ contracted on nodes $p$ and $q$; and
\item
for each admittance element $y_{pk}$ connecting node $p$ to some other node $k \neq q$, a new current source is inserted, directed from node $q$ to node $k$ and with value $y_{pk} V$.
\end{itemise}
\end{proposition}
Observe that for the case $m_p = 1$ where only one admittance joins source $V$ at node $p$, the procedure defined by Proposition \ref{prop:isovsrc} is just the conventional Th\'{e}venin--Norton transformation.

The equivalence in Proposition \ref{prop:isovsrc}, one may add, follows also from more elementary considerations.
The node $p$, connected to $m_p$ separate admittance branches, can be notionally split into $m_p$ copies of itself, each connecting a copy of the source $V$ to one admittance in a Th\'{e}venin configuration.
With each node being at the same potential, the split nodes and the original node $p$ are equivalent.
Converting each to an equivalent Norton representation then results in the same equivalent configuration as given by the construction in Proposition \ref{prop:isovsrc}.
In the event a current source is connected to node $p$, one likewise adds a notional copy of the source between nodes $p$ and $q$, which has no effect on the network due to being in parallel with the voltage source $V$.
The net effect of the two copies of the current source is that of a single source connected to node $q$ in place of node $p$.

Proposition \ref{prop:isovsrc} has the following immediate corollary:
\begin{proposition}
\label{prop:isovi}
Let a network $\nw$ contain a voltage source between two nodes $p$ and $q$ where the only other branches connecting to node $p$ (or alternatively node $q$) are current source branches.
Then the voltage source may be deleted and the network $\nw$ contracted on nodes $p$ and $q$ with no effect on the network behaviour.
\end{proposition}
As a special case of Proposition \ref{prop:isovi} one obtains the elementary fact that a voltage source in series with a current source is equivalent to the current source alone.

\subsection{Linear dependent sources}
\label{sec:depsource}

In passive networks as considered to this point, all sources are independent---their current or voltage is externally fixed.
But the admittance matrix framework can also be extended to consider (linear) sources whose current or voltage depends on a quantity elsewhere in the network.
Such `dependent' or `controlled' sources are capable of modelling more general network phenomena such as mutural coupling, ideal transformers, amplifiers or hybrid DC--AC networks.
In general, however, networks with dependent sources cease to be passive, and the general properties previously described for passive linear networks, such as reciprocity or monotonicity of driving-point impedances, no longer apply.
This is in essence because the equivalent admittance matrix $Y$ for such networks fails to have the symmetric structure with zero row and column sums assumed at the outset.

The dependent source most readily incorporated into the $Y v = i$ framework is the linear voltage-dependent current source: an element that, when placed between any two nodes, produces a current between its terminals
\begin{equation}
i = \mathcal{Y} \paren{v_p - v_q}
\label{eq:vdcs}
\end{equation}
proportional to the voltage difference between two designated nodes $p$ and $q$.
The effect of such a source, with its positive terminal at node $j$ and its negative terminal at node $k$, is easily modelled by adjusting elements of the admittance matrix $Y$ as follows:
\begin{itemise}
\item
elements $Y_{jp}$ and $Y_{kq}$ are adjusted by an amount $+\mathcal{Y}$; and
\item
elements $Y_{jq}$ and $Y_{kp}$ are adjusted by an amount $-\mathcal{Y}$.
\end{itemise}
One immediately sees that this breaks the symmetry of the matrix $Y$, except in the uninteresting cases $j = p$, $k = q$ (where it becomes an ordinary admittance element with value $\mathcal{Y}$) and $j = q$, $k = p$ (where it is an ordinary admittance with value $-\mathcal{Y}$).
In the general case, sources of the form (\ref{eq:vdcs}) can readily be used to generate counterexamples to other general properties of the $Y$ matrix such as diagonal dominance and positive semidefiniteness, and to corresponding properties of the underlying electrical network such as reciprocity and passivity.
Note however that the zero-sum property of rows and columns of $Y$ is unaffected, meaning that uniqueness of solutions (up to uniform adjustment of voltages) and the role of grounded solutions is preserved.

A source may be specified as current-dependent rather than voltage-dependent, but in the framework developed here the two are almost equivalent, since a general network branch is either an admittance element (with current proportional to voltage difference) or an independent source (with fixed current).
Consider a current-dependent current source $i = \mathcal{K} i_{\alpha}$ with transfer gain $\mathcal{K}$ and $i_{\alpha}$ the current in an admittance element $y$ directed from node $q$ to node $p$; then the rule $i = \mathcal{K} i_{\alpha}$ can be rendered into the form
\begin{equation}
i = \mathcal{K} y \paren{v_p - v_q}
\label{eq:cdcs}
\end{equation}
matching that for a voltage-dependent current source.
As above, the source is incorporated by adjusting elements of the admittance matrix by amounts $\pm \mathcal{K} y$, with the caveat that the adjustment has a dependency on the value of the admittance $y$ in the original network.

Dependent voltage sources can also be incorporated into the framework by transforming them to equvialent dependent current sources, using the theory leading up to Propositions \ref{prop:thevnort} and \ref{prop:isovsrc}.
The simplest case is a voltage-dependent voltage source
\begin{equation}
v = \mathcal{K} \paren{v_p - v_q}
\label{eq:vdvs}
\end{equation}
with transfer gain $\mathcal{K}$, that appears in series with an admittance element $y$.
By an algebraically valid extension of the Th\'{e}venin--Norton correspondence shown earlier, this is exactly equivalent to a dependent \emph{current} source with the same rule (\ref{eq:cdcs}) as above, connected in parallel with the admittance $y$ (and with the original voltage source removed by contraction).

Of course, when replacing (\ref{eq:vdvs}) with (\ref{eq:cdcs}) for a dependent voltage source, the value $y$ stands for something quite different than in the original current source: there it is equal to the admittance in the branch bearing the dependent current, while for the source (\ref{eq:vdvs}) $y$ is an assumed Th\'{e}venin admittance associated with the source element itself.
\emph{Both} types of associated admittance come into play when incorporating a current-dependent voltage source into the framework: this is a source producing voltage according to the rule
\begin{equation}
v = \mathcal{Z} i_{\alpha}
\label{eq:cdvs}
\end{equation}
where $i_{\alpha}$ is the current in a specified branch $\alpha$ directed from node $q$ to node $p$ (say) and with admittance $y_{\alpha}$.
Assuming this voltage source is itself connected in series with a Th\'{e}venin admittance $y_{\tau}$, this can likewise be rendered into a voltage-dependent current source
\begin{equation}
i = y_{\tau} \mathcal{Z} y_{\alpha} \paren{v_p - v_q}
\label{eq:cdvsequiv}
\end{equation}
where both the Th\'{e}venin admittance $y_{\tau}$ and the dependent branch admittance $y_{\alpha}$ act to scale the source's transfer gain $\mathcal{Z}$.
The adjustments to the admittance matrix are by amounts $\pm \mathcal{Z} y_{\tau} y_{\alpha}$.

It may come about that the dependent voltage source (\ref{eq:vdvs}) or (\ref{eq:cdvs}) is not in Th\'{e}venin configuration, in which case a generalisation of Proposition \ref{prop:isovsrc} may be developed.
Once again, this expresses the source as equivalent to a collection of dependent current sources with a common terminal, each expressible with a rule of the form (\ref{eq:vdcs}) or (\ref{eq:cdcs}) indicating how the source may be incorporated into the admittance matrix.

Overall then, the observations above regarding voltage-dependent current sources generalise to all linear dependent sources: the equivalent adjustments to the admittance matrix $Y$ that result will in general break the symmetry that applies with independent sources, and so properties such as reciprocity and monotonicity may not hold.
One may however still rely on linearity of the equations and solutions, and on the uniqueness and grounding properties that follow from $Y$ having vanishing row and column sums.

\bibliographystyle{unsrt}
\bibliography{systems}

\begin{thebibliography}{1}

\bibitem{j:mtoeam}
James~H. Jeans.
\newblock {\em The Mathematical Theory of Electricity and Magnetism}.
\newblock Cambridge, 1911.

\bibitem{bsst:doris}
R.~L. Brooks, C.~A.~B. Smith, A.~H. Stone, and W.~T. Tutte.
\newblock The dissection of rectangles into squares.
\newblock {\em Duke Mathematical Journal}, 7(1):312--340, 1940.

\bibitem{f:aioen}
R.~M. Foster.
\newblock The average impedance of an electrical network.
\newblock In {\em Reissner Anniversary Volume, Contributions in Applied
  Mechanics}, pages 333--340. J. W. Edwards, Ann Arbor, MI, 1948.

\bibitem{hj:ma}
Roger~A. Horn and Charles~R. Johnson.
\newblock {\em Matrix Analysis}.
\newblock Cambridge, 1985.

\bibitem{bs:di}
R.A. Brualdi and H.~Schneider.
\newblock Determinantal identities: {G}auss, {S}chur, {C}auchy, {S}ylvester,
  {K}ronecker, {J}acobi, {B}inet, {L}aplace, {M}uir, and {C}ayley.
\newblock {\em Linear Algebra and its Applications}, 52/53:769--791, 1983.

\bibitem{kr:rd}
D.J. Klein and M.~Randi\'{c}.
\newblock Resistance distance.
\newblock {\em Journal of Mathematical Chemistry}, 12:81--95, 1993.

\bibitem{k:rdsr}
Douglas~J. Klein.
\newblock Resistance-distance sum rules.
\newblock {\em Croatica Chemica Acta}, 75(2):633--649, 2002.

\end{thebibliography}

\end{document}